\documentclass[12pt]{amsart}
\usepackage{amsmath,amsthm,amsfonts,amssymb}
\usepackage{graphicx}
\usepackage{mathrsfs}
\usepackage{enumitem}
\usepackage{etoolbox}
\usepackage{xcolor}
\apptocmd{\sloppy}{\hbadness 10000\relax}{}{}
\apptocmd{\sloppy}{\vbadness 10000\relax}{}{}
\usepackage[pdfpagelabels]{hyperref}
\usepackage[letterpaper,margin=1.1in]{geometry}

\numberwithin{equation}{section}
\theoremstyle{plain}
\newtheorem{theorem}{Theorem}[section]
\newtheorem{proposition}[theorem]{Proposition}
\newtheorem{corollary}[theorem]{Corollary}

\newtheorem{lemma}[theorem]{Lemma}

\theoremstyle{definition}
\newtheorem{remark}[theorem]{Remark}
\newtheorem{definition}[theorem]{Definition}
\newtheorem{example}[theorem]{Example}
\newtheorem{algorithm}[theorem]{Algorithm}

\newcommand{\Span}{\mathop\mathrm{span}\nolimits}

\def\RR{\mathbb{R}}
\def\ZZ{\mathbb{Z}}
\def\N{\mathbb{N}}

\def\XX{{\mathbb{X}}}
\def\YY{{\mathbb{Y}}}

\newcommand{\Flat}{\mathsf{Flat}}

\newcommand{\diam}{\mathop\mathrm{diam}\nolimits}
\newcommand{\dist}{\mathop\mathrm{dist}\nolimits}

\newcommand{\Lip}{\mathop\mathrm{Lip}\nolimits}

\newcommand{\Haus}{\mathcal{H}}
\newcommand{\Start}{\mathop\mathsf{Start}\nolimits}
\newcommand{\End}{\mathop\mathsf{End}\nolimits}
\newcommand{\Domain}{\mathop\mathsf{Domain}\nolimits}
\newcommand{\Image}{\mathop\mathsf{Image}\nolimits}

\newcommand{\Diam}{\mathop\mathsf{Diam}\nolimits}

\newcommand{\Edge}{\mathop\mathsf{Edge}\nolimits}

\newcommand{\var}{\mathop\mathrm{var}\nolimits}
\newcommand{\flatepsilon}{\epsilon_2}
\newcommand{\flatarcsepsilon}{\epsilon_1}

\newcommand{\excess}{\mathop\mathrm{excess}\nolimits}
\newcommand{\gap}{\mathop\mathrm{gap}\nolimits}

\begin{document}

\title[Subsets of rectifiable curves in Banach spaces]{Subsets of rectifiable curves in Banach spaces I:\\ sharp exponents in traveling salesman theorems}

\author{Matthew Badger \and Sean McCurdy}
\thanks{M.~Badger was partially supported by NSF DMS grants 1650546 and 2154047.}
\date{August 23, 2022}

\address{Department of Mathematics\\ University of Connecticut\\ Storrs, CT 06269-1009}
\email{matthew.badger@uconn.edu}
\address{Department of Mathematics\\ National Taiwan Normal University\\ Taipei City\\ Taiwan (R.\,O.\,C.)}
\email{smccurdy@ntnu.edu.tw}

\begin{abstract} The Analyst's Traveling Salesman Problem is to find a characterization of subsets of rectifiable curves in a metric space. This problem was introduced and solved in the plane by Jones in 1990 and subsequently solved in higher-dimensional Euclidean spaces by Okikiolu in 1992 and in the infinite-dimensional Hilbert space $\ell_2$ by Schul in 2007. In this paper, we establish sharp extensions of Schul's necessary and sufficient conditions for a bounded set $E\subset\ell_p$ to be contained in a rectifiable curve from $p=2$ to $1<p<\infty$. While the necessary and sufficient conditions coincide when $p=2$, we demonstrate that there is a strict gap between the necessary condition and sufficient condition when $p\neq 2$. We also identify and correct technical errors in the proof by Schul. This investigation is partly motivated by recent work of Edelen, Naber, and Valtorta on Reifenberg-type theorems in Banach spaces and complements work of Hahlomaa and recent work of David and Schul on the Analyst's TSP in general metric spaces.\end{abstract}

\maketitle

\tableofcontents

\renewcommand{\thepart}{\Roman{part}}

\section{Introduction}\label{sec:intro}

Given a set in a path-connected metric space, we may ask whether or not the set is contained in a curve of finite length (also called a \emph{rectifiable curve}), and if so, ask how to find a curve containing the set that is (essentially) as short as possible. This problem was introduced and solved in the Euclidean plane by Jones \cite{Jones-TST} and is now commonly known as the \emph{Analyst's Traveling Salesman Problem}. While it is immediate that a set contained in a rectifiable curve is necessarily bounded and has finite one-dimensional Hausdorff measure $\Haus^1$, this pair of conditions is not sufficient. To decide when a set is contained in a rectifiable curve requires additional information about the local geometry of the set within the space (see below). Full solutions to the Analyst's Traveling Salesman Problem are currently available in $\RR^n$ \cite{Ok-TST}, Carnot groups \cite{Li-TSP}, infinite-dimensional Hilbert space $\ell_2$ \cite{Schul-Hilbert}, graph inverse limit spaces \cite{Guy-Schul}, and for Radon measures in $\RR^n$ \cite{BS3} and Carnot groups \cite{BLZ}; partial solutions are available in general metrics spaces \cite{Hah05, DS-metric} as well as for higher-dimensional curves \cite{BNV,Balogh-Zust} and surfaces \cite{DT, AS-TST,ENV-Banach,Ghinassi,AV-q,Villa-TST,Hyde-TST}. Beyond the intrinsic interest of the Analyst's TSP in metric geometry, finding tests to determine when a set is contained in a rectifiable curve or ``nice" surface has led to applications in complex analysis, dynamics and probability, geometric analysis, and harmonic analysis. For a sample of applications, see \cite{BJ, wiggly, frontier, Tolsa-P, shortcuts, logic-tst, AT15, NV-annals,Badger-survey,BV,Azzam-hmtst,BV-IFS,GM-wild,Naples20,Badger-Naples}.

In this paper, we establish sharp extensions of Schul's necessary and sufficient conditions for a bounded set $E\subset\ell_p$ to be contained in a rectifiable curve from $p=2$ to $1<p<\infty$ (see Theorems \ref{t:main} and \ref{t:main2}). As usual, $\ell_p$ denotes the (real) Banach space of $p$-summable sequences, $$x=(x_1,x_2,\dots)\in\ell_p\quad\text{if and only if}\quad x_1,x_2,\dots\in\RR\text{ and }|x|_p=\left(\sum_1^\infty |x_i|^p\right)^{1/p}<\infty.$$ While the necessary and sufficient conditions in Schul's theorem coincide when $p=2$, we demonstrate that there is a strict gap between the necessary condition and sufficient condition when $p\neq 2$. En route, we prove that the classes of rectifiable curves in the infinite-dimensional spaces $\ell_p$ and $\ell_q$ differ when $p\neq q$:

\begin{proposition}\label{t:pq} Let $1<p<\infty$. Every rectifiable curve in $\ell_p$ is a rectifiable curve in $\ell_q$ for all $q\geq p$. However, there exists a curve $\Gamma$ in $\ell_p$ such that $\Gamma$ is rectifiable (i.e.~$\Haus^1(\Gamma)<\infty$) in $\ell_q$ for every $q>p$, but $\Gamma$ is not rectifiable (i.e.~$\Haus^1(\Gamma)=\infty$) in $\ell_p$.\end{proposition}

Proposition \ref{t:pq}, Theorem \ref{t:main}, and Theorem \ref{t:main2} capture a special infinite-dimensional phenomenon. In particular, they imply that a solution of the Analyst's TSP in $\ell_p$ cannot be neatly derived from the solution in $\ell_2$. By contrast, bi-Lipschitz equivalence of finite-dimensional Banach spaces ensures that a set in $\RR^n$ is (a subset of) a rectifiable curve independent from the choice of underlying norm, even though the actual length of the curve depends on the norm. This paper serves to clarify the difference between the finite and infinite-dimensional settings. An essential reason for us to study (subsets of) rectifiable curves in $\ell_p$ for $1<p<\infty$ is that the spaces interpolate between $\ell_2$, where the Analyst's TSP is solved, and $\ell_\infty$, which contains an isometric copy of any separable metric space. Thus, a resolution of the Analyst's TSP in $\ell_p$ may provide insight into the Analyst's TSP in general metric spaces. For further discussion and description of related research, see \S\ref{sec:related}.

\subsection{Analyst's TSP in Euclidean space and \texorpdfstring{$\ell_2$}{l\_2}}
\label{sec:euclidean}

To solve the Analyst's TSP in the plane, Jones introduced unilateral linear approximation numbers, now universally called \emph{Jones' beta numbers}, which measure how ``flat" a set is in a given window. The Jones' beta numbers are well-defined in any Banach space. Let $\XX$ be a Banach space\footnote{All Banach spaces in this paper are real Banach spaces of dimension at least 2.}, let $E\subset \XX$ be a nonempty set, and let $Q\subset\XX$ be a set of positive diameter. If $E\cap Q\neq \emptyset$, define \begin{equation}\label{def:beta} \beta_E(Q) = \inf_{L} \sup_{x\in E\cap Q}\frac{\dist(x,L)}{\diam Q}\in[0,1],\end{equation} where the infimum ranges over all one-dimensional affine subspaces (lines) $L\subset \XX$; and, if $E\cap Q=\emptyset$, define $\beta_E(Q)=0$. At one extreme, if $\beta_E(Q)=0$, then $E\cap Q$ is contained in a line. At the other extreme,  if $\beta_E(Q) \gtrsim 1$, then the set $E\cap Q$ is uniformly far away from every line passing through $Q$. From the definition, it immediately follows that \begin{equation} \label{beta-monotone} \beta_E(R) \leq \frac{\diam Q}{\diam R}\, \beta_F(Q)\quad\text{for all }E\subset F\text{ and }R\subset Q.\end{equation}  In view of the fact that rectifiable curves (having parameterizations of bounded variation) admit tangent lines $\Haus^1$-a.e., one may expect that sets contained in a finite length curve have ``vanishing beta numbers" at typical points of those sets. The \emph{Analyst's Traveling Salesman Theorem} makes this idea precise and provides a characterization of subsets of rectifiable curves with an estimate on the shortest length of a curve containing the set:

\begin{theorem}[Jones \cite{Jones-TST} in $\RR^2$; Okikiolu \cite{Ok-TST} in $\RR^n$] \label{t:Jones} Let $n\geq 2$ and let $E\subset\RR^n$. Then $E$ is contained in a rectifiable curve if and only if \begin{equation}\label{e:S-finite} S_E(\RR^n):= \diam E + \sum_{Q\in\Delta(\RR^n)} \beta_E(3Q)^2 \diam Q<\infty,\end{equation} where the sum ranges over all dyadic cubes $Q$ in $\RR^n$ and $3Q$ denotes the concentric dilate of the cube $Q$ with scaling factor 3. More precisely, if $S_E(\RR^n)<\infty$, then $E$ is contained in a curve $\Gamma$ in $\RR^n$ with \begin{equation}\label{e:Jones} \Haus^1(\Gamma) \lesssim_n S_E(\RR^n).\end{equation} If $\Sigma\subset \RR^n$ is a connected set, then \begin{equation} \label{e:Okikiolu} S_{\Sigma}(\RR^n) \lesssim_n \Haus^1(\Sigma).\end{equation} The constant $3$ in \eqref{e:S-finite} can be replaced with any constant $A>1$. Then \eqref{e:Jones} and \eqref{e:Okikiolu} hold with implicit constants depending on $n$ and $A$.\end{theorem}

For refinements of \eqref{e:Jones}, \eqref{e:Okikiolu} for rectifiable Jordan arcs in $\RR^n$ and in Hilbert space, see \cite{Bishop-TST},  \cite{Krandel-Hilbert}. Time complexity of the Analyst's Traveling Salesman algorithm for constructing the rectifiable curve in \eqref{e:Jones} is investigated in \cite{TSP-complexity}.

\begin{remark} If $E$ is a subset of a rectifiable curve in $\RR^n$, then the Analyst's Traveling Salesman Theorem ensures  \begin{equation}\sum_{\stackrel{Q\in\Delta(\RR^n)}{\beta_E(3Q)\geq \epsilon}}\diam Q \leq \epsilon^{-2} S(E)<\infty\quad\text{for all }\epsilon>0.\end{equation} It follows that \begin{equation}\lim_{\stackrel{Q\in\Delta(\RR^n)}{Q\downarrow x}} \beta_E(3Q)=0\quad\text{at $\Haus^1$-a.e. }x\in E,\end{equation} or equivalently by \eqref{beta-monotone}, \begin{equation}\label{e:beta-vanish} \lim_{r\downarrow 0} \beta_E(B(x,r))=0\quad\text{at $\Haus^1$-a.e. }x\in E.\end{equation} Thus, subsets of rectifiable curves in $\RR^n$ have ``vanishing beta numbers" at typical points in the sense of \eqref{e:beta-vanish}. It is possible to construct examples of generalized von Koch snowflake curves (with carefully chosen angles) to show that \eqref{e:beta-vanish} is also satisfied by certain curves of infinite length (Proposition \ref{prop:2}). By contrast, the Analyst's Traveling Salesman Theorem guarantees that every curve $\Gamma$ in $\RR^n$ of infinite length satisfies $S_\Gamma(\RR^n)=\infty$. In other words, finiteness of $S_\Gamma(\RR^n)$ is a perfect test to determine rectifiability of a curve $\Gamma$ in $\RR^n$.\end{remark}

\begin{figure}\begin{center}\includegraphics[width=.8\textwidth]{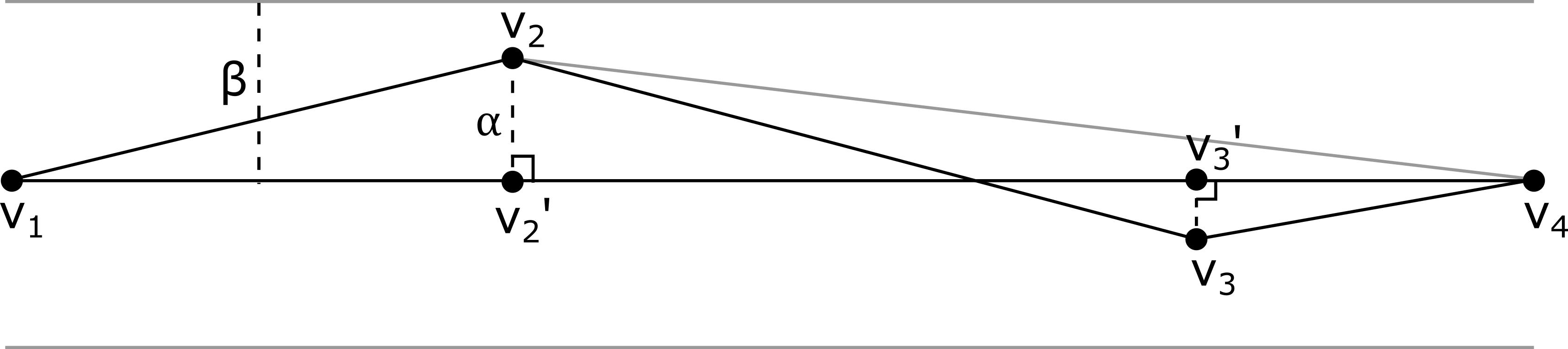}\end{center}\caption{To use the Pythagorean theorem to verify \eqref{e:ex-above} and \eqref{e:ex-below} (shown with $k=4$), first draw right triangles formed from line segments $[v_i,v_{i+1}]$ between consecutive points and the line passing through $v_1$ and $v_k$.}\label{fig:3}\end{figure}

\begin{remark}\label{r:pythag} Let $V$ be a finite set of points in $\RR^n$ (equipped with the Euclidean norm) that is \emph{1-separated} in the sense $|v-w|\geq 1$ for all distinct $v,w\in V$. Assume that $L$ is a line in $\RR^n$ such that $\dist(v,L)\leq \beta\ll 1$ for all $v\in V$. Then the set $V=\{v_1,\dots,v_k\}$ may be enumerated according to its orthogonal projection onto $L$. For simplicity, let us further assume that $v_1,v_k\in \ell$. By the triangle inequality and a simple computation with the Pythagorean theorem (see Figure \ref{fig:3}), \begin{equation}\label{e:ex-above} |v_1-v_k| \leq |v_1-v_2|+\dots+|v_{k-1}-v_k| \leq (1+C_1\beta^2)|v_1-v_k|\end{equation} for some universal constant $C_1$. Conversely, if $\dist(v_i,L)\geq \alpha$ for some $1\leq i\leq k$, then the Pythagorean theorem yields \begin{equation}\label{e:ex-below} |v_1-v_i|+|v_i-v_k| \geq  (1+C_2\alpha^2) |v_1-v_k|\end{equation} for some universal constant $C_2$. At a high level, the estimates \eqref{e:ex-above} and \eqref{e:ex-below} correspond to \eqref{e:Jones} and \eqref{e:Okikiolu} in Analyst's Traveling Salesman Theorem, respectively. Informally, we say that the Pythagorean theorem is responsible for the exponent 2 on $\beta_E(3Q)^2$ in \eqref{e:S-finite}.\end{remark}

The dependence on the ambient dimension  in the implicit constants in \eqref{e:Jones} and \eqref{e:Okikiolu} is ultimately a consequence of using dyadic cubes in the sum $S_E(\RR^n)$ to index ``all" locations and scales in $\RR^n$. By bi-Lipschitz equivalence, Theorem \ref{t:Jones} persists in every finite-dimensional Banach space with implicit constants that also depend on the norm as well as a choice of coordinates; see \S\ref{sec:finite-dimensional} for a detailed statement. To formulate a dimension independent version of the Analyst's Traveling Salesman Theorem in $\ell_2$, Schul \cite{Schul-Hilbert} replaced $S_E(\RR^n)$ with a sum $S_E(\mathscr{G})$ indexed over a multiresolution family $\mathscr{G}$ of ``all" locations and scales in the set $E$.

Let $\mathbb{X}$ be a metric space and let $E\subset \mathbb{X}$ be a nonempty set. For any $\rho>0$, a \emph{$\rho$-net} $X_\rho$ in $\mathbb{X}$ is a set such that $\dist(y,z)\geq \rho$ for all distinct $y,z\in X_\rho$ and $\dist(x,X_\rho)<\rho$ for all $x\in \mathbb{X}$. Following \cite{Schul-Hilbert}, we define a \emph{multiresolution family} $\mathscr{G}$ for $E$ with \emph{inflation factor} $A_{\mathscr{G}}>1$ to be a collection of closed balls of the form \begin{equation}\mathscr{G}=\{B(x,A_{\mathscr{G}} 2^{-k}): x\in X_k, k\in\ZZ\},\end{equation} where $(X_k)_{k\in\ZZ}$ is a nested family of $2^{-k}$-nets for $E$. Analogously, if each set $X_k$ has the property $\dist(y,z)\geq 2^{-k}$ for all distinct $y,z\in X_k$, but one or more of the sets $X_k$ are not $2^{-k}$-nets, then we call $\mathscr{G}$ a \emph{partial multiresolution family} for $E$.

For any nonempty set $E\subset \ell_2$ and (partial) multiresolution family $\mathscr{G}$ for $E$, define \begin{equation} \label{e:Jones-sum2}S_{E}(\mathscr{G}):=\diam E+\sum_{Q\in\mathscr{G}} \beta_E(Q)^2 \diam Q.\end{equation}

\begin{theorem}[Schul \cite{Schul-Hilbert}] \label{t:schul} If $E\subset\ell_2$ and $S_{E}(\mathscr{G})<\infty$ for some multiresolution family $\mathscr{G}$ for $E$ with inflation factor $A_\mathscr{G}> 200$, then $E$ is contained in a rectifiable curve $\Gamma$ in $\ell_2$ with \begin{equation}\label{e:suff-2ss}\Haus^1(\Gamma) \lesssim_{A_\mathscr{G}} S_{E}(\mathscr{G}).\end{equation} If $\Sigma\subset\ell_2$ is a connected set and $\mathscr{H}$ is a (partial) multiresolution family for $\Sigma$ with inflation factor $A_\mathscr{H} > 1$, then \begin{equation}\label{e:nec-2ss} S_{\Sigma}(\mathscr{H})\lesssim_{A_{\mathscr{H}}} \Haus^1(\Sigma).\end{equation}\end{theorem}

Note that the implicit constants in \eqref{e:suff-2ss} and \eqref{e:nec-2ss} depend on the inflation factor of the multiresolution family, but are otherwise dimension free. Once again, the Pythagorean theorem in $\ell_2$ determines the exponent 2 on $\beta_E(Q)^2$ in Theorem \ref{t:schul} \emph{\`a la} Remark \ref{r:pythag}. More generally, Theorem \ref{t:schul} holds in \emph{any} Hilbert space (including nonseparable spaces).

\subsection{Schul's theorem in \texorpdfstring{$\ell_p$}{l\_p}} \label{sec:new}

For any nonempty set $E\subset \ell_p$ and (partial) multiresolution family $\mathscr{G}$ for $E$, define \begin{equation} \label{e:Jones-sum}S_{E,r}(\mathscr{G}):=\diam E+\sum_{Q\in\mathscr{G}} \beta_E(Q)^r \diam Q\quad\text{for all }0<r<\infty.\end{equation} Note that in the notation of the previous section, $S_{E,2}(\mathscr{G})\equiv S_E(\mathscr{G})$. Moreover, by \eqref{beta-monotone}, \begin{equation}\label{S-monotone} S_{E,r}(\mathscr{G}) \leq S_{F,s}(\widehat{\mathscr{G}})\quad\text{for all }E\subset F \text{ and } s\leq r,\end{equation} where $\widehat{\mathscr{G}}$ is a multiresolution family for $F$ extending $\mathscr{G}$.

The following pair of theorems, extending Schul's theorem from $\ell_2$ to $\ell_p$ with $1<p<\infty$, constitute our main result. We emphasize that when $p\neq 2$, there is a strict gap between the necessary and sufficient conditions for a set to be contained in a rectifiable curve.

\begin{theorem}[sharp sufficient conditions in {$\ell_p$}] \label{t:main} Let $1<p<\infty$. If $E\subset\ell_p$ and $S_{E,\min(p,2)}(\mathscr{G})<\infty$ for some multiresolution family $\mathscr{G}$ for $E$ with inflation factor $A_\mathscr{G}\geq 240$, then $E$ is contained in a curve $\Gamma$ in $\ell_p$ with \begin{equation}\label{e:suff-p}\Haus^1(\Gamma) \lesssim_{p,A_\mathscr{G}} S_{E,\min(p,2)}(\mathscr{G}).\end{equation} The exponent $\min(p,2)$ on beta numbers in \eqref{e:suff-p} is sharp.
\end{theorem}

\begin{theorem}[sharp necessary conditions in {$\ell_p$}] \label{t:main2} Let $1<p<\infty$. If $\Sigma\subset\ell_p$ is a connected set and $\mathscr{H}$ is a (partial) multiresolution family for $\Sigma$ with inflation factor $A_\mathscr{H} > 1$, then \begin{equation}\label{e:nec-2} S_{\Sigma,\max(2,p)}(\mathscr{H})\lesssim_{p,A_\mathscr{H}} \Haus^1(\Sigma).\end{equation} The exponent $\max(2,p)$ on beta numbers in \eqref{e:nec-2} is sharp.
\end{theorem}

\begin{remark} In this paper, we prove Theorem \ref{t:main} in full and we reduce the proof of Theorem \ref{t:main2} to Theorem \ref{l:H2}; the proof of the latter theorem is deferred to Part II \cite{Badger-McCurdy-2}. More generally, we show that an analogue of Theorem \ref{t:main} holds in uniformly smooth Banach spaces (see Theorem \ref{t:smooth}) and an analogue of Theorem \ref{t:main2} holds in uniformly convex spaces (see Theorem \ref{t:convex}). For example, $L^p$ spaces of $p$-power integrable functions over a measure space, which include $\ell_p$, are uniformly smooth and uniformly convex for all $1<p<\infty$. In addition, we prove that a universal sufficient condition with exponent 1 is valid in arbitrary Banach spaces (see Theorem \ref{t:banach-suff}). Each of these results apply to both separable and nonseparable Banach spaces.\end{remark}

\begin{figure}\includegraphics[height=.28\textheight]{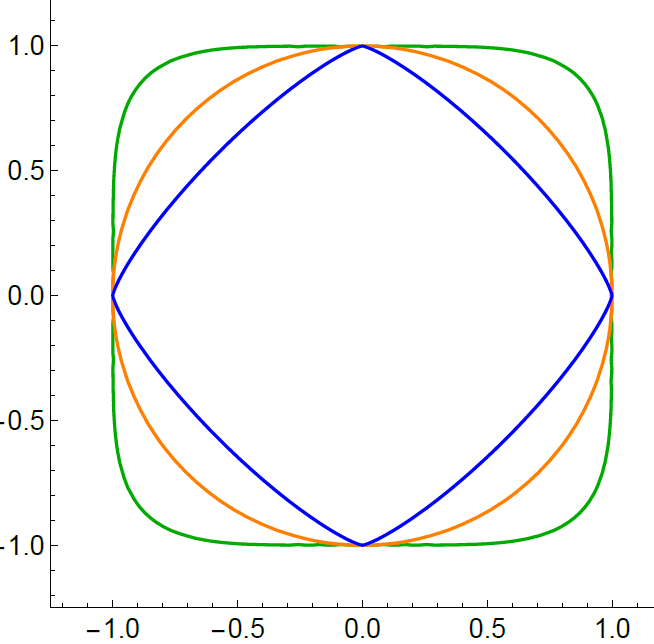}\hspace{.1\textwidth} \includegraphics[height=.3\textheight]{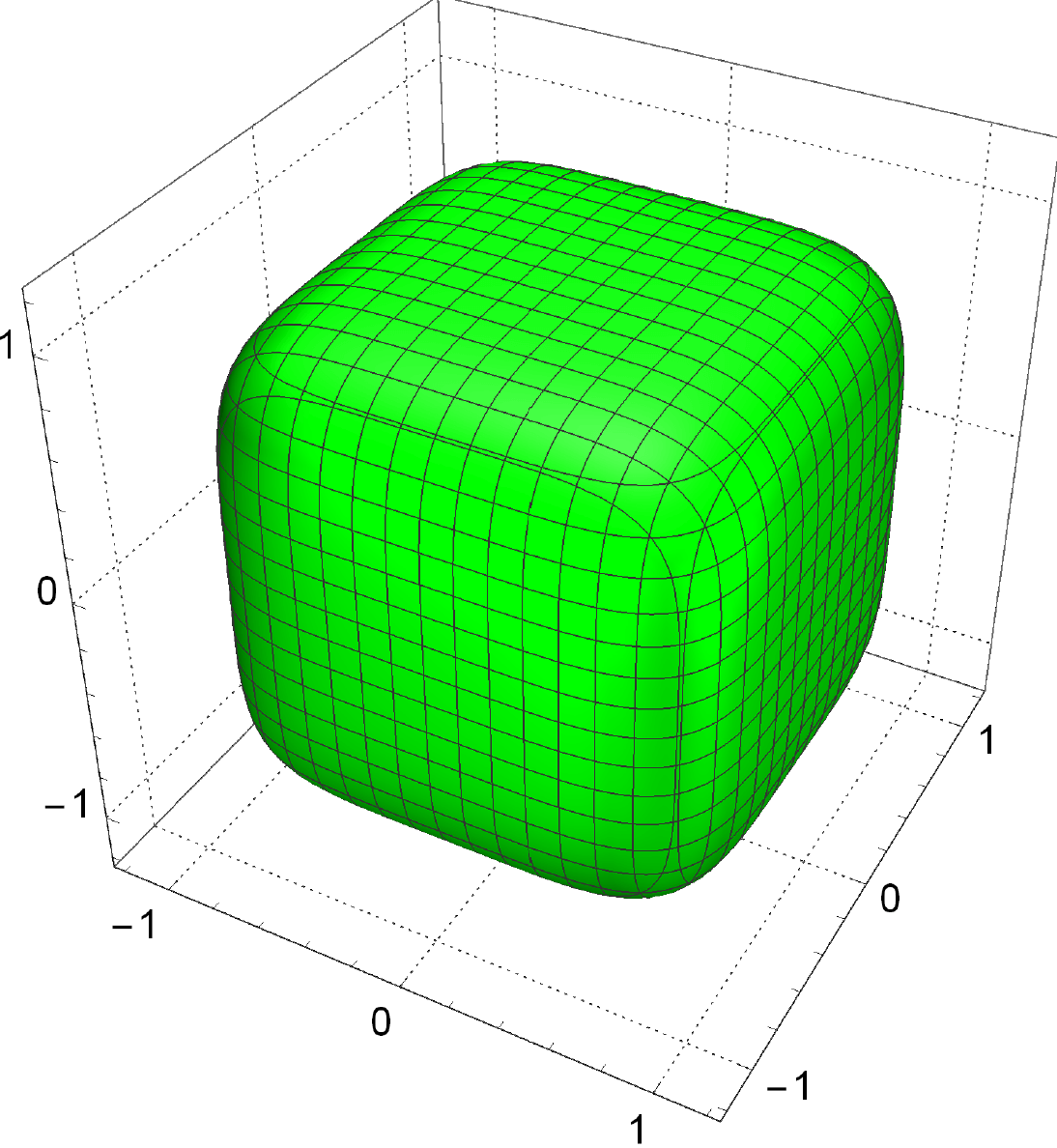} \caption{Unit balls in $\ell_{5/4}^2$ (blue), $\RR^2$ (orange), $\ell_{5}^2$ and $\ell_5^3$ (green), where $\ell^n_p$ denotes $(\RR^n,|\cdot|_p)$.}\label{fig:4}\end{figure}

An essential feature of $\ell_2$ is that the unit ball and induced distance are \emph{rotationally-invariant}. In particular, to compute the distance of a point $x$ to a line $L$ in Hilbert space, one may first translate and rotate so that $L=\Span(e_1)$ and $x\in\Span(e_1,e_2)$ if convenient. In $\ell_p$, when $p\neq 2$, rotational-invariance is no longer available and computation of the distance of a point to a line is sensitive to the position of the line and geometry of the unit ball. We remark that the gain in complexity witnessed when moving from $\ell_p^2$ to $\ell_p^3$ (e.g.~ consider the shape of slices of their unit balls, see Figure \ref{fig:4}) continues to increase when moving from the finite-dimensional spaces $\ell^n_p$ to the infinite-dimensional space $\ell_p$. For instance, although the norms in $\ell_2^n$ and $\ell_p^n$ are $C(n,p)$-bi-Lipschitz equivalent for each pair $n$ and $p$, the bi-Lipschitz constant $C(n,p)$ degenerates as $n\rightarrow\infty$ for each $p\neq 2$.

\begin{example} \label{ex:l5} To illustrate the essential idea behind the exponents in Theorems \ref{t:main} and \ref{t:main2}, let's estimate the length gain in $\ell_5^2=(\RR^2,|\cdot|_5)$ of isosceles triangles \begin{equation} \begin{split} T_h &: a_h=(0,0),\ b_h=(l/2,h),\ c_h=(l,0),\text{ and} \\  T_d &: a_d=(0,0),\ b_d=(2^{-6/5}l-2^{-1/5}h,2^{-6/5}l+2^{-1/5}h),\ c_d=(2^{-1/5}l,2^{-1/5}l)\end{split}\end{equation} with horizontal and diagonal bases $\overline{ac}$ of length $l$ and height $\dist(b,\overline{ac})=h$. On one hand, \begin{equation}\begin{split} |\overline{a_hb_h}|_5+|\overline{b_hc_h}|_5-|\overline{a_hc_h}|_5 &= 2\left((l/2)^5+h^5\right)^{1/5}-l \\ &= l\left((1+32(h/l)^5)^{1/5}-1\right)\simeq l(h/l)^5\end{split}\end{equation} for $h\ll l$ by Taylor's theorem for $x\mapsto x^{1/5}$ at $x=1$. On the other hand,  \begin{equation}\begin{split}|\overline{a_db_d}|_5+|\overline{b_dc_d}|_5-|\overline{a_dc_d}|_5 &= 2\left((2^{-6/5}l-2^{-1/5}h)^5+(2^{-6/5}l+2^{-1/5}h)^5\right)^{1/5}-l\\ &= l\left((1+40(h/l)^2+80(h/l)^4)^{1/5}-1\right)\simeq l(h/l)^2\end{split}\end{equation} when $h\ll l$ by Taylor's theorem, again. This computation indicates that on the axial directions the excess in the triangle inequality for ``flat" triangles is comparable to $l(h/l)^5$, while on the diagonal directions, the excess is comparable to $l(h/l)^2$. In fact, as long as the base of the triangle is neither horizontal nor vertical, the excess is comparable to $l(h/l)^2$ provided that $h\ll l$ (depending on proximity of the base to an axial direction). In higher dimensions, the excess in the triangle inequality also depends on the direction of the altitude of a triangle in addition to the direction of its base.
\end{example}

The strict gap in the exponents witnessed in Example \ref{ex:l5} is an essential feature of $\ell_p$ geometry when $p\neq 2$. To prove Theorems \ref{t:main} and \ref{t:main2}, we employ tools from functional analysis such as modulus of smoothness, modulus of convexity, and normalized duality mappings to carry out the estimates outlined in Remark \ref{r:pythag} and Example \ref{ex:l5} in the setting of $\ell_p$ in full generality. In fact, we work in arbitrary uniformly smooth and uniformly convex Banach spaces, which include the $\ell_p$ spaces when $1<p<\infty$. To establish the sufficient condition \eqref{e:suff-p}, we modify the ``parametric proof" of the sufficient condition in $\ell_2$ recently developed by Badger, Naples, and Vellis \cite{BNV}. This approach also yields new sufficient criteria for a set in a Banach space to be contained in a $(1/s)$-H\"older curve with $s>1$ (see Theorem \ref{t:param}). To verify the necessary condition \eqref{e:nec-2}, we follow the proof of the necessary condition in $\ell_2$ originally introduced by Schul \cite{Schul-Hilbert}, indicating which parts of the proof are metric, which parts are Banach, and which parts rely on the uniform convexity of the norm. The adaptation of the Hilbert space proofs to $\ell_p$ and related Banach spaces is non-trivial; see \S\ref{sec:smooth} (sufficient conditions) and \S\ref{sec:convex} (necessary conditions) for details. While revising an earlier draft of this manuscript, we encountered some mistakes in the proof of the necessary conditions in \cite{Schul-Hilbert}; see Remark \ref{the-mistake} and \cite[Appendix C]{Badger-McCurdy-2}. In \S\ref{ss:martingale}, we show how to correct the error in a minimal way, leaving the outline of virtually all of the original proofs intact; see Remark \ref{the-fix}. We complete the proof of Theorem \ref{t:main2} in \cite{Badger-McCurdy-2}.

To show that the exponents in Theorems \ref{t:main} and \ref{t:main2} are sharp and to prove Proposition \ref{t:pq}, we construct Koch-snowflake-like curves in \S\ref{sec:examples}. While the examples of sharpness when $\min(p,2)=2$ or $\max(2,p)=2$ can be built inside any two-dimensional subspace of $\ell_p$, the examples of sharpness when $\min(p,2)=p\neq 2$ or $\max(2,p)=p\neq 2$  require curves that extend outside of every finite-dimensional subspace of $\ell_p$. This is plausible, because we know the critical exponent in any finite-dimensional Banach space is 2 by the Analyst's Traveling Salesman Theorem in $\RR^n$ and bi-Lipschitz equivalence of norms in finite dimensions (see \S\ref{sec:finite-dimensional} for details). Thus, it is natural to expect examples showing sharpness of an exponent $p\neq 2$ to live in infinite dimensions.

Section \ref{sec:finite-dimensional} (the Analyst's TSP in finite-dimensional Banach spaces) and \S\ref{sec:examples} (examples) may be read independently of \S\S\ref{sec:smooth} and \ref{sec:convex} (proofs of the main theorems).

\subsection{Related work} \label{sec:related} The inception for this investigation is a recent paper of Edelen, Naber, and Valtorta \cite{ENV-Banach} that extends \emph{Reifenberg's Topological Disk Theorem} \cite{Reif} (also see \cite{DT}) from the Euclidean to infinite-dimensional Hilbert and Banach spaces. The original formulation of the theorem says that if a closed set $\Sigma\subset\RR^n$ is uniformly bilaterally $\delta(k,n)$-close\footnote{For a comparison of unilateral (Jones) versus bilateral (Reifenberg) flatness of a set, see \cite{lsa}.} to some $k$-dimensional affine plane at all locations in $\Sigma$ and on all sufficiently small scales, then $\Sigma$ is locally homoeomorphic to open subsets of $\RR^k$ (that is, $\Sigma$ is locally a topological disk). Reifenberg proved the Topological Disk Theorem to establish existence and regularity of the Plateau problem in arbitrary dimension and codimension. Underpinning the theorem is an algorithm that takes a collection of planes approximating the set $\Sigma$ and patches them together using orthogonal projections and partitions of unity to construct a parameterization. Edelen, Naber, and Valtorta solve the problem of how to implement this algorithm in a Banach space with dimension independent estimates. The main application of the Reifenberg algorithm in \cite{ENV-Banach} is a structure theorem for measures in Banach spaces, which we now briefly describe.

Following the convention used in \cite{ENV-Banach}, for every Borel regular measure $\mu$ on a Banach space $\XX$, location $x\in\XX$, and scale $r>0$, define the \emph{$k$-dimensional $L^2$ Jones beta number} $\beta_\mu^k(x,r)$ by \begin{equation}\beta_\mu^k(x,r)^2:= \inf_{p+V^k} r^{-k-2} \int_{B(x,r)} \dist(z,p+V)^2\,d\mu(z),\end{equation} where the infimum ranges over all $k$-dimensional affine subspaces of $\XX$. Beta numbers associated to measures were originally introduced by David and Semmes \cite{DS91,DS93} to build a bridge between singular integral operators and quantitative rectifiability of sets.

\begin{theorem}[Edelen, Naber, Valtorta {\cite[Theorem 2.1]{ENV-Banach}}] \label{t:env} Let $\XX$ be a Banach space, let $\mu$ be a finite Borel measure supported in $B(0,1)$, let $S\subset B(0,1)$ be a set with $\mu(B(0,1)\setminus S)=0$, and let $r:S\rightarrow(0,1)$. Assume that $\mu$ satisfies \begin{equation}\label{e:env} \int_{r(x)}^2 \beta_\mu^k(x,s)^\alpha\, \frac{ds}{s} \leq M^{\alpha/2}\quad\text{for all }x\in S,\end{equation} where the exponent $\alpha$ is determined as follows: \begin{itemize}
\item if $\XX$ is a generic Banach space, then $\alpha=1$;
\item if $\XX$ is a Hilbert space, then $\alpha=2$;
\item if $\XX$ is a uniformly smooth Banach space and $k=1$, then $\alpha$ is the smoothness power of $\XX$ (e.g.~ the smoothness power of $\ell_p$ is $\min(p,2)$ for $1<p<\infty$).
    \end{itemize}
Then there is a subset $S'\subset S$ so that we have the following packing/measure estimate: \begin{equation} \mu\left(B(0,1) \setminus \bigcup_{x \in S'} B(x,r(x))\right)\lesssim_{k,\rho_\XX} M\quad\text{and}\quad \sum_{x\in S'} r(x)^k \lesssim_{k,\rho_\XX} 1,\end{equation} where $\rho_{\XX}$ denotes the modulus of smoothness of $\XX$ (see \S\ref{sec:smooth} below).\end{theorem}

Roughly speaking, condition \eqref{e:env} allows one to control the tilt of the approximating planes in the Reifenberg algorithm and construct local bi-Lipschitz parameterizations. For our present discussion, the most interesting aspect of Theorem \ref{t:env} is the dependence of the exponent $\alpha$ on the geometry of the Banach space $\XX$ and the dimension $k$ of the approximating planes. In any given space, one would like to identify the largest possible exponent such that the theorem holds. The exponent $\alpha=1$ corresponds to the triangle inequality, which holds in any Banach space, and the exponent $\alpha=2$ corresponds to the Pythagorean theorem, which holds in any Hilbert space. In an intermediate scenario, Edelen, Naber, and Valtorta prove that when $\XX$ is a smooth Banach space \emph{and}  $k=1$, the exponent $\alpha$ can be taken to be the smoothness power of the Banach space. For example, $\alpha=p$ when $\XX=\ell_p$ and $1<p<2$. Furthermore, Edelen, Naber, and Valtorta prove that the restriction to $k=1$ is necessary to obtain $\alpha>1$ in non-Hilbert spaces. This is tied up with the existence of good projections onto lines when $k=1$ and the absence of good projections onto subspaces when $k\geq 2$; see \cite[\S3.6 and \S5.3]{ENV-Banach} for details.

A strength of the Reifenberg algorithm over the Analyst's Traveling Salesman Theorem is that it gives conditions to build parameterizations of \emph{every} dimension $k$. An advantage of the Analyst's Traveling Salesman Theorem over the Reifenberg algorithm is that it provides \emph{necessary and sufficient} conditions for parameterizations of dimension $k=1$. Edelen, Naber, and Valtorta's successful implementation of the Reifenberg algorithm in smooth Banach spaces with $k=1$ provided our initial motivation to look for an Analyst's Traveling Salesman Theorem in Banach spaces.

A separate vein of research by Hahlomaa \cite{Hah05} and David and Schul \cite{DS-metric} (also see \cite{Schul-AR,Hah08}) focuses on the Analyst's TSP in the setting of an arbitrary metric space. Because metric spaces are not necessarily path-connected, it is natural to reformulate the Analyst's TSP as stated above and instead ask which sets in a given metric space $\XX$ are contained in \emph{rectifiable curve fragments}, i.e.~ images of Lipschitz maps $f:S\rightarrow\XX$ from some set $S\subset[0,1]$.  Hahlomaa's original work in this direction established an analogue of the \emph{sufficient half} of the Analyst's Traveling Salesman Theorem by redefining Jones' beta numbers using Menger curvature, or equivalently, using the excess in the triangle inequality. For different perspectives on rectifiability in measure metric spaces, see e.g.~\cite{PT,Kirchheim,Ambrosio-Kirchheim,Bate,BCW,Bate-Li,Naples20}.

Following \cite{DS-metric}, for a given metric space $E$ and ball $Q=B(p,r)$, define $\beta^E_\infty(Q)$ by \begin{equation}\begin{split}\label{metric-beta} \beta^E_\infty(Q)^2&:= r^{-1}\sup\big\{ \dist(x,y)+\dist(y,z)-\dist(x,z):\\
&\qquad\qquad\quad\ \ x,y,z\in E\cap Q \text{ and } \dist(x,y)\leq \dist(y,z)\leq \dist(x,z)\big\}.\end{split}\end{equation} The \emph{metric beta number} $\beta^E_\infty(Q)$ measures the normalized excess in the triangle inequality among triples of points in $E\cap Q$. The exponent 2 on the left hand side of \eqref{metric-beta} is a convention that is imposed to make the statement of Theorem \ref{t:hah} look similar to Theorem \ref{t:Jones} when $E\subset\RR^N$ is endowed with the Euclidean metric.

\begin{theorem}[Hahlomaa {\cite[Theorem 5.3]{Hah05}}] \label{t:hah} Let $E$ be a metric space and let $\mathscr{G}$ be a multiresolution family for $E$ with inflation factor $A_\mathscr{G}\simeq 1$.  If \begin{equation}S^E_\infty(\mathscr{G}):=\diam E+ \sum_{Q\in\mathscr{G}} \beta^E_\infty(Q)^2\diam Q<\infty,\end{equation} then there exists a set $A\subset[0,1]$ and a surjective Lipschitz map $f:A\rightarrow E$ with Lipschitz constant $\Lip(f)\lesssim S^E_\infty(\mathscr{G})$.
\end{theorem}

\begin{remark} It is known that the converse to Hahlomaa's theorem is (quantitatively) false for certain rectifiable curves in $\ell^2_1=(\RR^2,|\cdot|_1)$; see \cite[Example 3.3.1]{Schul-survey}. A similar phenomenon occurs in graph inverse limit spaces; see \cite[\S7]{Guy-Schul}. This issue is not fully understood and merits further investigation.\end{remark}

David and Schul recently announced a partial converse to Hahlomaa's theorem, which is the first non-trivial necessary condition for the Analyst's TSP in a metric space. Together, Theorems \ref{t:hah} and \ref{t:DS-nec} are quite striking and indicate the rough shape that a full solution to the Analyst's TSP in a general metric space might take. Recall that a metric space is \emph{doubling} if every ball of radius $r$ can be covered by at most $D$ balls of radius $r/2$.

\begin{theorem}[David and Schul {\cite[Theorem A]{DS-metric}}]\label{t:DS-nec} Let $\Sigma$ be a connected, doubling metric space with doubling constant $D$ and let $\mathscr{H}$ be a multiresolution family for $\Sigma$ with inflation factor $A_\mathscr{H}>1$. For every $\epsilon>0$, \begin{equation} S^{\Sigma,\epsilon}_{\infty}(\mathscr{\mathscr{H}}):=\diam Q+\sum_{Q\in\mathscr{H}} \beta^\Sigma_\infty(Q)^{2+\epsilon}\diam Q\lesssim_{\epsilon,D,A_{\mathscr{H}}}\Haus^1(\Sigma).\end{equation}\end{theorem}

\begin{remark} The doubling assumption in Theorem \ref{t:DS-nec} allows the authors to simplify the overall proof of theorem. David and Schul conjecture (see \cite[Remark 1.6]{DS-metric}) that the doubling assumption can be dropped by implementing the techniques in \cite{Schul-Hilbert}.\end{remark}

David and Schul present several corollaries to Theorem \ref{t:DS-nec} with alternative definitions of the metric beta numbers $\beta^\Sigma_\infty(Q)$ (see \eqref{metric-beta}). In particular, they obtain necessary conditions for the Analyst's TSP in $\ell_p$, with $1<p<\infty$, using traditional Jones' beta numbers $\beta_\Sigma(Q)$ (see \eqref{def:beta}). Also see \cite[Corollary D]{DS-metric} for a more general statement on uniformly convex Banach spaces.

\begin{corollary}[David and Schul \cite{DS-metric}] \label{c:DS-nec}Let $1<p<\infty$, let $\Sigma\subset\ell_p$ be a connected set with doubling constant $D$, and let $\mathscr{H}$ be a multiresolution family for $\Sigma$ with inflation factor $A_\mathscr{H}>1$. For all $\epsilon>0$, \begin{equation}S_{E,\max(2,p)+\epsilon}(\mathscr{H})=\diam E + \sum_{Q\in\mathscr{H}} \beta_E(Q)^{\max(2,p)+\epsilon}\diam Q \lesssim_{\epsilon,D,p,A_\mathscr{H}} \Haus^1(\Sigma).\end{equation}\end{corollary}

In Theorem \ref{t:main2}, we remove the doubling assumption \emph{and} the error $\epsilon$ from Corollary \ref{c:DS-nec}. This is accomplished by following the strategy in \cite{Schul-Hilbert}.

\subsection*{Acknowledgements} The authors thank an anonymous referee for helpful comments that led us to improve the exposition in \S\ref{sec:convex} and include \S\ref{sec:finite-dimensional}. The referee initially brought the gap in \cite{Schul-Hilbert}, which we report in Remark \ref{the-mistake}, to our attention; for the correction, see Remark \ref{the-fix}. We thank Raanan Schul for encouraging us to write out the details in \S3 and the appendix in full and hope this will assist others in learning the proof. Finally, the first author would like to thank Raanan for his past mentorship and for imparting some of his good mathematical taste.

\section{Modulus of smoothness and proof of the sufficient conditions}\label{sec:smooth}

\subsection{Ordering flat sets in Banach spaces} \label{sec:order}

A simple, but important ingredient in all proofs of the Analyst's Traveling Salesman Theorem is that ``almost flat" sets of points can be linearly ordered.
To implement a generic Banach space version of the sufficient half of the Analyst's Traveling Salesman Theorem with universal constants (in the spirit of Hahlomaa \cite{Hah05}), we first develop an instance of this principle. The following lemma is modeled after \cite[Lemma 8.3]{BS3}.

\begin{lemma}[flatness implies order]\label{l:graph} Let $\XX$ be a Banach space. Suppose that $V\subset\XX$ is a $\delta$-separated set with $\#V\geq 2$ and there exist lines $L_1$ and $L_2$ and a number $0\leq \alpha<1/6$ such that \begin{equation}\dist(v,L_i)\leq \alpha\delta\quad\text{for all $v\in V$ and $i=1,2$}.\end{equation} Let $\pi_i:\XX\rightarrow L_i$ denote a metric projection onto $L_i$, i.e.~ any  map satisfying \begin{equation}\dist(x,L_i)=\dist(x,\pi_i(x))\quad\text{ for all }x\in\mathbb{X}.\end{equation} There exist compatible identifications of $L_1$ and $L_2$ with $\RR$ such that $\pi_1(v')\leq \pi_1(v'')$ if and only if $\pi_2(v')\leq \pi_2(v'')$ for all $v',v''\in V$. If $v_1,v_2\in V$, then \begin{equation}\label{e:graph1} (1+2\alpha)^{-1}|\pi_1(v_1)-\pi_1(v_2)|\leq |v_1-v_2|\leq (1+ 3\alpha)|\pi_1(v_1)-\pi_1(v_2)|.\end{equation}\end{lemma}

\begin{proof} Without loss of generality, it suffices to assume $\delta=1$. Let $V\subset\XX$ be a 1-separated set with at least two points. Assume that there exist one-dimensional affine subspaces $L_1$ and $L_2$ in $\XX$ and a number $0\leq \alpha < 1/6$ such that $$\dist(v,L_i)\leq \alpha\quad\text{for all $v\in V$ and $i=1,2$}.$$ Let $\pi_i$ denote a metric projection onto $L_i$. For any distinct pair of points $v_1,v_2\in V$, $$1\leq |v_1-v_2| \leq |\pi_{i}(v_1)-\pi_i(v_2)| +2\alpha,$$ because $V$ is 1-separated and the distance of points in $V$ to $L_i$ is bounded by $\alpha$. Hence \begin{equation}\label{e:graph3}|\pi_i(v_1)-\pi_i(v_2)|\geq |v_1-v_2|-2\alpha\geq 1-2\alpha> 2/3.\end{equation} In particular, $2\alpha\leq 3\alpha|\pi_i(v_1)-\pi_i(v_2)|$, and it follows that $$ |v_1-v_2| \leq (1+3\alpha)|\pi_i(v_1)-\pi_i(v_2)|.$$ This establishes the right half of \eqref{e:graph1}. Similarly, $$ |\pi_i(v_1)-\pi_i(v_2)|\leq |v_1-v_2|+2\alpha \leq (1+2\alpha)|v_1-v_2|,$$ which yields the left half of \eqref{e:graph1}. In particular, note that \begin{equation}\label{e:graph4} |v_1-v_2|\geq |\pi_i(v_1)-\pi_i(v_2)|-2\alpha.\end{equation}

Suppose for contradiction that there are identifications of $L_1$ and $L_2$ with $\RR$ and distinct points $v,v',v''\in V$ such that $\pi_1(v)<\pi_1(v')< \pi_1(v'')$, but $\pi_2(v')<\pi_2(v)<\pi_2(v'')$. Set \begin{align*}x&:=|v-v'|, \quad & y&:=|v''-v'|,\quad & z&:=|v''-v|, \\
x_1&:=|\pi_1(v)-\pi_1(v')|,\quad & y_1 &:=|\pi_1(v'')-\pi_1(v')|,\quad & z_1&:=|\pi_1(v'')-\pi_1(v)|, \\
x_2&:=|\pi_2(v)-\pi_2(v')|,\quad & y_2&:=|\pi_2(v'')-\pi_2(v')|,\quad & z_2&:=|\pi_2(v'')-\pi_2(v)|.\end{align*} Heuristically, since $\pi_1(v)<\pi_1(v')<\pi_1(v'')$, we have $z\approx x+y$, and since $\pi_2(v')<\pi_2(v)<\pi_2(v'')$, we have $y \approx x+z$. Hence $z\approx z+2x$, which yields a contradiction if $\alpha$ is sufficiently small. More precisely, by repeated application of \eqref{e:graph3} and \eqref{e:graph4}, \begin{equation*}\begin{split} z\geq z_1-2\alpha=x_1+y_1-2\alpha\geq x_1+y - 4\alpha &\geq x_1+y_2-6\alpha \\ &= x_1+x_2+z_2-6\alpha \geq x_1+x_2+z-8\alpha.\end{split}\end{equation*} Rearranging, we obtain $4/3<x_1+x_2\leq 8\alpha<8/6$, which is absurd. Therefore, under any choice of identifications of $L_1$ and $L_2$ with $\RR$, either $\pi_1(v)\leq \pi_1(v')$ if and only if $\pi_2(v)\leq \pi_2(v')$ for all $v,v'\in V$, or $\pi_1(v)\leq \pi_1(v')$ if and only if $\pi_2(v)\geq \pi_2(v')$ for all $v,v'\in V$. Thus, we can choose compatible identifications of $L_1$ and $L_2$ with $\RR$ such that $\pi_1(v')\leq \pi_1(v'')$ if and only if $\pi_2(v')\leq \pi_2(v'')$ for all $v',v''\in V$.
\end{proof}

\begin{corollary}[{cf.~\cite[Lemma 2.2]{BNV}}] \label{c:var} Let $\XX$ be a Banach space. Suppose that $V\subset\XX$ is a $\delta$-separated set with $\#V\geq 2$ and there exists a line $L$ and a number $0\leq \alpha< 1/6$ such that \begin{equation}\dist(v,L)\leq \alpha\delta\quad\text{for all $v\in V$.}\end{equation} Enumerate $V=\{v_1,\dots,v_n\}$ so that $v_{i+1}$ lies to the right of $v_i$ for all $1\leq i\leq n-1$. Then \begin{equation}\label{e:var2} \sum_{i=1}^{n-1} |v_{i+1}-v_i|^s \leq (1+3\alpha)^{2s}|v_1-v_n|^s\quad\text{for all $s\in[1,\infty)$.}\end{equation}
\end{corollary}

\begin{proof} Let $\pi$ denote a metric projection onto $L$. For all $1\leq i\leq n$, set $x_i:=\pi(v_i)$. Then $$(1+3\alpha)^{-1}|x_{i+1}-x_i|\leq |v_{i+1}-v_i| \leq (1+3\alpha)|x_{i+1}-x_i|\quad\text{for all } 1\leq i\leq n-1$$ by Lemma \ref{l:graph}. Assume $s\in[1,\infty)$ and $\#V\geq 2$. Then
\begin{equation*}\begin{split} \sum_{i=1}^{n-1} \frac{|v_{i+1}-v_i|^s}{(1+3\alpha)^s} &\leq \sum_{i=1}^{n-1} |x_{i+1}-x_i|^s \\ &\leq \left(\sum_{i=1}^{n-1}|x_{i+1}-x_i|\right)^s =|x_1-x_n|^s\leq (1+3\alpha)^s|v_1-v_n|^s,
\end{split}\end{equation*} because $s\geq 1$ and $x_1,\dots,x_n$ appear in the given order on the line $L$ by Lemma \ref{l:graph}. This proves \eqref{e:var2}.
\end{proof}

\begin{remark}In a general Banach space, the metric projection is not unique and may be norm-increasing. For example, in $\ell^2_\infty=(\RR^2,|\cdot|_\infty)$, consider the horizontal line $L$ through the origin (``the $x$-axis") and the point $v=(1,\alpha)$ for some $0<\alpha\leq 1$. Then $|v|_\infty=1$, $\dist(v,L)=\alpha$, and a metric projection from $v$ to $L$ can be any point on the line segment $[1-\alpha,1+\alpha]\times\{0\}$. In particular, if $\pi_L(v)=(1+\alpha,0)$, then $|\pi_L(v)|_\infty = (1+\alpha)|v|_\infty>|v|_\infty$. This shows that in Lemma \ref{l:graph} for an arbitrary Banach space, we cannot expect to replace the lower bound in \eqref{e:graph1} with a $1$-Lipschitz bound.\end{remark}

\subsection{Lipschitz and H\"older continuous Traveling Salesman parameterizations in Banach spaces} \label{sec:bnv}

Throughout this section, let $\XX$ denote an arbitrary Banach space.

\begin{definition}[doubling scales] Let $0<\xi_1\leq \xi_2<1$. A \emph{$(\xi_1,\xi_2)$-doubling sequence of scales} is a sequence $(\rho_k)_{k=0}^\infty$ of positive numbers such that $\rho_0=1$ and for all $k\geq 0$, $\xi_1\rho_k \leq \rho_{k+1} \leq \xi_2\rho_k$. In the special case when $\xi_1=\xi_2$, we may call $(\rho_k)_{k=0}^\infty$ a \emph{geometric sequence of scales}.\end{definition}

Following \cite{BNV}, let $\mathscr{V}=(V_k,\rho_k)_{k=0}^\infty$ be a sequence consisting of nonempty finite sets $V_k$ in $\XX$ and positive numbers $\rho_k$. Assume that there exist $x_0\in \XX$, $r_0>0$, $C^*\geq 1$, and $0<\xi_1\leq \xi_2<1$ such that $\mathscr{V}$ has the following properties. \begin{enumerate}
\item[(V0)] The numbers $(\rho_k)_{k=0}^\infty$ are a $(\xi_1,\xi_2)$-doubling sequence of scales.
\item[(V1)] When $k=0$, we have $V_0 \subset B(x_0, C^*r_0)$.
\item[(V2)] For all $k\geq 0$, we have $V_k\subset V_{k+1}$.
\item[(V3)] For all $k\geq 0$ and all distinct $v,v'\in V_k$, we have $|v-v'|\geq \rho_k r_0$.
\item[(V4)] For all $k \geq 0$ and all $v\in V_{k+1}$, there exists $v'\in V_k$ such that $|v-v_k|<C^*\rho_{k+1}r_0$.
\end{enumerate} With $C^*$ and $\xi_2$ given, define the associated constant \begin{equation}A^*:= \frac{C^*}{1-\xi_2}.\end{equation} In addition to (V0)--(V4), assume that for each $k\geq 0$ and $v\in V_k$, we are given a number $\alpha_{k,v}\geq 0$ and a line $L_{k,v}$ in $\XX$ such that  \begin{gather}\tag{V5} \sup_{x\in V_{k+1}\cap B(v,30A^*\rho_k r_0)} \dist (x,L_{k,v})\leq \alpha_{k,v}\rho_{k+1}r_0.\end{gather}

\begin{definition}[flat pairs, see \cite{BNV}] \label{def:flatpairs} Fix a parameter $\alpha_0\in(0,1/6)$. For all $k\geq 0$, define $\Flat(k)$ to be the set of pairs $(v,v')\in V_k\times V_k$ such that \begin{enumerate}
\item[(F1)] $\rho_k r_0 \leq |v-v'| < 14A^*\rho_k r_0$, and
\item[(F2)] $\alpha_{k,v}<\alpha_0$ and $v'$ is the first point in $V_k\cap B(v,14A^*\rho_kr_0)$ to the left or to the right of $v$ with respect to the ordering induced by $L_{k,v}$.\end{enumerate}
\end{definition}

Given a pair $(v,v')\in\Flat(k)$, let $V_{k+1}(v,v')$ denote the set of all points $x\in V_{k+1}\cap B(v,14A^*\rho_k r_0)$ such that $x$ lies between $v$ and $v'$ (inclusive) with respect to the ordering induced by $L_{k,v}$.

\begin{definition}[variation excess, see \cite{BNV}] For all $s\in[1,\infty)$, for all $k\geq 0$, and for all $(v,v')\in\Flat(k)$, define the \emph{$s$-variation excess} $\tau_s(k,v,v')$ by \begin{equation}\tau_s(k,v,v')|v-v'|^s = \max\left\{ \left(\sum_{i=1}^{n-1}|v_{i+1}-v_i|^s\right)-|v-v'|^s,0\right\},\end{equation} where $V_{k+1}(v,v')=\{v_1,\dots,v_n\}$ is enumerated so that $v_1=v$ and for all $1\leq i\leq n-1$, $v_{i+1}\in V_{k+1}(v,v')$ is the first point after $v_i$ in the direction from $v$ to $v'$ with respect to the ordering induced by $L_{k,v}$ (hence $v_n=v'$).\end{definition}

The following theorem extends Badger, Naples, and Vellis; see \cite[Theorem 5.1]{BNV}. In its original form, the theorem was stated for $\XX=\ell_2$ with a weaker restriction on $\alpha_1$, achieved through targeted use of the Pythagorean theorem.

\begin{theorem}[H\"older Traveling Salesman parameterizations for nets in Banach spaces] \label{t:param} Let $\XX$ be a Banach space. In addition to (V0)--(V5), assume that  \begin{equation}\alpha_0 \leq \alpha_1:=\frac{\xi_1(1-\xi_2)}{87C^*}.\end{equation}  If the sum \begin{equation} \label{net-sum} S^s_\mathscr{V}:= \sum_{k=0}^\infty \sum_{{(v,v')\in\Flat(k)}}\tau_s(k,v,v')\rho_k^s+\sum_{k=0}^\infty \sum_{{\stackrel{v\in V_k}{\alpha_{k,v}\geq \alpha_0}}}\rho_k^s<\infty,\end{equation} then there exists a $(1/s)$-H\"older continuous map $f:[0,1]\rightarrow\XX$ such that $f([0,1])$ contains $\bigcup_{k\geq 0}V_k$ and the $(1/s)$-H\"older constant of $f$ satisfies $H\lesssim_{s,C^*,\xi_1,\xi_2} r_0(1+S^s_\mathscr{V}).$
\end{theorem}

\begin{proof} Repeat the proof of \cite[Theorem 5.1]{BNV} (see \cite[\S\S2--5]{BNV}) \emph{mutatis mutandis}. Use Lemma \ref{l:graph} and Corollary \ref{c:var} above in place of \cite[Lemmas 2.1 and 2.2]{BNV}.\end{proof}

\begin{remark}[outline of key modifications] We specify some details to aid the reader with the proof of Theorem \ref{t:param}. The reader is first urged to read through \cite[\S2]{BNV}, followed by item (C0) in \cite[\S3.10]{BNV}, which indicates how key parameters are chosen.

In Hilbert space, the initial upper bound $1/16$ on the size of $\alpha_0$ is made so that $$1+3\alpha_0^2\leq1+3(1/16)^2<1.1.$$ See \cite[Lemma 2.3]{BNV}. In generic Banach space, let us initially require $\alpha_0\leq 1/31$ so that $1+3\alpha_0\leq 1+3(1/31)<1.1.$ Then using Lemma \ref{l:graph} instead of \cite[Lemma 2.1]{BNV}, the rest of the proof and conclusion of \cite[Lemma 2.3]{BNV} goes through as written.

The algorithm presented in \cite[\S3]{BNV} requires no changes in the Banach setting.

The principal estimates in the proof of the theorem occur in \cite[\S4]{BNV}. No changes are required until we reach the proof of \cite[Lemma 4.6]{BNV}, where we need to use Corollary \ref{c:var} instead of \cite[Lemma 2.2]{BNV}. This time, we require that $\alpha_0\leq 1/62$ so that we can replace the original estimate $1+3\alpha_0^2<1.1$ with the estimate $(1+3\alpha_0)^2<1.1.$ The next required change occurs at the end of the proof of \cite[Lemma 4.9]{BNV}. Using Corollary \ref{c:var} once again, we see the original requirement $1+3\alpha_0^2-\xi_1/14A^*\leq 1$ becomes $(1+3\alpha_0)^2-\xi_1/14A^*\leq 1$, or equivalently $6\alpha_0+9\alpha_0^2\leq \xi_1/14A^*$. With our \emph{a priori} bound $\alpha_0\leq 1/62$, this certainly holds provided $(6+9/62)\alpha_0 \leq \xi_1/14A^*$. Thus, after noting that  $(6+9/62)14=86.03...$, it suffices to take $$\alpha_1 \leq \frac{\xi_1}{87A^*}=\frac{\xi_1(1-\xi_2)}{87C^*}.$$ Note that $$\min\left\{\frac{1}{31},\frac{1}{62}, \frac{\xi_1(1-\xi_2)}{87C^*}\right\}=\frac{\xi_1(1-\xi_2)}{87C^*},$$ because $0<\xi_1\leq \xi_2<1$ and $C^*\geq 1$. There are two final uses of Lemma \ref{l:graph} instead of \cite[Lemma 2.1]{BNV} to estimate the separation of points after projection onto an approximating line $\ell_{k,v}$, once in the proof of \cite[Proposition 4.11]{BNV} and once in the proof of \cite[(4.3)]{BNV}. This change affects the value of the implicit constant in \cite[Proposition 4.11]{BNV}, but not dependencies of the constant.

To finish the proof of Theorem \ref{t:param}, repeat the argument in \cite[\S5.1]{BNV} verbatim.\end{remark}

\begin{corollary} \label{c:param1} Let $\XX$ be a Banach space. Assume  $\mathscr{V}=(V_k,\rho_k)_{k=0}^\infty$ satisfies (V0)--(V5) above. If the sum \begin{equation}S_\mathscr{V}:=\sum_{k=0}^\infty \sum_{v\in V_k} \alpha_{k,v} \rho_k<\infty,\end{equation} then there exists a rectifiable curve $\Gamma$ containing $\bigcup_{k\geq 0}V_k$ such that \begin{equation}\Haus^1(\Gamma)\lesssim_{C^*,\xi_1,\xi_2} r_0(1+S_\mathscr{V}).\end{equation}\end{corollary}

\begin{proof} Set $\alpha_0=\alpha_1$, which depends only on $\xi_1$, $\xi_2$, and $C^*$. By Corollary \ref{c:var} with $s=1$, we have $\tau_1(k,v,v') \leq 6\alpha_{k,v}+9\alpha_{k,v}^2\leq 7\alpha_{k,v}$ for every flat pair $(v,v')\in\Flat(k)$. Thus, \begin{align*}
S^1_\mathscr{V} &= \sum_{k=0}^\infty \sum_{(v,v')\in\Flat(k)}\tau_1(k,v,v')\rho_k+\sum_{k=0}^\infty \sum_{\stackrel{v\in V_k}{\alpha_{k,v}\geq \alpha_1}}\rho_k\\
                &\leq \sum_{k=0}\sum_{(v,v')\in\Flat(k)}7\alpha_{k,v} \rho_k+\alpha_1^{-1} \sum_{\stackrel{v\in V_k}{\alpha_{k,v}\geq \alpha_1}} \alpha_{k,v}\rho_k
               \leq \alpha_1^{-1} S_\mathscr{V}<\infty.
\end{align*}
By Theorem \ref{t:param}, there exists a Lipschitz map $f:[0,1]\rightarrow\XX$ such that $\Gamma:=f([0,1])$ contains $\bigcup_{k=0}^\infty V_k$ and $\Haus^1(\Gamma)\leq \Lip(f)\lesssim_{C^*,\xi_1,\xi_2} r_0(1+S^1_\mathscr{V})\lesssim_{C^*,\xi_1,\xi_2} r_0(1+S_\mathscr{V})$.
\end{proof}

For completeness, we show how to use Corollary \ref{c:param1} to derive a beta number criterion for a set in a Banach space to be contained in a rectifiable curve. The following theorem is best viewed as the Banach space analogue of Theorem \ref{t:hah}, expressed with the geometric Jones' beta numbers \eqref{def:beta} instead of metric beta numbers \eqref{metric-beta}. To recall the definition of the sum $S_{E,1}(\mathscr{G})$ of beta numbers over a multiresolution family $\mathscr{G}$ for $E$, see \eqref{e:Jones-sum}.

\begin{theorem}[sufficient half of Schul's theorem in arbitrary Banach spaces] \label{t:banach-suff} Let $\XX$ be a Banach space. If $E\subset\XX$ and $S_{E,1}(\mathscr{G})=\diam E+\sum_{Q\in\mathscr{G}} \beta_E(Q) \diam Q<\infty$ for some multiresolution family $\mathscr{G}$ for $E$ with inflation factor $A_\mathscr{G}\geq 240$, then $E$ is contained in a rectifiable curve $\Gamma$ in $\XX$ with \begin{equation}\label{e:suff-banach}\Haus^1(\Gamma) \lesssim_{A_\mathscr{G}} S_{E,1}(\mathscr{G}).\end{equation}\end{theorem}

\begin{proof} Let $\XX$ be a Banach space, let $E\subset\XX$, let $\mathscr{G}$ be a multiresolution family for $E$ with inflation factor $A_\mathscr{G}\geq 240$, and assume that $S_{E,1}(\mathscr{G})<\infty$. Then $E$ is bounded and there exists a unique integer $k_0\in\ZZ$ such that \begin{equation}\label{e:diam-choice} 2^{-k_0}\leq \diam E< 2\cdot 2^{-k_0}.\end{equation} For all $k\geq 0$, define $\rho_k=2^{-k}$ and $V_k=X_{k_0+k}$, where $(X_j)_{j\in\ZZ}$ are the $2^{-j}$-nets for $E$ used to define $\mathscr{G}$. Set parameters $C^*=2$, $\xi_1=\xi_2=\frac12$, and $r_0=2^{-k_0}$, and choose any $x_0\in V_0=X_{k_0}$. Then the sequence $\mathscr{V}=(V_k,\rho_k)_{k=0}^\infty$ satisfies properties (V0)--(V4) above. Note that $A^*=4C^*=8$ and $30A^*=240$, since $\xi_1=\xi_2=\frac12$. For each $k\geq 0$ and $v\in V_k$, set $\alpha_{k,v}= 8A_\mathscr{G} \beta_E(B(v,A_\mathscr{G} 2^{-(k_0+k)}))$ and choose $L_{k,v}$ to be any line such that \begin{equation} \sup_{x\in E\cap B(v,A_\mathscr{G} 2^{-(k_0+k)})} \dist(x,L_{k,v}) \leq 2 \beta_E(B(v,A_\mathscr{G} 2^{-(k_0+k)})) \diam B(v,A_\mathscr{G} 2^{-(k_0+k)}).\end{equation}
Because $A_\mathscr{G}\geq 240=30A^*$, it follows that for all $k\geq 0$ and $v\in V_k$, \begin{equation}\begin{split} \sup_{x\in V_{k+1}\cap B(v,30A^* \rho_k r_0)} \dist(x,L_{k,v}) &\leq 2 \beta_E(B(v,A_\mathscr{G} 2^{-(k_0+k)})) \diam B(v,A_\mathscr{G} 2^{-(k_0+k)})\\
&\leq 8A_\mathscr{G}\beta_E(B(v,A_\mathscr{G} 2^{-(k_0+k)}))2^{-(k_0+k+1)}\\ &=\alpha_{k,v}\rho_{k+1}r_0.\end{split}\end{equation} Thus, property (V5) is satisfied, as well. To proceed, observe that \begin{equation}S_\mathscr{V} = \sum_{k=0}^\infty \sum_{v\in V_k}\alpha_{k,v} \rho_k\\
=\sum_{k=0}^\infty \sum_{x\in V_k} 8A_\mathscr{G}\beta_E(B(v,A_\mathscr{G} 2^{-(k_0+k)})) 2^{-k}\leq 8A_\mathscr{G} 2^{k_0}S_{E,1}(\mathscr{G}).\end{equation} Since $S_\mathscr{V} \leq 8A_{\mathscr{G}}S_{E,1}(\mathscr{G})<\infty$, there is a rectifiable curve $\Gamma$ containing $\bigcup_{k=0}^\infty V_k$ such that \begin{equation}\Haus^1(E) \lesssim_{C^*,\xi_1,\xi_2} r_0(1+S_\mathscr{V}) \lesssim r_0(1+8A_\mathscr{G}r_0^{-1}S_{E,1}) \lesssim_{A_\mathscr{G}} S_{E,1}\end{equation} by Corollary \ref{c:param1} and \eqref{e:diam-choice}. Finally, note that since $(V_k)_{k=0}^\infty$ is a sequence of $2^{-(k_0+k)}$-nets for $E$ and $\Gamma$ is closed, $\Gamma$ contains the set $\overline{\bigcup_{k=0}^\infty V_k}\supset E$, as well.
\end{proof}

\begin{remark} The constant 240 in Theorem \ref{t:banach-suff} has not been optimized and can be at least somewhat reduced at the cost of growing the implicit constant in \eqref{e:suff-banach}. In the future event that a smaller constant is needed, the reader should first consult \cite[\S3.10]{BNV} \end{remark}

\subsection{Triangle inequality excess in uniformly smooth Banach spaces}\label{ss:smooth} Our goal in this section is to prove that in a uniformly smooth Banach space of power type $p\in(1,2]$, the exponent $1$ in Corollary \ref{c:param1} and Theorem \ref{t:banach-suff} may be replaced with the exponent $p$. In the process, we will verify \eqref{e:suff-p} in Theorem \ref{t:main}. The essential step is to improve the exponent in the bound \eqref{e:graph1} in Lemma \ref{l:graph}. To achieve this, we follow the strategy used by Edelen, Naber, and Valtorta \cite{ENV-Banach} in their proof of one-dimensional Reifenberg-type theorems in uniformly smooth Banach spaces. The approach utilizes a special projection operator, which is available in uniformly smooth Banach spaces.

\begin{definition}\label{modulus of smoothness} \label{def:modsmooth} Let $\XX$ be a Banach space. The \textit{modulus of smoothness $\rho_\XX$} of $\XX$ is the function $\rho_\XX: [0, \infty) \rightarrow [0, \infty)$ defined by
\begin{equation}
\rho_\XX(t) := \sup_{|x|=1, |y|=t} \frac{1}{2}(|x+y|+|x-y|)-1\quad\text{for all $t\in[0,\infty)$}.\end{equation}
\end{definition}

\begin{definition} We say that $\XX$ is \emph{uniformly smooth} if $\rho_\XX(t)/t=o(t)$ as $t\rightarrow 0$. In this case, we say that $\XX$ is \emph{smoothness power type }$p\in(1,2]$ if there exists $C>0$ such that $\rho_\XX(t)\leq C t^p$ for all $t>0$.
\end{definition}

\begin{remark}[essential facts]\label{rho remark} For general background on the modulus of smoothness and uniformly smooth spaces, including the following inequalities, see e.g.~ \cite[Chapter Three]{Diestel} or \cite[Chapter 1, \S e]{LT-II}. On any Banach space $\XX$, the modulus of smoothness $\rho_\XX$ is a non-decreasing convex function such that $\rho_\XX(0)=0$ and
\begin{align} \label{basic-rho-bounds}
\sqrt{1 + t^2} - 1 =\rho_{\ell_2}(t)\le \rho_{\XX}(t) \le t\quad\text{for all }t\geq 0.
\end{align} Since $\rho_\XX$ is convex and $\rho_\XX(0)=0$, \begin{equation} \frac{\rho_\XX(t_1)}{t_1} \leq \frac{\rho_\XX(t_2)}{t_2}\quad\text{for all } 0<t_1\leq t_2.\end{equation} Furthermore, there exists a constant $1<L_0<3.18$ (see Figiel \cite[Proposition 10]{Figiel}) such that  \begin{equation}\label{e:figiel}\frac{\rho_\XX(t_2)}{t_2^2} \leq L_0\frac{\rho_\XX(t_1)}{t_1^2}\quad\text{for all }0<t_1\leq t_2,\end{equation} The modulus of smoothness $\rho_\XX$ and modulus of convexity $\delta_{\XX}$ (see \S\ref{ss:convex}) are related by \begin{equation} \rho_{\XX^*}(t) = \sup\left\{\frac{1}{2}t\epsilon - \delta_\XX(\epsilon):0<\epsilon\leq 2\right\}\quad\text{for all }t>0.\end{equation}
Hence the dual $\XX^*$ of a uniformly convex Banach space $\XX$ (see \S\ref{ss:convex}) is uniformly smooth. Finally, every uniformly smooth Banach space is reflexive.
\end{remark}

\begin{example}\label{example:smooth} By Hanner's inequalities \cite{Hanner}, $\rho_{L^p}(t) = p^{-1}t^p+o(t^p)$ when $1<p\leq 2$; and $\rho_{L^p}(t)=\frac12(p-1)t^2+o(t^2)$ when $2\leq p<\infty$. In particular, the $L^p$ spaces are uniformly smooth with power type $\min(p,2)$ when $1<p<\infty$.\end{example}

We now present a class of 1-Lipschitz projections onto a line in a Banach space. Given a real Banach space $\XX$, let $\XX^*$ denote the dual of $\XX$ and let $J: \XX \rightarrow \XX^*$ denote a \textit{normalized duality mapping}, i.e.~ a (nonlinear) map satisfying
\begin{equation}\label{J-properties}
|J(x)|_{\XX^*} = |x| \quad\text{and}\quad \langle J(x), x \rangle = |x|^2\quad\text{for all }x\in \XX,
\end{equation} where $\langle f,x\rangle\equiv f(x)\in\RR$ denotes the natural pairing of $f\in \XX^*$ and $x\in\XX$. Alternatively, $J$ is a subgradient of the convex function $x\in\XX\mapsto \frac{1}{2}|x|^2$ (see \cite{Asplund,Kien02}). The norm on any (uniformly) smooth Banach space $\XX$ is Gateaux (uniformly Fr\'echet) differentiable, and thus, $J$ is uniquely determined (see e.g.~\cite[Chapter Two]{Diestel}) when $\XX$ is smooth. For example, when $\XX=\ell_p$ with $1<p<\infty$,  $$J(x)= |x|_{\ell_p}^{2-p}y\in \ell_p^*=\ell_{p'},$$ where $y=(|x_1|^{p-2}x_1,|x_2|^{p-2}x_2,\dots)$ and $p'$ is the conjugate exponent to $p$.

\begin{definition}[{\cite[Definition 3.31]{ENV-Banach}}] \label{def:J-proj}
Let $\XX$ be a Banach space and let $L$ be a one-dimensional linear subspace of $\XX$. Define the \emph{$J$-projection} $\Pi_L$ onto $L$ by
\begin{equation}\Pi_{L}(x) := \langle J(v), x \rangle v\quad\text{for all }x\in\XX,\end{equation} where $J$ is a normalized dual mapping and $v$ is a point in $L$ with $|v|=1$. When $L$ is a one-dimensional affine subspace of $\XX$, define $\Pi_L\equiv p+\Pi_{L-p}(\cdot-p)$ for any choice of $p\in L$. For all lines $L$, we also define $\Pi_L^{\perp} \equiv \mathrm{Id}_\XX-\Pi_L$.
\end{definition}

Let us record some elementary properties of $J$-projection.

\begin{lemma}\label{j-proj facts}
Let $\XX$ be a Banach space and let $L$ be a one-dimensional linear subspace of $\XX$.  The $J$-projection $\Pi_{L}$ satisfies each of the following properties.
\begin{enumerate}
\item For all $x \in L$, we have $\Pi_{L}(x) = x$.
\item For all $x \in L^{\perp} \equiv \Pi_{L}^{\perp}(\XX)$, we have $\langle Jv, x \rangle = 0$ and $\Pi_L(x)=0$.
\item For all $x \in \XX$, we have $|\Pi_{L}(x)| \le |x|.$
\item For all $x \in \XX$, we have $\dist(x,L)\leq |\Pi_L^\perp(x)| \le 2 \dist(x, L)$.
\end{enumerate}
\end{lemma}

\begin{proof}
Let $v \in L$ be the unit vector in the definition of $\Pi_{L}$.  If $x \in L$, say $x = cv$, then $$\Pi_L(x)= \langle Jv,cv\rangle v = c\langle Jv,v\rangle v = c|v|^2v=cv=x.$$ This gives the first point. To see the second point, for any $x \in L^{\perp}$, say $x=\Pi_L^\perp(y)$,
\begin{equation*}\begin{split}
\langle Jv, x \rangle  = \langle Jv, y - \langle Jv, y \rangle v \rangle &= \langle Jv,y \rangle - \langle Jv,y \rangle \langle Jv,v \rangle\\
& = \langle Jv,y \rangle - \langle Jv,y \rangle |v |^2 = \langle Jv,y \rangle - \langle Jv,y \rangle = 0.
\end{split}\end{equation*} Hence $\Pi_L(x)=\langle Jv,x\rangle v=0v=0$.
For the third point, observe that for any $x\in\XX$,
\begin{align*}
|\Pi_L(x)|=|\langle Jv, x \rangle v| \le |Jv |_{\XX^*} |x||v| = |v |^2 |x |  =|x|.
\end{align*} To see the last point, suppose that $x\in \XX$. Clearly, $|\Pi_L^\perp(x)|=|\Pi_L(x)-x|\geq \dist(x,L)$. Choose $y \in L$ such that $|x - y| = \dist(x, L).$ Then
\begin{equation*}\begin{split}
|\Pi_{L}(x) - x|  &\le |\Pi_{L}(x) - y| + |x -y| \\ &= |\Pi_{L}(x) - \Pi_{L}(y)|  + |x -y|\le 2|x - y| = 2 \dist(x, L). \qedhere \end{split}\end{equation*}
\end{proof}

We now check that in any Banach space, the $J$-projection $\Pi_L$ induces a well-defined order on a sufficiently flat, separated set of points, by checking compatibility with the order induced by a metric projection $\pi_L$. The importance of this fact for us is that in the definition of flat pairs in the traveling salesman algorithm (see Definition \ref{def:flatpairs}), it does not matter whether we order points by $\pi_L$ or $\Pi_L$.

\begin{lemma}[order compatibility for $\Pi_L$ and $\pi_L$] \label{l:order-p} Let $\XX$ be a Banach space. Suppose that $V$ is a $\delta$-separated set with $\#V\geq 2$ and there exists a line $L$ and a number $0\leq \alpha<1/8$ such that $\dist(x,L)\leq \alpha\delta$ for all $x\in V.$ Then the $J$-projection $\Pi_L$ induces an order on $V$ compatible with the order induced by a metric projection $\pi_L$ onto $L$.
\end{lemma}

\begin{proof} Let $V\subset\XX$ be a $\delta$-separated set with $\#V\geq 2$, let $L$ be a line in $\XX$, let $0\leq \alpha<1/6$, and assume that $\dist(x,L)\leq \alpha \delta$ for all $x\in V$. Without loss of generality, we may assume that $L$ is a linear subspace of $\XX$. Fix any metric projection $\pi_L$ onto $L$. The restriction on $\alpha$ ensures that $\pi_L$ induces a unique order on $V$ by Lemma \ref{l:graph}. If $x,y\in V$ are distinct, then \begin{equation}\label{Pi-L-separated} |x-y|\geq |\Pi_L(x)-\Pi_L(y)| \geq |x-y| - |\Pi_L^\perp(x)| - |\Pi_L^\perp(y)| \geq (1-4\alpha)\delta > \frac13\delta>0\end{equation}  by the triangle inequality and Lemma \ref{j-proj facts}.
We may now check that the $J$-projection $\Pi_L$ induces an order on $V$ that is compatible with the order induced by $\pi_L$. To that end, suppose that $x,y,z\in V$ are distinct points, write $x'=\pi_L(x)$, $y'=\pi_L(y)$, $z'=\pi_L(z)$ and $x''=\Pi_L(x)$, $y''=\Pi_L(y)$, $z''=\Pi_L(z)$, and suppose to get a contradiction that there exists identifications of $L$ with $\RR$ such that $x'<y'<z'$ and $y''<x''<z''$. Heuristically, because of the order of the triples on the line $L$, $$|x-z| \approx |x-y|+|y-z| \approx |x-y|+\underbrace{|x-y|+|x-z|},$$ which is impossible if $\alpha$ is sufficient small. More precisely, on one hand, by \eqref{e:graph3} and \eqref{e:graph4} from the proof of Lemma \ref{l:graph} (recalling that there we normalized $\delta=1$), \begin{equation*}|x-z| \geq |x'-z'| -2\alpha\delta =|x'-y'|+|y'-z'| - 2\alpha\delta \geq |x-y|+|y-z|-6\alpha\delta,\end{equation*} since $x'<y'<z'$. On the other hand, by \eqref{Pi-L-separated}, \begin{equation*} |y-z| \geq |y''-z''| = |y''-x''|+|x''-z''| \geq |y-x|+|x-z| - 8\alpha\delta,\end{equation*} since $y''<x''<z''$. Combining the previous two displayed equations and recalling $V$ is $\delta$-separated, we obtain $$|x-z| \geq 2|x-y|+|x-z|-14\alpha\delta\geq |x-z| + 2\delta - 14\alpha\delta,$$ or equivalently $(2-14\alpha)\delta \leq 0$. This is a contradiction, because $\alpha<1/8$ implies that $2-14\alpha > 2-14/8 =1/4>0$. It readily follows that there is an identification of $L$ with $\RR$ such that $\pi_L(x)\leq \pi_L(y)$ if and only if $\Pi_L(x)\leq \Pi_L(y)$ for all $x,y\in V$.
\end{proof}

\begin{figure}
\begin{center}\includegraphics[width=.4\textwidth]{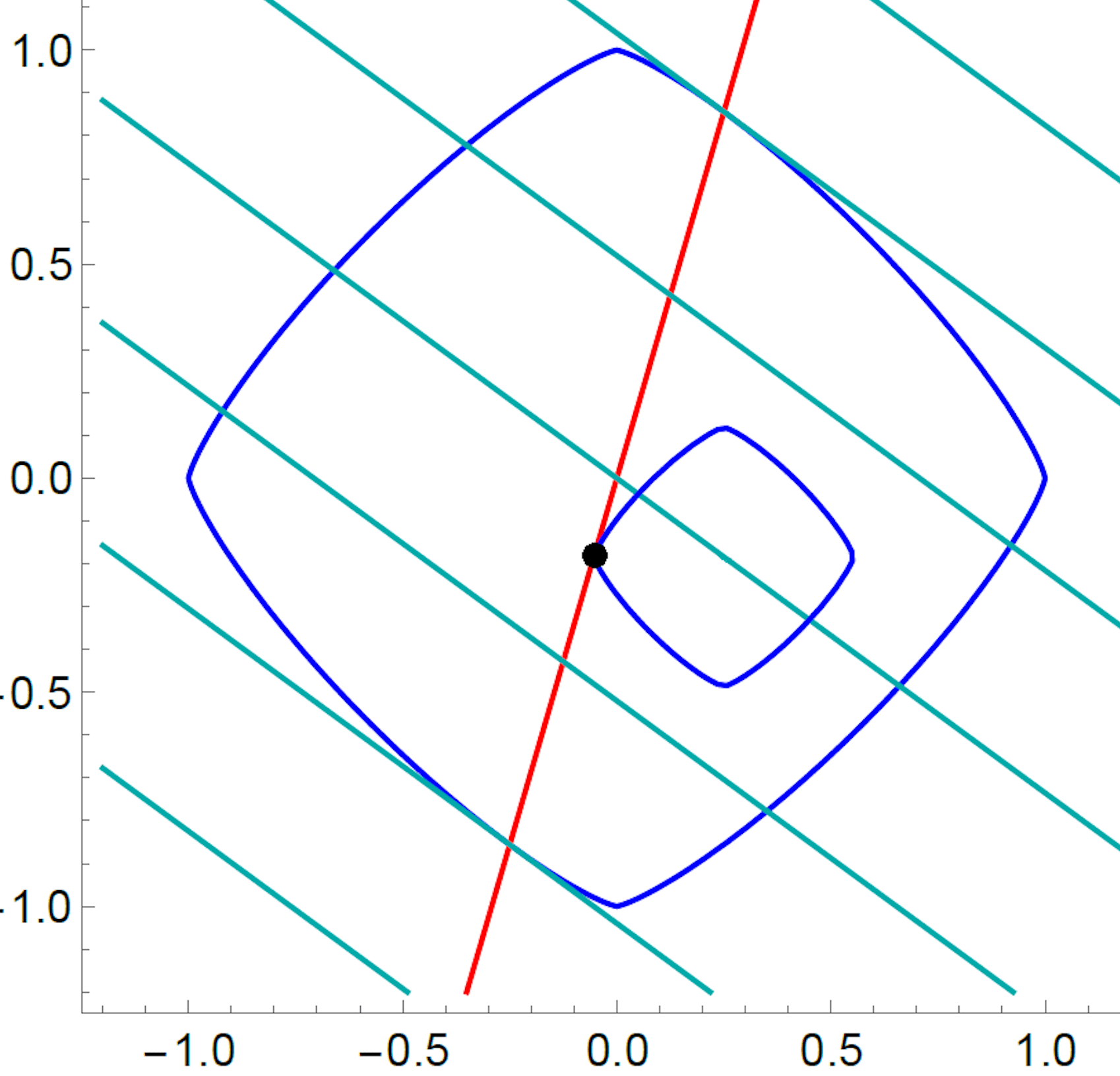}\end{center}
\caption{Fibers (cyan) of the $J$-projection $\Pi_{L}$ of $\XX=\ell_{5/4}^2$ onto the (red) line $L=\Span(v)$, where $v=(1/4,y_0)$, $(1/4)^{5/4}+y_0^{5/4}=1$. The dot (black) is the metric projection $\pi_L(w)$ onto $L$ of $w=(1/4,y_1)$ with $\Pi_L(w)=(0,0)$.}\label{fig:J-proj}
\end{figure}

\begin{remark}[geometric interpretation in smooth spaces] \label{remark J and tangent interpretation} Assume that $\XX$ is smooth. For a line $L$, spanned by a unit vector $v \in \XX$, the $J$-projection $\Pi_{L}$ admits the following geometric description. Let $T_{v}\partial B(0,1)$ denote the tangent hyperplane to $\partial B(0,1)$ at the point $v$ (which exists because $\XX$ is smooth). By Lemma \ref{j-proj facts}(2), a point $x\in\XX$ satisfies $\Pi_{L}(x) = cv$ if and only if $x \in (cv + T_{v}\partial B(0,1))$. See Figure \ref{fig:J-proj}.\end{remark}

\begin{lemma}\label{lemma:magic} If $\XX$ is a smooth Banach space, then $$\frac{d}{dt}|x+ty|^2= 2\langle J(x+ty),y\rangle$$ for all $x,y\in\XX$ and $t\in\RR$ with $x+ty\neq 0$ and $y\neq 0$. \end{lemma}

\begin{proof} Since $2J$ is the subgradient of $x\in\XX\mapsto |x|^2$, we have $$|z|^2-|x+ty|^2\geq \langle 2J(x+ty),z-(x+ty)\rangle\quad\text{for all }z\in\XX.$$ The claim follows by applying the inequality to the difference quotient $$\frac{|x+(t +h)y|^2-|x+ty|^2}{h}$$ along values $h\rightarrow 0^+$ and for $h\rightarrow 0^-$, where the limit exists by Gateaux differentiability of the norm.\end{proof}

The following estimate by Alber is crucial for our application below.

\begin{lemma}[see Alber {\cite[Remark 7.3]{Alber96} or \cite[(2.13)]{Alber2}}] \label{l:alber} Let $\XX$ be a uniformly smooth Banach space and write $h_\XX(t)\equiv \rho_\XX(t)/t$. Then \begin{equation}\label{e:alber} |J(x)-J(y)|_{\XX^*} \leq 4C_0 h_\XX(8C_0L_0|x-y|)\quad\text{for all }x,y\in\XX,\end{equation} where $C_0=2\max\left(1,\sqrt{\frac12(|x|^2+|y|^2))}\right)$ and $L_0$ is Figiel's constant \eqref{e:figiel}.\end{lemma}

The next estimate is similar to \cite[Lemma 3.27]{ENV-Banach}. Unfortunately, the proof of the lemma given in \cite{ENV-Banach} is partially based on \cite[Lemma 3.26]{ENV-Banach}, which appears to misstate Lemma \ref{l:alber}. Thus, we supply a proof of the estimate.

\begin{lemma}\label{pythagorean lemma}
Let $\XX$ be a uniformly smooth Banach space and let $L$ be a one-dimensional linear subspace of $\XX$. If $x\in \XX$, $|x|\leq 1$, and $\dist(x,L) \le \alpha |x|$, then
\begin{align} \label{s gain est}
|x| \leq |\Pi_{L}(x)| + \frac{8\alpha}{1-4\alpha} h_\XX(51 \alpha |x|),
\end{align} where $h_\XX(t)\equiv \rho_\XX(t)/t$.
\end{lemma}

\begin{proof} For brevity, write $y=\Pi_L(x)$ and $z=\Pi_L^\perp(x)=x-y$. If $\dist(x,L)\leq \alpha|x|$, then $|y|\leq |x|$ and $|z|\leq 2\alpha|x|$ by Lemma \ref{j-proj facts}. Thus, by Lemma \ref{lemma:magic} and the fundamental theorem of calculus,
\begin{align*}
|x| - |y|  &= \int_0^1 \frac{d}{dt}|y+tz|\, dt = \int_0^1 \frac{1}{2}|y+tz|^{-1}\frac{d}{dt}|y+tz|^2\,dt  \\ &= \int_0^1 |y+tz|^{-1} \langle J(y+tz), z\rangle\, dt
\leq \frac{1}{(1-4\alpha)|x|}\int_0^1 \langle J(y+tz),z\rangle\,dt,\end{align*} where we used the rough estimate $|y+tz| \geq |y|-|z| \geq |x| - 2|z| \geq (1-4\alpha)|x|$.

Now, by Lemma \ref{j-proj facts}, $\langle J(y), z\rangle=0$. Therefore, by Alber's inequality \eqref{e:alber},
\begin{align*}
 \int_0^1 \langle J(y+tz), z \rangle\, dt  &=   \int_0^1  \langle J(y+tz) - J(y), z \rangle\, dt  \le \int_0^1  |J(y+tz)-J(y)|_{\XX^*} |z|\, dt \\
& \le 4C_0|z|\int_0^1 h_\XX(8C_0L_0|tz|)\, dt \leq 4C_0|z| h_\XX(8C_0L_0|z|),\end{align*} where $C_0= \sup_{0\leq t\leq 1}\max\{1, \sqrt{\frac{1}{2}(|y|^2 + |y+tz|^2)} \} \leq 1$ (since $|x|\leq 1$) and $1<L_0<3.18$.
Recall that $h_{\XX}(t)$ is non-decreasing (see Remark \ref{rho remark}) and $|z|\leq 2\alpha |x|$. Thus, $$\int_0^1 \langle J(y+tz), z \rangle\, dt   \leq 8\alpha|x| h_\XX(51\alpha|x|),$$ as $16L_0<50.88$.
Combining the displayed estimates yields (\ref{s gain est}).
\end{proof}

We may now improve Lemma \ref{l:graph} in uniformly smooth Banach spaces.

\begin{lemma}[cf.~{Lemma \ref{l:graph}}] \label{l:graph-p} Let $\XX$ be a uniformly smooth Banach space. Suppose that $V$ is a $\delta$-separated set with $\#V\geq 2$ and there exists a line $L$ and a number $0\leq \alpha\leq 43/1224=0.0351...$ such that $\dist(x,L)\leq \alpha\delta$ for all $x\in V.$ If $v_1,v_2\in V$, then \begin{equation}\label{e:graph-p} |\Pi_L(v_1)-\Pi_L(v_2)| \leq |v_1-v_2| \leq (1+\rho_\XX(102\alpha)) |\Pi_L(v_1)-\Pi_L(v_2)|.\end{equation}
\end{lemma}

\begin{proof} Let $V\subset\XX$ be a $\delta$-separated set with $\#V\geq 2$, let $L$ be a line in $\XX$, let $\alpha\geq 0$, and assume that $\dist(x,L)\leq \alpha \delta$ for all $x\in V$. Fix any pair of distinct points $v_1,v_2\in V$. By applying two translations and invoking the triangle inequality, we may assume that $L$ is a linear subspace of $\XX$, $0=v_1\in L$, and $v:=v_2$ satisfies $|v|\geq \delta$ and $\dist(v,L)\leq 2\alpha\delta$. Applying a dilation, we may further assume that $|v|=1$. Then $\dist(v,L) \leq 2\alpha\delta \leq 2\alpha|v|$ and by Lemma \ref{pythagorean lemma}, \begin{equation*}
 |v| \leq |\Pi(v)| + \frac{16\alpha}{1-8\alpha}h_\XX(102\alpha)= |\Pi(v)|+\frac{8}{51(1-8\alpha)} \rho_\XX(102\alpha).
\end{equation*} Recall that $\rho_\XX(t)\leq t$ in any Banach space; see \eqref{basic-rho-bounds}. Hence \begin{equation*} |\Pi(v)| \geq |v| -\frac{8}{51(1-8\alpha)} \rho_\XX(102\alpha) \geq 1 - \frac{16\alpha}{1-8\alpha} = \frac{1-24\alpha}{1-8\alpha}.\end{equation*} We now require that $$\frac{8}{51(1-8\alpha)} \leq \frac{1-24\alpha}{1-8\alpha},$$ or equivalently, $\alpha\leq 43/1224$. Then combining the displayed equations yields the right hand side of \eqref{e:graph-p}. The left hand side of \eqref{e:graph-p} follows immediately from Lemma \ref{j-proj facts}.
\end{proof}

\begin{corollary}[cf.~{Corollary \ref{c:var}}] \label{c:var-p} Let $\XX$ be a uniformly smooth Banach space. Suppose that $V\subset\XX$ is a $\delta$-separated set with $\#V\geq 2$ and there exists a line $L$ and a number $0\leq \alpha\leq 43/1224$ such that $\dist(x,L)\leq \alpha\delta$ for all $x\in V$. Enumerate $V=\{v_1,\dots,v_n\}$ so that $v_{i+1}$ lies to the right of $v_i$ for all $1\leq i\leq n-1$, relative to the ordering induced by the metric projection $\pi_L$ or the $J$-projection $\Pi_L$ (see Lemma \ref{l:order-p}). Then \begin{equation}\label{e:var2-p} \sum_{i=1}^{n-1} |v_{i+1}-v_i|^s \leq (1+\rho_\XX(102\alpha))^s|v_1-v_n|^s\quad\text{for all $s\in[1,\infty)$.}\end{equation}
\end{corollary}

\begin{proof} Repeat the proof of Corollary \ref{c:var} \emph{mutatis mutandis}. Use Lemma \ref{l:graph-p} instead of Lemma \ref{l:graph}.\end{proof}

\begin{theorem}[Analyst's Traveling Salesman parameterizations for nets in uniformly smooth Banach spaces, cf.~{Corollary \ref{c:param1}}] \label{c:param2} Let $\XX$ be a uniformly smooth Banach space and let $\rho_\XX$ denote its modulus of smoothness (see Definition \ref{def:modsmooth}). Assume that $\mathscr{V}=(V_k,\lambda_k)_{k=0}^\infty$ satisfies (V0)--(V5) in \S\ref{sec:bnv}. If the sum \begin{equation}S_{\rho_\XX}(\mathscr{V}):=\sum_{k=0}^\infty\sum_{v\in V_k} \rho_\XX(102\alpha_{k,v})\lambda_k<\infty,\end{equation} then there exists a rectifiable curve $\Gamma$ containing $\bigcup_{k=0}^\infty V_k$ such that \begin{equation}\Haus^1(\Gamma)\lesssim_{C^*,\xi_1,\xi_2} r_0(1+S_{\rho_\XX}(\mathscr{V})).\end{equation}\end{theorem}

\begin{proof} As in the proof of Corollary \ref{c:param1}, set $\alpha_0=\alpha_1$, which only depends on $C^*$, $\xi_1$, and $\xi_2$. Note that $\alpha_1\leq 1/384 < 43/1224$. Thus, by Corollary \ref{c:var-p} with $s=1$ and the bound $\rho_\XX(102\alpha_1) \geq \rho_{\ell_2}(102\alpha_1)\simeq \alpha_1^2\simeq_{C^*,\xi_1,\xi_2}1$ (see \eqref{basic-rho-bounds}),\begin{align*}
S^1_\mathscr{V} &= \sum_{k=0}^\infty \sum_{(v,v')\in\Flat(k)}\tau_1(k,v,v')\lambda_k+\sum_{k=0}^\infty \sum_{\stackrel{v\in V_k}{\alpha_{k,v}\geq \alpha_1}}\lambda_k\\
                &\leq \sum_{k=0}\sum_{(v,v')\in\Flat(k)}\rho_\XX(102\alpha_{k,v})\lambda_k+\rho_{\ell_2}(102\alpha_1)^{-1} \sum_{\stackrel{v\in V_k}{\alpha_{k,v}\geq \alpha_1}} \rho_{\XX}(102\alpha_{k,v})\lambda_k \\
               &\lesssim_{C^*,\xi_1,\xi_2} S_{\rho_\XX}(\mathscr{V})<\infty.\end{align*}
By Theorem \ref{t:param}, there exists a Lipschitz map $f:[0,1]\rightarrow\XX$ such that $\Gamma:=f([0,1])$ contains $\bigcup_{k=0}^\infty V_k$ and $\Haus^1(\Gamma)\leq \Lip(f)\lesssim_{C^*,\xi_1,\xi_2} r_0(1+S^1_\mathscr{V})\lesssim_{C^*,\xi_1,\xi_2} r_0(1+S_{\rho_\XX}(\mathscr{V}))$.\end{proof}

\begin{theorem}[sufficient half of Schul's theorem in uniformly smooth Banach spaces] \label{t:smooth} Let $\XX$ be a uniformly smooth Banach space of power type $p\in(1,2]$. If $E\subset\XX$ and $S_{E,p}(\mathscr{G})=\diam E + \sum_{Q\in\mathscr{G}}\beta_{E}(Q)^p\diam Q<\infty$ for some multiresolution family $\mathscr{G}$ for $E$ with inflation factor $A_\mathscr{G}\geq 240$, then $E$ is contained in a rectifiable curve $\Gamma$ in $\XX$ with \begin{equation}\label{e:suff-smooth}\Haus^1(\Gamma) \lesssim_{\rho_\XX,A_\mathscr{G}} S_{E,p}(\mathscr{G}).\end{equation}\end{theorem}

\begin{proof} Repeat the proof of Theorem \ref{t:banach-suff} \emph{mutatis mutandis}. Use Theorem \ref{c:param2} in lieu of Corollary \ref{c:param1}.\end{proof}

Because the Banach space $\ell_p$ is uniformly smooth of power type $\min(p,2)$ when $1<p<\infty$, the sufficient condition \eqref{e:suff-p}  in Theorem \ref{t:main} follows immediately from Theorem \ref{t:smooth}.

\section{Modulus of convexity and proof of the necessary conditions}\label{sec:convex}

\subsection{Canonical parameterization of finite length continua and beta numbers associated to a parameterization}\label{sec:wazewski}
At the heart of the proof of necessary conditions in Analyst's Traveling Salesman Theorems is the existence of parameterizations of finite length continuum by Lipschitz curves. An excellent source for the essential background is the recent paper \cite{AO-curves} by Alberti and Ottolini.

Given a continuous map $f:[0,1]\rightarrow \XX$ into a metric space and a closed, non-degenerate interval $I\subset[0,1]$, the \emph{variation} $\var(f,I)$ of $f$ over $I$ is defined by \begin{equation} \var(f,I):= \sup_{a_1\leq \dots \leq a_{n+1}} \sum_{i=1}^{n}|f(a_{i+1})-f(a_i)|\in[0,\infty],\end{equation} where the supremum ranges over all finite increasing sequences in $I$. Associated to $f$, define the \emph{multiplicity function} $m(f,\cdot):X\rightarrow[0,\infty]$, $m(f,x)=\#f^{-1}(x)$ for all $x\in \XX$. The following proposition records the well-known connection between the variation of $f$ (intrinsic length) and the Hausdorff measure of the image of $f$ (extrinsic length).

\begin{proposition}[{\cite[Proposition 3.5]{AO-curves}}] Let $f:[0,1]\rightarrow \XX$ be a continuous map into a metric space and let $I\subset[0,1]$ be a closed, non-degenerate interval. If $\var(f,I)<\infty$, then $m(f|_I,\cdot)$ is a Borel function and \begin{equation}\label{density def}\var(f,I)= \int_{\XX} m(f|_I,x)\,d\Haus^1(x).\end{equation}\end{proposition}

A theorem of Wa\.{z}ewski \cite{Wazewski} asserts that every connected, compact metric space $\Sigma$ with finite one-dimensional Hausdorff measure $\Haus^1$ admits a Lipschitz parameterization by the interval $[0,1]$ with Lipschitz constant $\Lip(f)=\sup_{x\neq y} |f(x)-f(y)|/|x-y|$ at most $2\Haus^1(\Sigma)$.\footnote{This fails dramatically for higher-dimensional curves, see e.g.~the ``Cantor ladders" in \cite[\S9.2]{BNV}.} Alberti and Ottolini have recently proved the following refinement of Wa\.{z}ewski's theorem (in particular, that $f$ has degree zero). Property (2) says that the parameterization $f$ of $\Sigma$ is \emph{essentially 2-to-1}.

\begin{theorem}[{\cite[Theorem 4.4]{AO-curves}}] \label{t:ao} Let $\Sigma$ be a connected, compact metric space with $\Haus^1(\Sigma)<\infty$. Then there exists a continuous function $f:[0,1]\rightarrow \Sigma$ such that \begin{enumerate}
\item $f$ is closed, Lipschitz, surjective, and has degree zero (see \cite{AO-curves});

\item $m(f,x)=2$ for $\Haus^1$-a.e.~$x\in \Sigma$, and $\var(f,[0,1])=2\Haus^1(\Sigma)$; and,

\item $f$ has constant speed equal to $2\Haus^1(\Sigma)$.
\end{enumerate}\end{theorem}

In any Banach space $\XX$, a connected set $\Sigma\subset\XX$ has the property that $\Haus^1(\Sigma)=\Haus^1(\overline{\Sigma})$, where $\overline{\Sigma}$ denotes the closure of $\Sigma$ in $\XX$. Moreover, if $\Sigma\subset\XX$ is closed, connected, and $\Haus^1(\Sigma)<\infty$, then $\Sigma$ is compact. The proofs of these facts are simple exercises with the definitions, using convexity of $\XX$; see e.g.~ \cite[\S5]{Schul-Hilbert} (although stated there for $\XX=\ell_2$, the proofs there hold in any Banach space).

\begin{corollary}\label{c:ao} Let $\XX$ be a Banach space. If $\Sigma\subset\XX$ is connected and $\Haus^1(\Sigma)<\infty$, then there exists an essentially 2-to-1  Lipschitz surjection $f:[0,1]\rightarrow \overline{\Sigma}$ with $\Lip(f)= 2\Haus^1(\Sigma)$.\end{corollary}

For the remainder of \S3, we fix a connected set $\Sigma$ in a Banach space $\XX$ with $\Haus^1(\Sigma)<\infty$ and we fix a parameterization $f:[0,1]\rightarrow\overline{\Sigma}$ given by Corollary \ref{c:ao}. Following \cite{Ok-TST,Schul-Hilbert}, we refer to subcurves of $f$ as ``arcs''; note that we do not require arcs be 1-to-1.

\begin{definition}[arcs and associated quantities] An \emph{arc}, $\tau=f|_{[a,b]}$, of $\overline{\Sigma}$ is the restriction of $f$ to some interval $[a,b]\subset[0,1]$. Given an arc $\tau:[a,b]\rightarrow\overline{\Sigma}$, define \begin{itemize}
\item $\Domain(\tau)=[a,b]$, $\Image(\tau)=\tau([a,b])=f([a,b])$, $\Diam(\tau)=\diam\Image(\tau)$,
\item $\Start(\tau)=\tau(a)=f(a)$, $\End(\tau)=\tau(b)=f(b)$, and
\item $\Edge(\tau)=[f(a),f(b)]$, i.e~ $\Edge(\tau)$ is the line segment in $\XX$ from $f(a)$ to $f(b)$.\end{itemize}\end{definition}

\begin{remark}\label{r:Diam-is-continuous}On any interval $[a,b]\subset[0,1]$, the map $t\mapsto \Diam(f|_{[a,a+t]})$ is a continuous function on $[0,b-a]$. Thus, by the intermediate value theorem, any arc $\tau$ with $\Diam(\tau)>\alpha$ can be partitioned into a finite number of arcs $\sigma$, all of which save one have $\Diam(\sigma)=\alpha$, and the final of which has $\Diam(\sigma)$ between $\alpha$ and $2\alpha$. \end{remark}

\begin{definition}[\cite{Ok-TST,Schul-Hilbert}] Given an arc $\tau$ of $\overline{\Sigma}$, we define the \emph{arc beta number} \begin{equation}\label{beta-tilde}\tilde \beta (\tau) := \sup_{x\in \Image(\tau)} \frac{\dist\big(x,\Edge(\tau)\big)}{\Diam(\tau)}\in[0,1].\end{equation}\end{definition}

\begin{example}\label{not-monotone} Unlike the Jones' beta numbers, which satisfy \eqref{beta-monotone}, the arc beta numbers are highly non-monotone. To see this, we construct a family of simple examples in the Euclidean plane $\RR^2$. Consider the three collinear points $$x=(0,0),\quad y=(\epsilon,0),\quad z=(1,0)\quad\text{for some }\epsilon\in(0,1).$$ Let $\tau'$ be an arc whose image traces a piecewise linear path $$x\rightarrow y\rightarrow z\rightarrow y.$$ Let $\tau$ be an arc with $\Domain(\tau')\subset\Domain(\tau)$ such that the image of $\tau$ traces a piecewise linear path $$z\rightarrow \underbrace{x\rightarrow y\rightarrow z\rightarrow y}_{\Image(\tau')}.$$ We have $\Image(\tau')=\Image(\tau)=[0,1]\times\{0\}$, whence $\Diam(\tau')=\Diam(\tau)=1$. However, $\Edge(\tau')=[x,y]$, while $\Edge(\tau)=[y,z]$, so that $$\tilde\beta(\tau')=\dist(z,[x,y])=1-\epsilon\quad\text{and}\quad \tilde\beta(\tau)=\dist(x,[y,z])=\epsilon.$$ Thus, $\tilde\beta(\tau')/\tilde\beta(\tau)=(1-\epsilon)/\epsilon$ can be arbitrarily large or arbitrarily small by choosing the parameter $\epsilon$ sufficiently near $0$ or $1$, respectively. Similar examples can be made with 1-to-1 arcs by allowing the trace of $\tau'$ and $\tau$ to lie inside $[0,1]\times[-h,h]$ with $h\ll \epsilon$.\end{example}

\begin{remark}\label{the-mistake} In \cite{Schul-Hilbert} (at the very end of the statement of Lemma 3.13), the author mistakenly asserts that $\Domain(\tau')\subset\Domain(\tau)$ and $\Diam(\tau)\lesssim \Diam(\tau')$ imply that $\tilde\beta(\tau)\gtrsim\tilde\beta(\tau')$. This claim is then applied in the proofs of Lemmas 3.14 and 3.16 of that paper. By modifying the definition of almost flat arcs (see Remark \ref{the-fix} below), one can avoid this trap entirely.\end{remark}

Under special circumstances, the arc beta numbers are almost monotone. We do not use Lemma \ref{almost-monotone-I} or \ref{almost-monotone-II} below, but include them for completeness.

\begin{lemma}[almost monotonicity I] \label{almost-monotone-I} If $\tau'$ and $\tau$ are arcs in $\overline{\Sigma}$ with $\Domain(\tau')\subset\Domain(\tau)$, then \begin{equation}\tilde\beta(\tau')\Diam(\tau')\leq \tilde\beta(\tau)\Diam(\tau)+\excess(\Edge(\tau'),\Edge(\tau)),\end{equation} where $\excess(A,B)=\sup_{x\in A}\inf_{y\in B}|x-y|$ denotes the excess of $A$ over $B$.

In particular, if $\epsilon\in(0,1)$ and $\tau$ is obtained by concatenating an arc of (image) diameter no more than $\epsilon\tilde\beta(\tau')\Diam(\tau')$ to each endpoint of $\tau'$, then \begin{equation}\tilde\beta(\tau')\leq (1-\epsilon)^{-1}\left(\frac{\Diam(\tau)}{\Diam(\tau')}\right)\tilde\beta(\tau).\end{equation} \end{lemma}

\begin{proof} Let $z\in\Image(\tau')$ and $w\in\Edge(\tau')$ be any pair of points such that $\tilde\beta(\tau')\Diam(\tau')=\dist(z,\Edge(\tau'))=|z-w|$. By definition of the excess and compactness of line segments, there exists a point $v\in \Edge(\tau)$ such that $|v-w|\leq \excess(\Edge(\tau'),\Edge(\tau))=:E$. Hence $\tilde\beta(\tau)\Diam(\tau)\geq |z-v| \geq |z-w|-|w-v| \geq \tilde\beta(\tau')\Diam(\tau')-E.$

Let $\epsilon\in(0,1)$ and suppose that $\tau$ is obtained from $\tau'$ by concatenating arcs of diameter no more than $\eta:=\epsilon\tilde\beta(\tau')\Diam(\tau')$ to the endpoints of $\tau'$. Then $|\Start(\tau)-\Start(\tau')|$ and $|\End(\tau)-\End(\tau')|$ are bounded by $\eta$. Given $p\in\Edge(\tau')$, there is $0\leq t\leq 1$ such that $p=t\Start(\tau')+(1-t)\End(\tau')$. Then $q=t\Start(\tau)+(1-t)\End(\tau)\in\Edge(\tau)$ and $$\dist(p,\Edge(\tau)) \leq |p-q| \leq t|\Start(\tau)-\Start(\tau')|+(1-t)|\End(\tau)-\End(\tau')|\leq \eta.$$ As $p$ was an arbitrary point in $\Edge(\tau')$, it follows that $E\leq \eta$. Thus, \begin{equation*}(1-\epsilon)\tilde\beta(\tau')\Diam(\tau') \leq \tilde\beta(\tau')\Diam(\tau')-E\leq \tilde\beta(\tau)\Diam(\tau). \qedhere\end{equation*}
\end{proof}

\begin{lemma}[almost monotonicity II] \label{almost-monotone-II}Let $\tau'$ and $\tau$ be arcs of $\overline{\Sigma}$ with $\Domain(\tau')\subset\Domain(\tau)$, let $x'=\Start(\tau')$, let $y'=\End(\tau')$, and let $z'\in\Image(\tau')$ be a point such that \begin{equation}\tilde\beta(\tau')\Diam(\tau')\leq \lambda \dist(z',\Edge(\tau')).\end{equation} Let $x,y,z$ denote (possibly non-unique) points in $\Edge(\tau)$, which realize the  distance of $x',y',z'$ to $\Edge(\tau)$, respectively. If $z$ lies between $x$ and $y$ in $\Edge(\tau)$, then \begin{equation}\label{arc-beta-monotone} \tilde\beta(\tau')\leq 2\lambda \left(\frac{\Diam(\tau)}{\Diam(\tau')}\right)\tilde\beta(\tau).\end{equation}\end{lemma}

\begin{proof} Assign $\delta:=\max_{w'\in\{x',y',z'\}}\dist(w',\Edge(\tau))$. If $z\in[x,y]$, then we can write $z=tx+(1-t)y$ for some $t\in[0,1]$. Then $$\dist(z,\Edge(\tau'))\leq |z-(tx'+(1-t)y')|\leq t|x-x'|+(1-t)|y-y'|\leq \delta.$$ Hence \begin{equation*}\begin{split}\lambda^{-1}\beta(\tau')\Diam(\tau')&\leq\dist(z',\Edge(\tau'))\\ &\leq |z'-z|+\dist(z,\Edge(\tau'))\leq 2\delta\leq 2\tilde\beta(\tau)\Diam(\tau),\end{split}\end{equation*} where the final inequality holds because $x',y',z'\in\Image(\tau)$.\end{proof}

Let us now recall a key element in the proof of the necessary conditions in Theorems \ref{t:Jones} and \ref{t:schul}, first introduced by Okikiolu and later formalized by Schul.

\begin{definition}[\cite{Ok-TST,Schul-Hilbert}] \label{def:filter} A \emph{filtration} $\mathscr{F}=\bigcup_{n=n_0}^\infty\mathscr{F}_n$ is a family of arcs in $\overline{\Sigma}$ with the following properties. \begin{enumerate}
\item \emph{Tree structure:} If $\tau'\in \mathscr{F}_{n+1}$, then there exists a unique arc $\tau\in \mathscr{F}_n$ such that $\Domain(\tau')\subset\Domain(\tau)$.
\item \emph{Geometric diameters:} For every $\tau\in \mathscr{F}_{n}$, $\underline{A}\rho^{-n}\le \Diam(\tau)\le A\, \rho^{-n}$ for some constants $\rho>1$ and $0<\underline{A}<A<\infty$ independent of $\tau$.
\item \emph{Trivial overlaps:} For all $\tau,\tau'\in\mathscr{F}_n$, either $\tau=\tau'$, or $\Domain(\tau)$ and $\Domain(\tau')$ intersect in at most one point.
\item \emph{Partitioning:} $\bigcup_{\tau\in\mathscr{F}_n}\Domain(\tau)=\bigcup_{\tau\in\mathscr{F}_{n_0}}\Domain(\tau)$ for every $n\geq n_0$.
\end{enumerate}
\end{definition}

\begin{lemma}[\cite{Ok-TST,Schul-Hilbert}] \label{hilbert-filter} Suppose that $\XX$ is a Hilbert space. If $\mathscr{F}=\bigcup_{n=n_0}^\infty\mathscr{F}_{n}$ is a filtration, then
\begin{equation}\label{e:filter}
\sum_{\tau \in \mathscr{F}} \tilde\beta(\tau)^2\Diam(\tau) \lesssim_{(A/\underline{A}),\rho} \Haus^1\left(\textstyle\bigcup_{\tau\in\mathscr{F}_{n_0}}\Image(\tau)\right).
\end{equation}
\end{lemma}

The exponent $2$ in \eqref{e:filter} is a consequence of the Pythagorean theorem or parallelogram law in Hilbert space. With Lemma \ref{hilbert-filter} in hand, the remainder of the proof of necessary conditions in the Analyst's Traveling Salesman Theorem in $\RR^n$ or $\ell_2$ is essentially metric, with a strong harmonic analysis flavor. We outline these last steps in \S\ref{ss:martingale}.

\subsection{Okikiolu's filtration lemma in uniformly convex spaces}\label{ss:convex}

We now develop an analogue of Lemma \ref{hilbert-filter} in uniformly convex spaces by following the proof in $\ell_2$ from \cite{Schul-Hilbert} (which is based on \cite{Ok-TST}) and replacing the parallelogram law in Hilbert space with a suitable inequality in uniformly convex spaces from \cite{DS-metric}.

\begin{definition} Let $\XX$ be a Banach space. The \textit{modulus of convexity $\delta_\XX$} of $\XX$ is the function $\delta_\XX:[0, 2] \rightarrow [0, 1]$ defined by
\begin{equation}\delta_\XX(\epsilon) := \inf\left\{ 1 - \frac{|x+y |}{2}:|x|=|y|=1\text{ and }|x-y|\geq \epsilon\right\}.\end{equation}
\end{definition}

\begin{definition}
A Banach space is called \textit{uniformly convex} if $\delta(\epsilon) > 0$ for all $\epsilon\in(0,2]$. In this case, we say that $\XX$ is \emph{convexity power type} $p\in[2,\infty)$ if there exists $c>0$ such that $\delta(\epsilon) \geq c \epsilon^p$ for all $\epsilon\in(0,2]$.
\end{definition}

\begin{remark}[essential facts]\label{delta remark} For general background on the modulus of convexity and uniformly convex spaces, we again refer the reader to \cite[Chapter Three]{Diestel} or \cite[Chapter 1, \S e]{LT-II}. On any Banach space $\XX$, the modulus of convexity satisfies the inequality \begin{equation}\label{best-delta} \delta_\XX(\epsilon) \leq \delta_{\ell_2}(\epsilon)= 1 - \sqrt{1-\epsilon^2/4}.\end{equation} Moreover, for all $\epsilon\in(0,2]$, \begin{equation}\label{e:delta-2}\delta_\XX(\epsilon) = \inf\left\{ 1 - \frac{|x+y |}{2}:|x|\leq 1,\ |y|\leq 1,\text{ and }|x-y|\geq \epsilon\right\};\end{equation} see e.g.~ \cite[p.~60]{LT-II}. In addition, \begin{equation} \label{delta-frac-increase} \frac{\delta_\XX(\epsilon_1)}{\epsilon_1} \leq \frac{\delta_\XX(\epsilon_2)}{\epsilon_2}\quad\text{for all }0<\epsilon_1\leq \epsilon_2\leq 2;\end{equation} see e.g.~\cite[Lemma 1.e.8]{LT-II}. In contrast with the modulus of smoothness $\rho_\XX$ (see \S\ref{ss:smooth}), the modulus of convexity $\delta_\XX$ is not necessarily a convex function.
The dual $\XX^*$ of a uniformly smooth Banach space $\XX$ is uniformly convex and every uniformly convex Banach space is reflexive.
\end{remark}

\begin{example}\label{example:convex} By Hanner's inequalities \cite{Hanner}, $\delta_{L^p}(\epsilon) = \frac18(p-1)\epsilon^2+o(\epsilon^2)$ when $1<p\leq 2$; and $\delta_{L^p}(\epsilon)=p^{-1}2^{-p}\epsilon^p+o(\epsilon^p)$ when $2\leq p<\infty$. In particular, the $L^p$ spaces are uniformly convex with power type $\max(2,p)$ when $1<p<\infty$.\end{example}

In \cite{DS-metric}, David and Schul observed that the modulus of convexity on a uniformly smooth space of power type $p$ can be used to control the triangle inequality excess from below. Because of its centrality to the proof of \eqref{e:nec-2}, we include a proof of their estimate here for reference. Actually, we provide a slightly stronger statement. There is a large literature on related substitutes for the Pythagorean theorem and parallelogram law in Banach spaces; see e.g.~ \cite{BynumDrew72,Bynum76,ChengRoss15,Cheng-optimal}.

\begin{lemma}[see {\cite[Lemma 8.2]{DS-metric}}]\label{mod of convexity bound} Suppose  $\XX$ is a uniformly convex Banach space. If $x,y,z\in\XX$, then
\begin{equation}\label{e:ds-bound}
|x-y| + |y-z| - |x-z| \geq 2r\,\delta_\XX\left(\frac{\dist(y,[x,z])}{r}\right)\quad\text{for all }r\geq\max\{|x-y|,|y-z|\}.
\end{equation}
\end{lemma}

\begin{proof} By \eqref{best-delta}, $\delta_\XX(1)\leq 1-(\sqrt{3}/2)<1$. Thus, in the degenerate case $x=z$, we have $|x-y|+|y-z|-|x-z| = 2|x-y| > 2\delta_\XX(1)= 2\delta_\XX(\dist(y,[x,z])/\max\{|x-y|,|y-z|\})$. This establishes \eqref{e:ds-bound} with $r=\max\{|x-y|,|y-z|\}$; the general case follows from  \eqref{delta-frac-increase}. For the remainder of the proof, we assume that $x\neq z$.

If $y\in [x,z]$, then \eqref{e:ds-bound} holds trivially (both sides vanish). Thus, we may also assume $\dist(y,[x,z])>0$.

Because $x\neq z$, the function $g(w) = |w-x|/|x - z|$ is continuous along $[x, z]$ with $g(x)=0$ and $g(z)=1$. Hence there exists $y_0 \in [x, z]$ such that
\begin{equation}\label{ratio1}
\frac{|y_0 - x|}{|x - z|} = \frac{|y- x|}{|x - y| + |y - z|}.
\end{equation}  Because $y_0\in[x,z]$, we have $|x-y_0|+|y_0-z|=|x-z|$ and it follows that
\begin{equation}\label{ratio2}
\frac{|y_0 - z|}{|x - z|} = \frac{|y- z|}{|x - y| + |y - z|}.
\end{equation} Rearranging \eqref{ratio1}, we see that
\begin{align*}
|y_0 - x| = |y - x| \frac{|x- z|}{|x - y| + |y - z|} \le |y - x|
\end{align*} by the triangle inequality. By a parallel argument, starting from \eqref{ratio2},  $|y_0 - z| \le |y - z|$.  Therefore, $y, y_0 \in B(x,|y-x|)\cap B(z,|y-z|)$.

Now, let $y' = (y+ y_0)/2$ and $h = |y - y_0| = 2 |y - y'|.$  Invoking \eqref{e:delta-2} on a scaled and translated copy of $B(x,|y-x|)$ and similarly on $B(z,|y-z|)$, we obtain
\begin{equation*}
\delta_\XX\left(\frac{h}{|x-y|}\right)  \le 1 - \frac{|y' - x|}{|y-x|}\quad\text{and}\quad \delta_\XX\left(\frac{h}{|z-y|}\right)  \le 1 - \frac{|y' - z|}{|y-z|}.
\end{equation*}
Therefore, if $|x-y|\leq r$ and $|y-z|\leq r$,
\begin{align*}
|x-y| + |y-z| - |x-z| & \ge |x-y| + |y-z| - |x-y'| - |z - y'|\\
& \ge |x - y|\, \delta_\XX\left(\frac{h}{|x - y|}\right) + |z- y|\, \delta_\XX\left(\frac{h}{|z - y|}\right) \\
& \ge 2 r \delta_\XX(h/r)
\end{align*} by \eqref{delta-frac-increase}.
Since $y_0 \in [x,z]$, we have $\dist(y,[x,z])\leq |y-y_0|=h$ and \eqref{e:ds-bound} follows.
\end{proof}

\begin{lemma}[filtration lemma, {cf.~Lemma \ref{hilbert-filter}}]\label{banach-filter}
Let $\XX$ be a uniformly convex Banach space of power type $p\in[2,\infty)$, say $\delta_\XX(\epsilon)\geq c\epsilon^p$ for all $\epsilon\in(0,2]$. If $\mathscr{F}=\bigcup_{n=n_0}^\infty \mathscr{F}_n$ is a filtration in the sense of Definition \ref{def:filter}, then
\begin{equation} \begin{split} \label{e:filter-b}
\sum_{\tau \in \mathscr{F}} \tilde\beta(\tau)^p\Diam(\tau) &\lesssim_{c,p,(A/\underline{A}),\rho} \sum_{\tau\in\mathscr{F}_{n_0}}\var(f,\Domain(\tau))-\sum_{\tau\in\mathscr{F}_{n_0}}\Haus^1(\Edge(\tau))\\ &\lesssim_{c,p,(A/\underline{A}),\rho}\Haus^1\left(\textstyle\bigcup_{\tau\in\mathscr{F}_{n_0}}\Image(\tau)\right). \end{split}
\end{equation} The implicit constant is of the form $c^{-1}(A/\underline{A})^{p-1}C(\rho,p)$ and blows up as either $c\downarrow 0$ or $p\uparrow\infty$; see \eqref{I-bound}.
\end{lemma}

\begin{proof} We mimic the proof of Lemma \ref{hilbert-filter} in \cite{Schul-Hilbert}, invoking Lemma \ref{mod of convexity bound} at a critical juncture to replace estimates depending on Hilbert space geometry.
For every $\tau \in \mathscr{F}_{n}$ and $k\in\N$, we let $\mathscr{F}_{\tau,k}$ denote the $k$th generation descendants of $\tau$,
\begin{equation}
\mathscr{F}_{\tau,k} := \{\tau' \in \mathscr{F}_{n+k} : \Domain(\tau') \subset \Domain(\tau) \}.
\end{equation} Also, for every $\tau\in\mathscr{F}$, define
\begin{equation} \label{Delta-def} \begin{split}
\Delta(\tau):= \left(\sum_{\tau' \in \mathscr{F}_{\tau, 1}} |\Start(\tau')-\End(\tau')|\right) - |\Start(\tau)-\End(\tau)|\in[0,\infty).
\end{split}\end{equation} We immediately see that for each $n_1\geq n_0$,
\begin{equation*}\begin{split}\sum_{\tau\in\mathscr{F}_{n_0}} |\Start(\tau)-\End(\tau)|+\sum_{n=n_0}^{n_1} \sum_{\tau\in\mathscr{F}_n} \Delta(\tau) &=\sum_{\tau'\in\mathscr{F}_{n_1+1}} |\Start(\tau')-\End(\tau')| \\&\leq \sum_{\tau\in\mathscr{F}_{n_0}}\var(f,\Domain(\tau)).
\end{split}\end{equation*} Thus, because $|\Start(\tau)-\End(\tau)|=\Haus^1(\Edge(\tau))$ for each arc and $f$ is essentially 2-to-1, recalling (\ref{density def}), we obtain, \begin{equation}\begin{split}\label{Delta-sum} \sum_{\tau\in\mathscr{F}} \Delta(\tau) &\leq \sum_{\tau\in\mathscr{F}_{n_0}} \var(f,\Domain(\tau))-\sum_{\tau\in\mathscr{F}_{n_0}}\Haus^1(\Edge(\tau))\\
 &\leq \sum_{\tau\in\mathscr{F}_{n_0}} \int_\XX m(\tau,x)\,d\Haus^1(x) \leq 2\Haus^1\left(\textstyle\bigcup_{\tau\in\mathscr{F}_{n_0}}\Image(\tau)\right). \end{split}\end{equation}

Next, for every $\tau \in \mathscr{F}$, we define the \emph{discretized edge distance} $d_\tau$ by
\begin{equation}
d_{\tau}:= \sup_{\tau' \in \mathscr{F}_{\tau, 1}}\sup_{x\in\Edge(\tau')} \dist(x, \Edge(\tau)).
\end{equation}In addition, for every arc $\tau$ and $k\in\N$, choose an arc $\tau_k \in \mathscr{F}_{\tau, k}$ such that $d_{\tau_k}$ is maximal among all arcs in $\mathscr{F}_{\tau, k}.$ We also write $\tau_0\equiv \tau$. We claim that \begin{equation}\label{d-sum} \tilde\beta(\tau)\Diam(\tau) \leq \sum_{k=0}^\infty d_{\tau_k}.\end{equation} To see this, choose an auxiliary sequence $\tau^k\in\mathscr{F}_{\tau,k}$ inductively so that $\tau^0\equiv \tau$ and $\tilde\beta(\tau^{k+1})\Diam(\tau^{k+1})$ is maximal over all arcs in $\mathscr{F}_{\tau^k,1}$. Then
\begin{align*}
\tilde\beta(\tau)\Diam(\tau) = \sum_{k=0}^{\infty}\left( \tilde\beta(\tau^k)\Diam(\tau^k) - \tilde\beta(\tau^{k+1})\Diam(\tau^{k+1})\right),
\end{align*} because the series is telescoping and absolutely convergent by our assumption that the arcs have geometrically decaying diameters. Moreover, $$ \tilde \beta(\tau^k)\Diam(\tau^k) - \tilde\beta(\tau^{k+1})\Diam(\tau^{k+1})\leq d_{\tau^k}$$ by the triangle inequality and definition of the sequence $\tau^k$. Hence
\begin{align*}
\tilde\beta(\tau)\Diam(\tau) \le \sum_{k=0}^{\infty} d_{\tau^k} \leq \sum_{k=0}^\infty d_{\tau_k}
\end{align*} by maximality of the distance $d_{\tau_k}$ among all arcs in $\mathscr{F}_{\tau,k}$. This verifies \eqref{d-sum}.

The proof up until this point is valid in any Banach space. To continue, we now suppose that $\delta_\XX$ is convexity power type $p\in[2,\infty)$, say $\delta_\XX(\epsilon)\geq c\epsilon^p$ for all $\epsilon\in(0,2]$. By Lemma \ref{mod of convexity bound} and the triangle inequality,
\begin{equation} \label{d-triangle-excess} 2c\, \frac{d_\tau^p}{\Diam(\tau)^{p-1}} \leq \left(\sum_{\tau' \in \mathscr{F}_{\tau, 1}} |\Start(\tau')-\End(\tau')|\right) - |\Start(\tau)-\End(\tau)| =\Delta(\tau) \end{equation} for any arc $\tau\in\mathscr{F}$. Now, by \eqref{d-sum} and Minkowski's inequality for $\ell_p$,
\begin{align*}
\underbrace{\left(\sum_{\tau\in\mathscr{F}} \tilde\beta(\tau)^p \Diam(\tau)\right)^{\frac{1}{p}}}_I & \leq \left(\sum_{\tau\in\mathscr{F}}\left(\sum_{k=0}^\infty d_{\tau_k}\right)^p \Diam(\tau)^{1-p}\right)^{\frac{1}{p}}  \leq \sum_{k=0}^{\infty} \left(\sum_{\tau\in\mathscr{F}} \frac{d_{\tau_k}^p}{\Diam(\tau)^{p-1}}\right)^{1/p}.
\end{align*} If $\tau\in\mathscr{F}$, say $\tau\in\mathscr{F}_n$, and $k\in\N$, then \begin{align*}
\Diam(\tau_{k}) &\leq A\, \rho^{-(n+k)} = A\, \rho^{-k} \rho^{-n} \leq (A/\underline{A})\, \rho^{-k} \Diam(\tau).\end{align*} Hence
\begin{align*}
I \leq \sum_{k=0}^{\infty} \big((A/\underline{A})\rho^{-k}\big)^{\frac{p-1}{p}} \left(\sum_{\tau\in\mathscr{F}} \frac{(d_{\tau_k})^p}{\Diam(\tau_k)^{p-1}} \right)^{\frac{1}{p}}.\end{align*}
Thus, by \eqref{d-triangle-excess}, \begin{equation}\begin{split}\label{I-bound} I &\leq
\sum_{k=0}^{\infty} (2c)^{-\frac{1}{p}}\big((A/\underline{A})\rho^{-k}\big)^{\frac{p-1}{p}} \left(\sum_{\tau\in\mathscr{F}}  \Delta(\tau_k)\right)^{\frac{1}{p}}
\\ &\leq \sum_{k=0}^{\infty} (2c)^{-\frac{1}{p}}\big((A/\underline{A})\rho^{-k}\big)^{\frac{p-1}{p}} \left(\sum_{\tau\in\mathscr{F}}  \Delta(\tau)\right)^{\frac{1}{p}} = \frac{(2c)^{-\frac{1}{p}}(A/\underline{A})^{\frac{p-1}{p}}}{1-\rho^{\frac{1-p}{p}}} \left(\sum_{\tau\in\mathscr{F}}  \Delta(\tau)\right)^{\frac{1}{p}}.
\end{split}\end{equation} Therefore, combining \eqref{Delta-sum} and \eqref{I-bound}, we obtain \eqref{e:filter-b}.\end{proof}

\begin{remark}The proof of Lemma \ref{banach-filter} used the standing assumption that $f:[0,1]\rightarrow\overline{\Sigma}$ is essentially at most 2-to-1. Instead, if one starts with a parameterization $f$ that is essentially at most $m$-to-1 for some $m\geq 1$, then \eqref{e:filter-b} still holds, but with the implicit constant in the second inequality multiplied by a factor $m/2$.\end{remark}

\subsection{Filtration design: bounding sums of $\beta_\Sigma(Q)^p\diam Q$ from above}\label{ss:martingale}

The following decoration of a lemma by Schul is an essential tool for constructing filtrations from certain families of arcs, which we think of as \emph{prefiltrations}. Although originally stated in Hilbert space, it is clear upon reading the proof that the lemma is valid in any metric space. Unfortunately, the statement in \cite{Schul-Hilbert} contains a mistake, claiming erroneously that $\Diam(\tau)\lesssim\Diam(\tau')$ implies $\tilde\beta(\tau)\gtrsim\tilde\beta(\tau')$, which is false (even in $\RR^2$) by Example \ref{not-monotone}. Also, Schul claims that one may transform a prefiltration into $2CJ$ or fewer filtrations with $J\gtrsim 1$ when $\underline{A}=1$ and $\rho=2$, but after writing down the details it seems to us that one must break apart the prefiltration into a larger number of families depending on $A$ and take $J\gtrsim_{A} 1$ in order to verify \eqref{new-arc-diameters} and \eqref{new-arc-properties}. Thus, we supply a corrected statement and detailed proof, whose outline is due to Schul. We defer the proof to Appendix \ref{ss:prefiltration}.

\begin{lemma}[prefiltration lemma, cf.~{\cite[Lemma 3.13]{Schul-Hilbert}}]\label{build-filtrations} Let $\XX$ be a metric space and let $f:[0,1]\rightarrow\overline{\Sigma}$ be a continuous parameterization of a set $\overline{\Sigma}\subset\XX$. Assume that $\rho>1$, $0<\underline{A}<A<\infty$, and $J\geq 1$ is any integer such that $\rho^{J}>6A/\underline{A}$. Then for every family $\mathscr{F}^0=\bigcup_{n=n_0}^\infty\mathscr{F}_n^0$ of arcs in $\overline{\Sigma}$ with $\mathscr{F}^0_{n_0}\neq\emptyset$ satisfying \begin{enumerate}
\item bounded overlap: for every arc $\tau\in\mathscr{F}^0_n$, there exists no more than $C$ arcs $\tau'\in\mathscr{F}^0_n$ such that $\Domain(\tau)\cap\Domain(\tau')\neq\emptyset$ for some constant $C$ independent of $\tau$,
\item geometric diameters: for every arc $\tau\in\mathscr{F}^0_n$, we have $\underline{A}\rho^{-n}\leq \Diam(\tau)\leq A\rho^{-n}$,
\end{enumerate} we can construct $5(A/\underline{A})CJ$ or fewer filtrations $\mathscr{F}^1=\bigcup_{n=n_1}^\infty \mathscr{F}^1_n$, $\mathscr{F}^2=\bigcup_{n=n_2}^\infty \mathscr{F}^2_n$, \dots, with starting index $n_j\in \{n_0,n_0+1,\dots,n_0+J-1\}$ for all $j$ and \begin{equation}\label{new-arc-diameters}\frac{1}{4}\left(\underline{A}\rho^{(J-1)n_j}\right)\rho^{-Jn} \leq \Diam(\tau)< 2\left(A\rho^{(J-1)n_j}\right) \rho^{-Jn}\quad\text{for all $j$, $\tau\in\mathscr{F}^j_n$, $n\geq n_j$,}\end{equation} such that for every index $n\geq n_0$ and arc $\tau'\in\mathscr{F}^0_{n}$, there exists $\mathscr{F}^j$ (in the list of filtrations), an index $N$ with $n-n_j=J(N-n_j)$, and an arc $\tau\in\mathscr{F}^j_N$ such that \begin{equation}\label{new-arc-properties} \Domain(\tau')\subset \Domain(\tau)\quad\text{and}\quad \Diam(\tau) < 2\Diam(\tau').\end{equation} The assignment $(n,\tau')\mapsto (\mathscr{F}^j,N,\tau)$ is injective.
\end{lemma}

\begin{remark}\label{multiple-arcs}In Lemma \ref{build-filtrations}, we may allow an arc to appear in some $\mathscr{F}^0_n$ several times, so long as the bounded overlap condition in the hypothesis is computed with multiplicity. During preprocessing (see the proof), each instance of a repeated arc $\tau'\in\mathscr{F}^0_n$ will be assigned to a different intermediate family $\mathscr{D}^j_n$ and their extensions $\tau$ will land in different filtrations $\mathscr{F}^j_N$. We use this observation in the proofs of Lemmas \ref{l:H1} and \ref{l:H3} below.\end{remark}

Because the $\ell_p$ spaces are uniformly convex of power type $\max(2,p)$ when $1<p<\infty$, the necessary condition \eqref{e:nec-2} in Theorem \ref{t:main2} is an immediate consequence of the following theorem.

\begin{theorem}[necessary half of Schul's theorem in uniformly convex Banach spaces] \label{t:convex} Let $\XX$ be a uniformly convex Banach space of power type $p\in [2, \infty)$. If $\Sigma\subset\XX$ is connected  and $\mathscr{H}$ is a (partial) multiresolution family for $\Sigma$ with inflation factor $A_\mathscr{H}>1$, then \begin{equation}\label{e:nec-convex}S_{\Sigma,p}(\mathscr{H})=\diam \Sigma+\sum_{Q\in\mathscr{H}} \beta_\Sigma(Q)^p \diam Q\lesssim_{p,\delta_\XX,A_\mathscr{H}} \Haus^1(\Sigma).\end{equation}\end{theorem}

In the remainder of this section, we outline the proof of Theorem \ref{t:convex} in detail. Because \eqref{e:nec-convex} is trivial when $\Haus^1(\Sigma)=\infty$, we may continue to assume that $\Haus^1(\Sigma)<\infty$ and work with the essentially 2-to-1 Lipschitz parameterization $f:[0,1]\rightarrow\overline{\Sigma}$ fixed above in \S \ref{sec:wazewski}. Also, because $\Sigma$ is connected, $\diam\Sigma \leq \Haus^1(\Sigma)$. Thus, the essential task is to bound the beta number sum from above in terms of $\Haus^1(\Sigma)$. To carry this out, we modify the proof from \cite[\S3]{Schul-Hilbert} with the correction noted in Remark \ref{the-fix} below. Most of the argument works in any Banach space and we make sure to explicitly state wherever we need uniform convexity. We work with two classes of almost flat arcs and one class of non-almost-flat arcs, which we call dominant arcs. Flatness is measured with respect to a best approximating line for the image of the arc.

\begin{definition}\label{arc-types} For any ball $Q\in\mathscr{H}$ and scaling factor $\lambda\in\{1,5,7\}$, let \begin{equation}\Lambda(\lambda Q):=\left\{f|_{[a,b]}:\begin{array}{c}[a,b]\text{ is a connected component of }f^{-1}(\overline{\Sigma}\cap 2\lambda Q)\\ \text{such that }\lambda Q\cap f([a,b])\neq\emptyset\end{array}\right\}.\end{equation} The elements in $\Lambda(\lambda Q)$ are arcs in $2\lambda Q$ that touch $\lambda Q$. Agree to write $\beta_{\Lambda(\lambda Q)}(2\lambda Q)$ as shorthand for $\beta_{\bigcup\{\Image(\tau):\tau\in\Lambda(\lambda Q)\}}(2\lambda Q)$.

An arc $\tau\in\Lambda(\lambda Q)$ is called \emph{$*$-almost flat} if \begin{equation}\label{e:star-flat} \beta(\tau):= \inf_{L}\sup_{z\in\Image(\tau)} \frac{\dist(z,L)}{\Diam(\tau)} \leq 50\flatepsilon \beta_{\Lambda(\lambda Q)}(2\lambda Q),\end{equation} where $L$ ranges over all lines in $\XX$ and $0<\flatepsilon\ll 1$ is a parameter, ultimately chosen to depend on at most the inflation factor $A_\mathscr{H}$ of $\mathscr{H}$ (shortly after Lemma \ref{l:H2}). Denote the set of $*$-almost flat arcs in $\Lambda(\lambda Q)$ by $S^*(\lambda Q)$.

An arc $\tau\in\Lambda(\lambda Q)$ is called \emph{almost flat} if
\begin{equation}\label{e:almost-flat-def}  \beta(\tau)\leq \flatepsilon \beta_{\Sigma}(Q).\end{equation} Denote the set of almost flat arcs in $\Lambda(\lambda Q)$ by $S(\lambda Q)$. An arc $\tau\in\Lambda(\lambda Q)\setminus S(\lambda Q)$ is called \emph{dominant}.\end{definition}

\begin{remark}\label{the-fix} We have modified the definition of almost flat arcs from \cite{Ok-TST,Schul-Hilbert}, which instead required $\tilde\beta(\tau)\leq \epsilon_2\beta_{\Sigma}(Q)$. Note that the quantity $\beta(\tau)$ is nothing other than the Jones' beta number $\beta_{\Image(\tau)}(\Image(\tau))$ of the image of $\tau$ in its own window. By \eqref{beta-monotone}, it follows that $\Domain(\tau')\subset\Domain(\tau)$ and $\Diam(\tau)\lesssim \Diam(\tau')$ imply $\beta(\tau)\gtrsim\beta(\tau')$. Further, it is apparent that $\beta(\tau)\leq \tilde\beta(\tau)$ for every arc. By using the Jones' beta number $\beta(\tau)$ instead of the arc beta number $\tilde\beta(\tau)$ wherever possible in the proof below, we avoid the issue in Remark \ref{the-mistake}. We only work with $\tilde\beta(\tau)$ in two spots, at the end of proofs of Lemmas \ref{l:H1} and \ref{l:H3}, when we must invoke Lemma \ref{banach-filter}.\end{remark}

\begin{remark}\label{arc-choices} For every $Q\in\mathscr{H}$ and $\lambda\in\{1,5,7\}$, $S(\lambda Q)\subset S^*(\lambda Q)\subset \Lambda(\lambda Q)$ by \eqref{beta-monotone}. Domains of distinct arcs in $\Lambda(\lambda Q)$ are disjoint, although the images of distinct arcs in $\Lambda(\lambda Q)$ may coincide, because we do not require $f$ be 1-to-1. If $\tau\in\Lambda(\lambda Q)$ is the ``only arc'' in $\Lambda(\lambda Q)$ in the sense that the image of each arc in $\Lambda(\lambda Q)$ is contained in $\Image(\tau)$ and $\beta_\Sigma(Q)>0$, then $\tau$ is dominant, since $\Diam(\tau)\leq \diam 14Q$ and $\flatepsilon\ll 1$ imply \begin{equation}\beta_\Sigma(Q)=\beta_{\Image(\tau)}(Q)\leq 14\beta(\tau)<(1/\flatepsilon) \beta(\tau).\end{equation} One should think of the latter situation as being the infinitesimally generic case. Note that almost flat arcs are defined relative to $\beta_\Sigma(Q)$ and $*$-almost flat arcs are defined relative to $\beta_{\Lambda(\lambda Q)}(2\lambda Q)$. This is intentional. The choice of the scaling factors $\lambda\in\{1,5,7\}$ and definition of $*$-almost flat arcs are made to implement the proof of Lemma \ref{l:H3}. Lastly, we have defined $\Lambda(\lambda Q)$ as arcs in $2\lambda Q$ touching $\lambda Q$ so that arcs $\tau\in\Lambda(\lambda Q)$ have uniformly large diameter: $\Diam(\tau)\geq \lambda A_{\mathscr{H}}2^{-k}$ whenever $Q=B(x,A_{\mathscr{H}}2^{-k})$, $x\in X_k$, and $\Sigma\setminus 2\lambda Q\neq\emptyset$. \end{remark}

To proceed, we categorize the balls $Q$ in the multiresolution family $\mathscr{H}$ according the behavior of the associated arcs $\Lambda(Q)\cup\Lambda(5Q)\cup\Lambda(7Q)$. First, let $\mathscr{H}_0$ denote the collection of all balls $Q\in\mathscr{H}$ such that \begin{equation}\beta_\Sigma(Q)=0\quad\text{or}\quad\Sigma\subset 14Q,\end{equation} where $\lambda Q=B(x,\lambda r)$ denotes the concentric dilate of the ball $Q=B(x,r)$ by $\lambda>0$. Next, we group $\mathscr{H}\setminus\mathscr{H}_0$ into three (overlapping) subfamilies $\mathscr{A}$, $\mathscr{B}$, and $\mathscr{C}$, as follows.

Let $(X_k)_{k\in\ZZ}$ denote the sequence of $2^{-k}$-separated sets for $\Sigma$ that is used to define the (partial) multiresolution family $\mathscr{H}$. For each $k\in\ZZ$, let $\mathscr{N}_k=\{B(x,(1/3)2^{-k}):x\in X_k\}$ denote the collection of \emph{net balls} of level $k$. By the triangle inequality, \begin{equation}\label{fat-separated} \gap(B_1,B_2)\geq \frac{1}{3}2^{-k}\quad\text{ for all }B_1\neq B_2\text{ in }\mathscr{N}_k,\end{equation} where the reader may recall that $\gap(S,T)=\inf\{\dist(s,t):s\in S,t\in T\}$.  Suppose $Q=B(x,A_\mathscr{H} 2^{-k})$ for some $k\in\ZZ$ and $x\in X_k$. Let $\lambda\in\{1,5,7\}$. We say that $Q\in\mathscr{H}_1^\lambda$ if \begin{equation}\text{there exists }\tau\in \Lambda(\lambda Q)\setminus S(\lambda Q)\text{ such that }\Image(\tau)\cap B(x,(1/3)2^{-k})\neq\emptyset.\end{equation} That is, $Q\in\mathscr{H}_1^\lambda$ if there exists a dominant arc in $\Lambda(\lambda Q)$ that touches the net ball at the center of $Q$. Otherwise, $Q\not\in\mathscr{H}^\lambda_1$ and every arc in $\Lambda(\lambda Q)$ that passes through the net ball at the center of $Q$ is an almost flat arc in $\Lambda(\lambda Q)$; in this case, we assign $Q$ to either $\mathscr{H}_2^\lambda$ or $\mathscr{H}_3^\lambda$, depending on global geometry of the $*$-almost flat arcs in $\Lambda(\lambda Q)$. Following a convention from \cite{Schul-Hilbert} (see page 345), we write $\beta_{S^*(\lambda Q)}(2\lambda Q)$ as shorthand for $\beta_{\bigcup\{\Image(\tau):\tau\in S^*(\lambda Q)\}}(2\lambda Q).$ Fix a constant $0<\flatarcsepsilon\ll 1$ depending on at most $A_\mathscr{H}$ to be specified below. We say that $Q\in\mathscr{H}_2^\lambda$ if \begin{equation}Q\not\in\mathscr{H}_1^\lambda\quad\text{and}\quad \beta_{S^*(\lambda Q)}(2\lambda Q)> \flatarcsepsilon \beta_{\Lambda(\lambda Q)}(2\lambda Q).\end{equation} We say that $Q\in\mathscr{H}_3^\lambda$ if \begin{equation} Q\not\in \mathscr{H}^\lambda_1\quad\text{and}\quad \beta_{S^*(\lambda Q)}(2\lambda Q)\leq \flatarcsepsilon \beta_{\Lambda(\lambda Q)}(2\lambda Q).\hbox{\quad\ \ }\end{equation} When a ball $Q\in\mathscr{H}^\lambda_2$, $\beta_{S^*(\lambda Q)}(Q)$ dominates $\flatarcsepsilon\beta_{\Lambda(\lambda Q)}(2\lambda Q)$. By contrast, when a ball $Q\in\mathscr{H}^\lambda_3$, the $*$-almost flat arcs in $\Lambda(\lambda Q)$ are collectively much flatter inside of the window $2\lambda Q$ than the union of all of the arcs in $\Lambda(\lambda Q)$.

We now define
 \begin{align} \mathscr{A}&:=\{Q\in\mathscr{H}\setminus\mathscr{H}_0: Q\in\mathscr{H}^\lambda_1\text{ for $\lambda=1$, $\lambda=5$, or $\lambda=7$}\},\\
\mathscr{B}&:=\{Q\in\mathscr{H}\setminus\mathscr{H}_0: Q\in\mathscr{H}^\lambda_2\text{ for $\lambda=1$ or $\lambda=5$}\},\\
\mathscr{C}&:=\{Q\in\mathscr{H}\setminus\mathscr{H}_0: Q\in\mathscr{H}^\lambda_3\text{ for $\lambda=1$ and $\lambda=5$, and $Q\not\in\mathscr{H}^7_1$}\}.
\end{align} (We neither use $\mathscr{H}^7_2$ nor $\mathscr{H}^7_3$.) Note that $\mathscr{H}\setminus\mathscr{H}_0\subset\mathscr{A}\cup\mathscr{B}\cup\mathscr{C}$. While the family $\mathscr{C}$ is disjoint from $\mathscr{A}\cup\mathscr{B}$, some balls in $\mathscr{H}\setminus\mathscr{H}_0$ could belong to both $\mathscr{A}$ and $\mathscr{B}$.

\begin{remark}Our classification of balls in the family $\mathscr{H}$ is slightly different than in \cite{Schul-Hilbert}, but roughly speaking our $\mathscr{A}$, $\mathscr{B}$, and $\mathscr{C}$ balls correspond to Schul's $\mathscr{G}_1$, $\mathscr{G}_2$, and $\mathscr{G}_3$ balls. See Figure 4 on page 346 of \cite{Schul-Hilbert} for an illustration of the different families. More specifically, our class $\mathscr{B}$ corresponds to Schul's $\Delta_1$ and $\Delta_{2.1}$ balls. Schul also defines $\Delta_{2.2}$ balls in $\mathscr{G}_2$, but our definition includes these in $\mathscr{A}$. We introduced the net balls to consolidate the estimates for $\mathscr{G}_1$ and $\Delta_{2.2}$ balls.\end{remark}

With the families $\mathscr{H}_0$, $\mathscr{A}$, $\mathscr{B}$, and $\mathscr{C}$ now defined, the proof of Theorem \ref{t:convex} reduces to establishing an estimate like \eqref{e:nec-convex} for each category.

\begin{lemma}[counting]\label{l:H0} If $\XX$ is an arbitrary Banach space, then for all $s>1$, \begin{equation}\sum_{Q\in\mathscr{H}_0}\beta_\Sigma(Q)^s\diam (Q) \lesssim_{s} A_\mathscr{H}\Haus^1(\Sigma),\end{equation} where the implicit constant blows up as $s\downarrow 1$.\end{lemma}

\begin{proof} Modify the proof of \cite[Lemma 3.9]{Schul-Hilbert}. Set $D=\diam \Sigma$. Suppose that $k\in\ZZ$, $x\in X_k$, and $Q=B(x,A_\mathscr{H}2^{-k})\in\mathscr{H}$. If $\diam 14Q=28A_\mathscr{H}2^{-k}<D$, then $\Sigma\not\subset 14Q$. Hence if $Q\in\mathscr{H}_0$ and $\beta_{\Sigma}(Q)>0$, then $28A_\mathscr{H}2^{-k}\geq D$. Thus, letting $\mathscr{H}(k)$ denote all $Q\in\mathscr{H}$ of radius $A_\mathscr{H}2^{-k}$, we have $$\sum_{Q\in\mathscr{H}_0} \beta_\Sigma(Q)^s\diam Q
\leq 2A_\mathscr{H}\sum_{k=-\infty}^{k_0}\sum_{Q\in\mathscr{H}(k)} \beta_\Sigma(Q)^s\, 2^{-k},$$ where $k_0$ is the unique integer such that $$28A_\mathscr{H}2^{-(k_0+1)}<D\leq 28A_\mathscr{H}2^{-k_0}.$$ By considering any line passing through the ball's center, we obtain the trivial bound $\beta_\Sigma(Q)\leq (\diam \Sigma)/(\diam Q)=D/2A_\mathscr{H}2^{-k}$ for all $Q\in\mathscr{H}(k)$. It follows that $$\sum_{Q\in\mathscr{H}_0} \beta_\Sigma(Q)^s\diam Q \leq (2A_\mathscr{H})^{1-s}D^s\sum_{k=-\infty}^{k_0} 2^{k(s-1)}\#\mathscr{H}(k).$$ As $\#\mathscr{H}(k)=\# X_k \leq \#X_{k_0}$ for all $k\geq k_0$, it suffices to bound $\#X_{k_0}$. We picked $k_0$ so that $14A_\mathscr{H}2^{-k_0}<D$, which more than certainly implies $(1/2)2^{-k_0}<D/2$. In particular, $\Haus^1(\Sigma\cap U(x,(1/2)2^{-k_0}))\geq (1/2)2^{-k_0}$, because $x\in\Sigma$, $\Sigma$ is connected, and $\Sigma$ is not trapped inside the ball. By definition of a net, the open balls $U(x,(1/2)2^{-k_0})$ centered on points $x\in X_{k_0}$ are pairwise disjoint. Ergo, $$(1/2)2^{-k_0}\#X_{k_0}
\leq \sum_{x\in X_{k_0}} \Haus^1(\Sigma\cap U(x,(1/2)2^{-k_0})) \leq \Haus^1(\Sigma).$$ Therefore, $\#X_{k_0} \leq 2^{k_0+1}\Haus^1(\Sigma) \leq 56A_\mathscr{H}\Haus^1(\Sigma)/D$. Assembling all of the pieces, we have \begin{align*}
\sum_{Q\in\mathscr{H}_0} \beta_\Sigma(Q)^s\diam Q &\leq 56A_{\mathscr{H}}\Haus^1(\Sigma)(D/2A_{\mathscr{H}})^{s-1}\sum_{k=-\infty}^{k_0} 2^{k(s-1)}\\
&=56A_{\mathscr{H}}\Haus^1(\Sigma)(D/2A_{\mathscr{H}})^{s-1}\frac{2^{k_0(s-1)}}{1-2^{1-s}}
\leq 56A_\mathscr{H}\Haus^1(\Sigma)\frac{14^{s-1}}{1-2^{1-s}}, \end{align*} where in the final inequality we used once again that $2^{k_0}\leq 28A_\mathscr{H}/D$.
\end{proof}

\begin{lemma}[filtrations I]\label{l:H1} If $\XX$ is uniformly convex of power type $p\in[2,\infty)$, then \begin{equation} \label{A-sum} \sum_{Q\in\mathscr{A}} \beta_\Sigma(Q)^p\diam Q\lesssim_{p,\delta_\XX,A_\mathscr{H},\flatepsilon} \Haus^1(\Sigma).\end{equation} \end{lemma}

\begin{proof} Modify the proof of \cite[Lemma 3.14]{Schul-Hilbert}. We start with an auxiliary observation. Let $\XX$ be any Banach space. Suppose that $\xi$ is an arc of $\overline{\Sigma}$ such that $\Image(\xi)$ intersects $N$ net balls in $\mathscr{N}_k$. We claim that if $N$ is sufficiently large, then $\beta(\xi)\Diam(\xi)\geq (1/18)2^{-k}$. To see this, choose points $v_1,\dots,v_N\in\Image(\xi)$, one inside each net ball intersecting $\xi$. Then $v_1,\dots,v_N$ are $\delta=(1/3)2^{-k}$ separated by \eqref{fat-separated}. Suppose that $L$ is a line such that $\dist(v_i,L)<(1/6)\delta$ for all $i$. By Corollary \ref{c:var} with $s=1$, after reordering $v_1,\dots,v_N$, $$\frac{1}{3}(N-1)2^{-k}=(N-1)\delta \leq \sum_{i=1}^{N-1}|v_i-v_{i+1}|<(1+3/6)^2|v_1-v_N|\leq \frac{9}{4}\Diam(\xi),$$ whence $N<1+(27/4)\Diam(\xi)/2^{-k}$. Thus, if $N\geq 1+(27/4)\Diam(\xi)/2^{-k}$, then $\beta(\xi)\Diam(\xi)\geq \inf_L \sup_{i} \dist(v_i,L)\geq (1/6)\delta=(1/18)2^{-k}$. In particular, suppose that \begin{equation}\Diam(\xi)\leq 28A_\mathscr{H}2^{-k}\quad\text{and}\quad N\geq 190 A_\mathscr{H},\end{equation} which guarantees $N\geq 1+189A_\mathscr{H}= 1+(27/4)28A_\mathscr{H}\geq 1+(27/4)\Diam(\xi)/2^{-k}$. Then \begin{equation}\beta(\xi) \geq (1/18)2^{-k}/\Diam(\xi)\geq 1/504A_\mathscr{H}.\end{equation}

Let $Q\in\mathscr{A}$. Then there exists $k\in\ZZ$, $x\in X_k$, and $\lambda\in\{1,5,7\}$ such that $\lambda Q=B(x,\lambda A_\mathscr{H}2^{-k})\in\mathscr{H}^\lambda_1$. Thus, we may choose a dominant arc $\gamma_{Q}\in\Lambda(\lambda Q)\setminus S(\lambda Q)$ such that $\Image(\gamma_{Q})$ intersects the net ball $B(x,(1/3)2^{-k})$ and pick a point $y_{Q}\in\Domain(\gamma_{Q})$ such that $f(y_{Q})\in B(x,(1/3)2^{-k})$. Let $N_{Q}$ be the number of net balls $R\in\mathscr{N}_k$ such that $R$ intersects $\Image(\gamma_{Q})$. We know that $N_{Q}\geq 1$, but if $\XX$ is infinite-dimensional, then $N_{Q}$ could be arbitrarily large. We now define an arc $\tau'_Q$ according to one of two alternatives. If $N_Q<\lceil 190A_\mathscr{H}\rceil $, then we set $\tau'_Q=\gamma_Q$ and have $\beta(\tau'_Q)>\flatepsilon\beta_\Sigma(Q)$, since $\gamma_Q$ is dominant. Otherwise, if $N_Q\geq \lceil 190A_\mathscr{H}\rceil$, then we choose $\tau'_Q$ to be any arc such that $y_Q\in\Domain(\tau'_Q)\subset \Domain(\gamma_Q)$ and $\Image(\tau'_Q)$ touches precisely $\lceil 190A_\mathscr{H}\rceil$ of the net balls from $\mathscr{N}_k$ that intersect $\Image(\gamma_Q)$; then $\beta(\tau'_Q)\geq 1/504A_\mathscr{H}\geq (1/504A_\mathscr{H})\beta_\Sigma(Q)$ by the auxiliary observation from above. In both cases, the arc $\tau'_Q$ satisfies $y_Q\in\Domain(\tau'_Q)$, $\beta(\tau'_Q)\gtrsim_{A_\mathscr{H},\flatepsilon} \beta_\Sigma(Q)$, and $\tau'_Q$ intersects no more than $191A_\mathscr{H}$ net balls in $\mathscr{N}_k$. Include an instance of the arc $\tau'_Q$ in the set $\mathscr{F}^0_k$. (In the unlikely event that $\tau'_R=\tau'_Q$ for some $R\neq Q$, we treat $\mathscr{F}^0_k$ as a multiset; see Remark \ref{multiple-arcs}.)

We claim that the family $\mathscr{F}^0=\bigcup_{k=k_0}^\infty \mathscr{F}^0_{k}$ of arcs is a prefiltration (see Lemma \ref{build-filtrations}), where $k_0$ is the smallest integer such that $\mathscr{A}$ contains a ball of radius $A_\mathscr{H}2^{-k_0}$. To see this, first note that if $\tau'_Q\in\mathscr{F}^0_k$, then $$\frac{1}{3}2^{-k}\leq \Diam(\tau'_Q)\leq 28A_\mathscr{H} 2^{-k}.$$ The upper bound follows since $\tau'_Q$ is contained in $14Q$. The lower bound holds, because either $\tau'_Q=\gamma_Q$ and $\Diam(\gamma_Q)\geq \lambda A_\mathscr{H}2^{-k}$ by Remark \ref{arc-choices}, or $\tau'_Q$ touches at least 2 balls in $\mathscr{N}_k$ and \eqref{fat-separated} is in effect. Next, let's confirm that the arcs in $\mathscr{F}^0_k$ have uniformly bounded overlap for each $k$. Fix $k\geq k_0$ and $\tau'_Q\in\mathscr{F}^0_k$. Consider the set $$Y=\{y_R:\tau'_R\in\mathscr{F}^0_k\text{ and }\Domain(\tau'_R)\cap\Domain(\tau'_Q)\neq\emptyset\}\subset[0,1].$$ Because each $y_R$ belongs to the net ball at the center of $R$ and \eqref{fat-separated} holds, there is a bijection between $Y$ and arcs $\tau'_R\in\mathscr{F}^0_k$ such that $\tau'_R$ and $\tau'_Q$ intersect in their domains. Let $\tau'_S,\tau'_T\in\mathscr{F}^0_k$ be the unique arcs such that $y_S=\min Y$ and $y_T=\max Y$. Observe that $$Y\subset \Domain(\tau'_S)\cup\Domain(\tau'_Q)\cup\Domain(\tau'_T),$$ because the latter set is an interval containing $y_S$ and $y_T$. Hence $\#Y$ is bounded from above by the number of net balls in $\mathscr{N}_k$ intersected by $\tau'_S$, $\tau'_Q$, or $\tau'_T$. Therefore, we know $\tau'_Q$ intersects at most $\#Y\leq 3\cdot 191A_\mathscr{H}=573A_\mathscr{H}$ arcs $\tau'_R\in\mathscr{F}^0_k$.

Invoking Lemma \ref{build-filtrations} with parameters $\rho=2$, $A=28A_\mathscr{H}$, $\underline{A}=1/3$, $C=573A_\mathscr{H}$, and $J=1+\lceil \log_2(504A_\mathscr{H})\rceil$ (so that $2^J>6A/\underline{A}$), we can find $O(A_\mathscr{H}^2\log(A_\mathscr{H}))$ filtrations $$\mathscr{F}^1=\bigcup_{k=k_1}^\infty \mathscr{F}^1_k,\quad \mathscr{F}^2=\bigcup_{k=k_2}^\infty \mathscr{F}^2_k,\quad \mathscr{F}^3=\bigcup_{k=k_3}^\infty \mathscr{F}^2_k,\quad \dots$$ with starting indices $k_0\leq k_j \leq k_0+J-1$ such that $$\frac{1}{12}2^{(J-1)k_j}2^{-Jk} \leq \Diam(\tau) <56A_\mathscr{H} 2^{(J-1)k_j} 2^{-Jk}\quad\text{for all }\tau\in\mathscr{F}^j_k.$$ Note that the ratio of the upper and lower bounds for $\Diam(\tau)$ depends only on $A_\mathscr{H}$. Moreover, for all $\tau'_Q\in\mathscr{F}^0_k$, there exists $\mathscr{F}^j$ (in the list of filtrations), an index $K$ with $k-k_j=J(K-k_j)$, and an arc $\tau_Q\in\mathscr{F}^j_K$ such that $$\Domain(\tau'_Q)\subset\Domain(\tau)\quad\text{and}\quad\Diam(\tau_Q)<2\Diam(\tau'_Q).$$ The assignment $(k,\tau'_Q)\mapsto (\mathscr{F}^j,K,\tau_Q)$ is injective. Now, for any $\tau'_Q\in\mathscr{F}^0_k$, we have $$\beta_\Sigma(Q)\lesssim_{A_\mathscr{H},\flatepsilon} \beta(\tau'_Q)\leq 2\beta(\tau_Q)\leq 2\tilde\beta(\tau_Q),$$ where the first inequality holds by our choice of $\tau'_Q$, the second inequality holds by \eqref{beta-monotone}, because $\Diam(\tau)\leq 2\Diam(\tau')$, and the third inequality holds generally. Thus, \begin{align*}\sum_{Q\in\mathscr{A}} \beta_\Sigma(Q)^p\diam Q &\lesssim_{p,A_\mathscr{H},\flatepsilon} \sum_{\tau'\in\mathscr{F}^0}\beta(\tau')^p\Diam(\tau')
\\ &\lesssim_p \sum_j\sum_{\tau\in\mathscr{F}^j} \beta(\tau)^p\Diam(\tau)\leq \sum_j\sum_{\tau\in\mathscr{F}^j} \tilde\beta(\tau)^p\Diam(\tau).\end{align*} Up until this point, the proof is valid in any Banach space and for any $1\leq p<\infty$. Finally, if $\XX$ is uniformly convex of power type $p\in[2,\infty)$, then by Lemma \ref{banach-filter} and the bound on the total number of filtrations, $$\sum_{Q\in\mathscr{A}} \beta_\Sigma(Q)^p\diam Q \lesssim_{p,A_\mathscr{H},\flatepsilon} \sum_j\sum_{\tau\in\mathscr{F}^j} \tilde\beta(\tau)^p\Diam(\tau) \lesssim_{p,\delta_{\XX},A_\mathscr{H},\flatepsilon} \Haus^1(\Sigma).$$ This verifies \eqref{A-sum}.
\end{proof}

\begin{lemma}[filtrations II]\label{l:H3} Assume that $\flatarcsepsilon$ is sufficiently small depending only on $A_\mathscr{H}$; $\flatarcsepsilon=1/126A_\mathscr{H}$ will suffice. If $\XX$ is uniformly convex of power type $p\in[2,\infty)$, then \begin{equation}\label{C-sum} \sum_{Q\in\mathscr{C}} \beta_\Sigma(Q)^p\diam Q\lesssim_{p,\delta_\XX,A_\mathscr{H},\flatepsilon} \Haus^1(\Sigma).\end{equation}\end{lemma}

\begin{proof} Modify the proof of \cite[Lemma 3.16]{Schul-Hilbert}. We initially assume that $\XX$ is any Banach space.
Let $Q\in\mathscr{C}$, say $Q=B(x_Q,A_\mathscr{H}2^{-k})$ for some $x_Q\in X_k$. Because $$\beta_{S^*(Q)}(2Q)\leq \epsilon_1 \beta_{\Lambda(Q)}(2Q),$$ we have $\Lambda(Q)\setminus S^*(Q)\neq\emptyset$ provided that $\flatarcsepsilon<1$. Choose any $\tau'_Q\in\Lambda(Q)\setminus S^*(Q)$. Then $$A_\mathscr{H}2^{-k}\leq \Diam(\tau'_Q)\leq 4A_\mathscr{H}2^{-k}\quad\text{and}\quad \beta_\Sigma(Q)\leq 2\beta_{\Lambda(Q)}(2Q)<\frac{1}{25\flatepsilon}\beta(\tau'_Q)$$ by definition of $\Lambda(Q)$ and $S^*(Q)$, respectively. Include an instance of $\tau'_Q$ in the set $\mathscr{F}^0_k$.

We claim that the family $\mathscr{F}^0=\bigcup_{k=k_0}^\infty \mathscr{F}^0_{k}$ of arcs is a prefiltration (see Lemma \ref{build-filtrations}), where $k_0$ is the smallest integer such that $\mathscr{C}$ contains a ball of radius $A_\mathscr{H}2^{-k_0}$. We already noted that arcs in $\mathscr{F}^0_k$ have good bounds on their diameters. To check the bounded overlap property, fix $Q$ as above. Using $\beta_{S^*(5Q)}(10Q)\leq \flatarcsepsilon \beta_{\Lambda(5Q)}(10Q)\leq \flatarcsepsilon$ for balls $Q\in\mathscr{C}$, we may choose a line $L_Q$ in $\XX$ such that \begin{equation}\label{G3-line} \dist(x,L_Q)< (21/20)\flatarcsepsilon\diam 10Q\leq 21\epsilon_1 A_\mathscr{H} 2^{-k}\quad\text{for all }x\in\bigcup_{\tau\in S^*(5Q)}\Image(\tau).\end{equation} Suppose that $R\in\mathscr{C}$, $R=B(x_R,A_{\mathscr{H}}2^{-k})$ for some $x_R\in X_k$, and $\tau'_R$ and $\tau'_Q$ intersect in their domains. Since $\tau'_Q$ is contained in $2Q$, $\tau'_R$ is contained in $2R$, and the arcs intersect, the triangle inequality yields $x_R\in 4Q$ and $x_Q\in 4R$. Hence $R\subset 5Q$ and $Q\subset 5R$, as well. Our strategy is to show that $x_R$ is close to $L_Q$ relative to the scale $2^{-k}$ of separation between points in $X_k$. By \eqref{G3-line}, it suffices to exhibit an arc $\xi_R\in S^*(5Q)$ containing $x_R$ and demand that $\flatarcsepsilon$ be sufficiently small relative to $A_\mathscr{H}$. To find $\xi_R$, first let $\gamma_R\in \Lambda(7R)$ be any arc containing $x_R$. Because $x_R$ lies in the net ball $B(x_R,(1/3)2^{-k})$ and $R\not\in(\mathscr{H}_0\cup\mathscr{A})$, we know that $\gamma_R\in S(7R)$, i.e.~$\beta(\gamma_R)\leq \flatepsilon\beta_\Sigma(R)$. Choose $\xi_R$ to be a subarc of the arc $\gamma_R$ such that $x_R\in\Image(\xi_R)$ and $\xi_R\in\Lambda(5Q)$, which exists because $x_R\in\Image(\gamma_R)\cap 4Q\subset 5Q$ and $10Q\subset 14R$. On the one hand, since $x_R\in\Image(\xi_R)$, $x_R\in 4Q$, and the endpoints of $\xi_R$ are contained in the boundary of $10Q$, we have $\Diam(\xi_R)\geq 6A_\mathscr{H}2^{-k}$. On the other hand, since $\Image(\gamma_R)\subset 14R$, $\Diam(\gamma_R)\leq 28A_\mathscr{H}2^{-k}$. Thus, as $28/6<5$ and $\Sigma\cap R\subset \Sigma\cap 5Q\subset\bigcup_{\tau\in\Lambda(5Q)}\Image(\tau)$, \begin{equation}\beta(\xi_R) \leq 5\beta(\gamma_R) \leq 5\flatepsilon\beta_\Sigma(R)\leq 25\flatepsilon\beta_\Sigma(5Q)\leq 50\flatepsilon\beta_{\Lambda(5Q)}(10Q)\end{equation} by three applications of \eqref{beta-monotone} and the observation that $\gamma_R$ is almost flat. This computation confirms that the arc $\xi_R$ belongs to $S^*(5Q)$. (In fact, this computation is the \emph{raison d'\^{e}tre} for $*$-almost flat arcs!) By \eqref{G3-line}, it follows that $\dist(x_R,L_Q)<(1/6)2^{-k}$ provided that $\flatarcsepsilon\leq 1/126A_\mathscr{H}$. For concreteness, we specify that $\flatarcsepsilon=1/126A_\mathscr{H}$. Let $\{x_1,\dots,x_N\}$ be an enumeration of the centers of balls $R\in\mathscr{C}$ of the same generation as $Q$ such that $\tau'_R$ intersects $\tau'_Q$ in the domain. Each $x_i\in 4Q$ and satisfies $\dist(x_i,L_Q)<(1/6)2^{-k}$. By Corollary \ref{c:var}, after reordering $x_1,\dots, x_N$, we have $$(N-1)2^{-k}\leq \sum_{i=1}^{N-1}|x_i-x_{i+1}| \leq (1+3/6)^2|x_1-x_N| \leq (9/4)\cdot 8A_\mathscr{H}2^{-k}=18A_\mathscr{H}2^{-k}.$$ We conclude that $N\leq 1+18A_\mathscr{H}\leq 19A_\mathscr{H}$. Therefore, $\mathscr{F}^0$ is a prefiltration.

Invoking Lemma \ref{build-filtrations} with parameters $\rho=2$, $A=4A_\mathscr{H}$, $\underline{A}=A_\mathscr{H}$, $C=19A_\mathscr{H}$, and $J=5$ (so that $2^J>6A/\underline{A}$), we can find $O(A_\mathscr{H})$ filtrations $$\mathscr{F}^1=\bigcup_{k=k_1}^\infty \mathscr{F}^1_k,\quad \mathscr{F}^2=\bigcup_{k=k_2}^\infty \mathscr{F}^2_k,\quad \mathscr{F}^3=\bigcup_{k=k_3}^\infty \mathscr{F}^2_k,\quad \dots$$ with starting indices $k_0\leq k_j \leq k_0+4$ such that $$\frac{1}{4}A_\mathscr{H}2^{(J-1)k_j}2^{-Jk} \leq \Diam(\tau) <8A_\mathscr{H} 2^{(J-1)k_j} 2^{-Jk}\quad\text{for all }\tau\in\mathscr{F}^j_k.$$ Note that the ratio of the upper and lower bounds for $\Diam(\tau)$ is universal. Moreover, for all $\tau'_Q\in\mathscr{F}^0_k$, there exists $\mathscr{F}^j$ (in the list of filtrations), an index $K$ with $k-k_j=J(K-k_j)$, and an arc $\tau_Q\in\mathscr{F}^j_K$ such that $$\Domain(\tau'_Q)\subset\Domain(\tau)\quad\text{and}\quad\Diam(\tau_Q)<2\Diam(\tau'_Q).$$ The assignment $(k,\tau'_Q)\mapsto (\mathscr{F}^j,K,\tau_Q)$ is injective. Now, for any $\tau'_Q\in\mathscr{F}^0_k$, we have that $$\beta_\Sigma(Q)\lesssim_{\flatepsilon}\beta(\tau'_Q)\leq 2\beta(\tau_Q)\leq 2\tilde\beta(\tau_Q),$$ where the first inequality holds by our choice of $\tau'_Q$, the second inequality holds by \eqref{beta-monotone}, because $\Diam(\tau)\leq 2\Diam(\tau')$, and the third inequality holds generally. Also, we have $\diam Q\simeq\Diam(\tau'_Q)\simeq \Diam(\tau_Q)$. Thus, \begin{align*}\sum_{Q\in\mathscr{C}} \beta_\Sigma(Q)^p\diam Q &\lesssim_{p,\flatepsilon} \sum_{\tau'\in\mathscr{F}^0}\beta(\tau')^p\Diam(\tau')
\\ &\lesssim_p \sum_j\sum_{\tau\in\mathscr{F}^j} \beta(\tau)^p\Diam(\tau)\leq \sum_j\sum_{\tau\in\mathscr{F}^j} \tilde\beta(\tau^p)\Diam(\tau).\end{align*} Up until this point, the proof is valid in any Banach space and for any $1\leq p<\infty$. Finally, if $\XX$ is uniformly convex of power type $p\in[2,\infty)$, then by Lemma \ref{banach-filter} and the bound on the total number of filtrations, $$\sum_{Q\in\mathscr{C}} \beta_\Sigma(Q)^p\diam Q \lesssim_{p,\flatepsilon} \sum_j\sum_{\tau\in\mathscr{F}^j} \tilde\beta(\tau)^p\Diam(\tau) \lesssim_{p,\delta_{\XX},A_\mathscr{H},\flatepsilon} \Haus^1(\Sigma).$$ This verifies \eqref{C-sum}.
\end{proof}

\begin{theorem}[geometric martingales] \label{l:H2} Assume that $\flatepsilon$ is sufficiently small depending only on $A_\mathscr{H}$ and $\flatarcsepsilon$; the value $\flatepsilon=2^{-55}\flatarcsepsilon/A_\mathscr{H}$ will suffice. If $\XX$ is any Banach space, then for all $q>0$, \begin{equation}\label{B-sum}\sum_{Q\in\mathscr{B}} \beta_\Sigma(Q)^q\diam Q \lesssim_{q,A_\mathscr{H},\flatarcsepsilon} \Haus^1(\Sigma),\end{equation} where the implicit constant blows up as $q\downarrow 0$.\end{theorem}

Note the arbitrary power $q$ in \eqref{B-sum}! This indicates that $Q\in\mathscr{B}$ (imagine $Q$ is centered at the intersection of two or more crossing line segments) occurs relatively infrequently. This result is universal insofar as it is valid in an \emph{arbitrary} Banach space! We defer the (rather long) proof of Theorem \ref{l:H2} to the sequel \cite{Badger-McCurdy-2}.

Putting it all together, on a uniformly convex space of power type $p\in[2,\infty)$, choose parameters $0<\flatarcsepsilon\ll 1$ and $0<\flatepsilon\ll 1$ with $\flatarcsepsilon\simeq A_{\mathscr{H}}^{-1}$ and $\flatepsilon\simeq A_{\mathscr{H}}^{-2}$ so that Lemma \ref{l:H3} and Theorem \ref{l:H2} are valid. Then $$\sum_{Q\in\mathscr{H}} \beta_\Sigma(Q)^p\diam Q \leq I + II + III + IV \lesssim_{p,\delta_\XX,A_\mathscr{H}}\Haus^1(\Sigma),$$ where $$I:= \sum_{Q\in\mathscr{H}_0} \beta_\Sigma(Q)^p\diam Q \leq \sum_{Q\in\mathscr{H}_0}\beta_\Sigma(Q)^2\diam Q\lesssim_{A_\mathscr{H}} \Haus^1(\Sigma)\qquad\text{ by Lemma \ref{l:H0}},$$
\begin{equation*}\begin{split}II:=  \sum_{Q\in\mathscr{A}} \beta_\Sigma(Q)^p\diam Q \lesssim_{p,\delta_\XX,A_\mathscr{H},\flatepsilon} \Haus^1(\Sigma) \qquad\text{by Lemma \ref{l:H1}},\end{split}\end{equation*}
\begin{equation*}\begin{split}III:=  \sum_{Q\in\mathscr{B}} \beta_\Sigma(Q)^p\diam Q \leq \sum_{Q\in\mathscr{B}} \beta_\Sigma(Q)^2\diam Q\lesssim_{A_\mathscr{H},\flatarcsepsilon} \Haus^1(\Sigma) \qquad \text{by Theorem \ref{l:H2}},\end{split}\end{equation*} \begin{equation*}\begin{split}\text{and}\quad IV:=  \sum_{Q\in\mathscr{C}} \beta_\Sigma(Q)^p\diam Q \lesssim_{p,\delta_\XX,A_\mathscr{H},\flatepsilon} \Haus^1(\Sigma) \qquad\text{by Lemma \ref{l:H3}.}\end{split}\end{equation*} This completes the proof of Theorem \ref{t:convex}.

\section{Analyst's TSP in finite-dimensional Banach spaces}\label{sec:finite-dimensional}

In this section, we record some additional observations about the Analyst's TSP in arbitrary finite-dimensional Banach spaces. It is well known that diameter, distance of a point to a set, and Hausdorff measures are bi-Lipschitz metric invariants in the sense that if $f:\XX\rightarrow\YY$ is a $C$-bi-Lipschitz map between metric spaces $\XX$ and $\YY$, then each of the quantities in the source and target spaces are quantitatively equivalent with multiplicative constants determined by $C$. As a consequence, the Jones' beta numbers in a normed vector space are bi-Lipschitz invariant in the following sense.

\begin{lemma}[Bi-Lipschitz equivalence of Jones' beta numbers] \label{beta-numbers-bilip} Let $V$ be a vector space and let $|\cdot|_\XX$ and $|\cdot|_\YY$ be equivalent norms on $V$, say for some constant $1\leq C<\infty$ that $C^{-1}|x|_\XX\leq |x|_\YY \leq C|x|_\XX$ for all $x\in V$. Let $\beta_{E,\XX}(Q)$ and $\beta_{E,\YY}(Q)$ denote the Jones' beta numbers in $\XX=(V,|\cdot|_\XX)$ and $\YY=(V,|\cdot|_\YY)$, respectively. That is, for all nonempty sets $E\subset V$ and positive diameter sets $Q\subset V$ with $E\cap Q\neq\emptyset$,
$$\beta_{E,\XX}(Q) = \inf_{L} \sup_{x\in E\cap Q}\frac{\dist_\XX(x,L)}{\diam_\XX Q}\quad\text{and}\quad\beta_{E,\YY}(Q) = \inf_{L} \sup_{x\in E\cap Q}\frac{\dist_\YY(x,L)}{\diam_\YY Q},$$ where $L$ runs over all lines in $V$.  Then $C^{-2} \beta_{E,\XX}(Q) \leq \beta_{E,\YY}(Q)\leq C^2\beta_{E,\XX}(Q)$ for all admissible $E$ and $Q$.\end{lemma}

\begin{proof} Let $E\subset V$ be nonempty and let $Q\subset V$ have positive diameter. By convention, $\beta_{E,\XX}(Q)=\beta_{E,\YY}(Q)=0$ if $E\cap Q=\emptyset$; in this case, the conclusion is trivial. Thus, we may assume that $E\cap Q\neq\emptyset$. By equivalence of the norms, $\diam_\XX Q=\sup\{|x-y|_\XX:x,y\in Q\}$ and $\diam_\YY Q=\sup\{|x-y|_\YY:x,y\in Q\}$ are also equivalent: \begin{equation}\label{diam-equivalent} C^{-1}\diam_\XX Q \leq \diam_\YY Q \leq C\diam_\XX Q.\end{equation} Similarly, for any nonempty set $S\subset V$ and $x\in V$, $\dist_\XX(x,S)=\inf\{|x-s|_\XX:s\in S\}$ and $\dist_\YY(x,S)=\inf\{|x-s|_\YY:s\in S\}$ are equivalent: \begin{equation}\label{dist-equivalent} C^{-1}\dist_\XX (x,S) \leq \dist_\YY (x,S) \leq C\dist_\XX (x,S).\end{equation} The set of all lines (i.e.~1-dimensional affine subspaces) in $V$ is independent of a choice of norm. Fix a line $L$ in $V$. By definition of $\beta_{E,\XX}(Q)$, \eqref{diam-equivalent}, and \eqref{dist-equivalent}, $$\beta_{E,\XX}(Q) \leq \sup_{x\in E\cap Q}\frac{\dist_\XX(x,L)}{\diam_\XX Q}
\leq C^2\sup_{x\in E\cap Q}\frac{\dist_\YY(x,L)}{\diam_\YY L}.$$ Taking the infimum over all possible lines, we obtain $\beta_{E,\XX}(Q)\leq C^2 \beta_{E,\YY}(Q)$. Interchanging the roles of $\XX$ and $\YY$ yields $\beta_{E,\YY}(Q)\leq C^2 \beta_{E,\XX}(Q)$.\end{proof}

A basic fact in functional analysis is that every pair of finite-dimensional normed vector spaces of the same dimension have equivalent norms. Therefore, the original formulation of the Analyst's Traveling Salesman Theorem due to Jones \cite{Jones-TST} ($\dim V=2$) and Okikiolu \cite{Ok-TST} ($\dim V\geq 3$) persists in any finite-dimensional vector space.

\begin{theorem}[Jones' and Okikiolu's theorems in finite-dimensional spaces] \label{t:Jones-again} Let $\XX$ be a finite-dimensional Banach space with $\dim \XX\geq 2$, let $\Delta(\XX)$ be a system of dyadic cubes with respect to some set of coordinates on $\XX$, and let $E\subset \XX$. Then $E$ is contained in a rectifiable curve if and only if \begin{equation}\label{e:S-finite-again} S_E(\XX):= \diam E + \sum_{Q\in\Delta(\XX)} \beta_{E}(3Q)^2 \diam Q<\infty,\end{equation} where $3Q$ denotes the concentric dilate of the cube $Q$ with scaling factor 3 and $\beta_E(3Q)$ denotes the Jones' beta number with respect to $\XX$. More precisely, if $S_E(\XX)<\infty$, then $E$ is contained in a curve $\Gamma$ in $\XX$ with \begin{equation}\label{e:Jones-again} \Haus^1(\Gamma) \lesssim_{\XX,\Delta(\XX),\dim \XX} S_E(\XX).\end{equation} If $\Sigma\subset \XX$ is a connected set, then \begin{equation} \label{e:Okikiolu-again} S_{\Sigma}(\XX) \lesssim_{\XX,\Delta(\XX),\dim\XX} \Haus^1(\Sigma).\end{equation} The constant $3$ in \eqref{e:S-finite-again} can be replaced with any constant $A>1$. Then \eqref{e:Jones-again} and \eqref{e:Okikiolu-again} hold with implicit constants depending on the norm of $\XX$, $\Delta(\XX)$ (i.e.~a choice of coordinates), $\dim\XX$, and $A$.\end{theorem}

\begin{proof} Suppose that $\dim \XX=n$. After fixing coordinates on $\XX$, we may identify $\XX$ with $\RR^n$ and identify $\Delta(\XX)$ with $\Delta(\RR^n)$. Then $(\XX,|\cdot|)=(\RR^n,|\cdot|)$ and $(\RR^n,|\cdot|_2)$ are bi-Lipschitz equivalent, where $|\cdot|_2$ denotes the standard Euclidean norm. (The bi-Lipschitz constant depends on the norm $|\cdot|$ and on the choice of coordinates for $\XX$.) Each of the quantities $\diam E$, $\beta_{E}(3Q)$, and $\Haus^1(\Gamma)$ defined relative to $|\cdot|$ are bi-Lipschitz equivalent to the respective quantities defined relative to $|\cdot|_2$. Thus, rectifiability of a curve $\Gamma$ containing $E$ is independent of the choice of norm (although the length of $\Gamma$ depends on the norm) and the theorem follows immediately from Theorem \ref{t:Jones}.\end{proof}

\begin{example}Given linearly independent vectors $\vec v$ and $\vec w$ in $\RR^2$, let $\Delta(\vec v,\vec w)$ denote the system of dyadic parallelograms corresponding to the lattice generated by $\vec v$ and $\vec w$; e.g., $P=\{s\vec v+t\vec w:0\leq s,t\leq 1\}$ and $P'=\{s\vec v+t\vec w:0\leq s,t\leq 1/2\}$ belong to $\Delta(\vec v,\vec w)$ and $P'$ is one of four children of $P$. A bounded set $E\subset \RR^2$ is contained in a rectifiable curve if and only if $$\sum_{P\in\Delta(\vec v,\vec w)} \beta_E(3P)^2\diam P<\infty.$$ To see this, note that every $P\in\Delta(\vec v,\vec w)$ is a dyadic square in $\RR^2$ with respect to the coordinates induced by $\vec v$ and $\vec w$ and apply Theorem \ref{t:Jones-again}.\end{example}

In a finite-dimensional Banach space, the Analyst's Traveling Salesman Theorem can be formulated using cubes, as in Theorems \ref{t:Jones} and \ref{t:Jones-again}, and using multiresolution families, as in Theorems \ref{t:schul}, \ref{t:main}, and \ref{t:main2}. It is always possible to pass between one and the other, at the expense of growing implicit constants. For a more general formulation of this principle, see \cite[Appendix B]{Bishop-TST}.

\begin{lemma}\label{change-resolutions} Let $\XX$ be a finite-dimensional Banach space with $\dim \XX\geq 2$, let $\Delta(\XX)$ be a system of dyadic cubes with respect to some set of coordinates on $\XX$, let $a>1$ be a scaling factor for cubes, let $E\subset \XX$, and let $\mathscr{G}$ be a multiresolution family for $E$ with inflation factor $A_\mathscr{G}>1$. For all $0< p<\infty$, \begin{equation}\label{cubes-balls}\sum_{Q\in\Delta(\XX)} \beta_E(aQ)^p\diam Q \simeq_{p,a,A_\mathscr{G},\dim\XX} \sum_{B\in\mathscr{G}} \beta_E(B)^p\diam B.\end{equation}\end{lemma}

\begin{proof}
The key point is that any bounded set in $\XX$ is contained in some dilate $aQ$ of a cube $Q\in\Delta(\XX)$ with comparable diameters. As long as the bounded set intersects $E$, that set is also contained in some ball $B\in\mathscr{G}$ with comparable diameters. We will use this in conjunction with monotonicity of Jones' beta numbers and volume doubling in finite-dimensional spaces.

To normalize scales, let $D=\diam Q_0$ for any choice of $Q_0\in\Delta(\XX)$. Fix $Q\in \Delta(\XX)$, say with $\diam Q=D2^{-j}$ for some $j\in\ZZ$. Fix $k\in\ZZ$ to be determined. If $aQ\cap E=\emptyset$, then $\beta_E(aQ)=0$, so the cube $Q$ is irrelevant to the sum on the left hand side of \eqref{cubes-balls}. Suppose that $aQ\cap E\neq\emptyset$ and fix any $z\in aQ\cap E$. Choose $x\in X_k$ (the $2^{-k}$-net appearing in the definition of $\mathscr{G}$) such that $|x-z|<2^{-k}$ and set $B=B(x,A_\mathscr{G}2^{-k})\in \mathscr{G}$. Then for any $y\in aQ$, $$|y-x|\leq |y-z|+|z-x|<  \diam aQ + 2^{-k}=aD2^{-j}+2^{-k} \leq A_{\mathscr{G}}2^{-k}$$ and $aQ\subset B$ provided that $aD2^{-j}\leq (A_\mathscr{G}-1)2^{-k}$. We now specify that $k$ is the unique integer such that \begin{equation}\label{cube-scale-comparison}(A_\mathscr{G}-1)2^{-(k+1)}<aD2^{-j}\leq (A_\mathscr{G}-1)2^{-k}.\end{equation} As noted, this ensures that $aQ\subset B$. Furthermore, $$\diam aQ\leq \diam B \leq 2A_\mathscr{G}2^{-k} \leq \frac{4A_\mathscr{G}}{A_\mathscr{G}-1}\diam aQ.$$ Thus, $\beta_E(aQ) \lesssim_{A_\mathscr{G}} \beta_E(B)$ by \eqref{beta-monotone}.  Now, by volume doubling, the sets of the form $aR$, where $R\in\Delta(\XX)$ is a dyadic cube of the same generation as $Q$, have bounded overlap determined by $a$ and $\dim \XX$. Also, by \eqref{cube-scale-comparison}, each $k\in\ZZ$ is associated to an unique $j\in\ZZ$. (This is a simple expression of the fact that $2^{-k}$-nets and dyadic cubes have the same scaling ratios.) All together, we conclude that $$\sum_{Q\in\Delta(\XX)} \beta_E(aQ)^p\diam Q \lesssim_{p,a,A_\mathscr{G},\dim\XX} \sum_{B\in\mathscr{G}} \beta_E(B)^p\diam B.$$ The reverse inequality holds by similar considerations.
\end{proof}

\section{Sharpness of the exponents via examples}\label{sec:examples}

Our goal in this section is to verify the sharpness of the exponents on beta numbers in Theorems \ref{t:main} and \ref{t:main2}. To do so, we build Koch-snowflake-like curves $\Gamma$, for which we can estimate beta number sums over \emph{arbitrary} multiresolution families. This type of construction is not new, see e.g.~\cite{BJ, Rohde-snowflakes, BNV, ENV-Banach} for motivating examples, but the details are subtle.

The organization is as follows. In \S\ref{ss:exponent-2}, we verify sharpness of the exponent 2 in \eqref{e:suff-p} when $2\leq p<\infty$ and in \eqref{e:nec-2} when $1<p\leq 2$ by building curves in the Euclidean plane.  In \S\ref{ss:exponent-p}, we verify sharpness of the exponent $p$ in \eqref{e:suff-p} when $1<p\leq 2$ and in \eqref{e:nec-2} when $2\leq p<\infty$ by building curves in infinite-dimensional Banach spaces. Finally, in \S\ref{ss:pq}, we carry out additional estimates to record a proof of Proposition \ref{t:pq}.

\subsection{Examples with critical exponent 2}\label{ss:exponent-2}
The main results of this subsection are the following two propositions. Recall that $\ell_p^n$ denotes $(\RR^n,|\cdot|_p)$.

\begin{proposition}\label{prop:2} There exists a curve $\Gamma$ in $\ell_2^2$ such that $\Haus^1(\Gamma)=\infty$ and $S_{\Gamma,2+\epsilon}(\mathscr{G})<\infty$ for every multiresolution family $\mathscr{G}$ for $\Gamma$ and every $\epsilon>0$. In particular, the exponent 2 in \eqref{e:suff-p} when $2\leq p<\infty$ is sharp.\end{proposition}

\begin{proposition}\label{prop:4} There exists a curve $\Gamma$ in $\ell_2^2$ such that $\Haus^1(\Gamma)<\infty$ and $S_{\Gamma,2-\epsilon}(\mathscr{G})=\infty$ for every multiresolution family $\mathscr{G}$ for $\Gamma$ and every $\epsilon>0$. In particular, the exponent 2 in \eqref{e:nec-2} when $1<p\leq 2$ is sharp.\end{proposition}

Sharpness of the exponent 2 in \eqref{e:suff-p} when $2<p<\infty$ and in \eqref{e:nec-2} when $1<p<2$ follows from the case $p=2$, because $\ell_p$ contains a subspace isomorphic to $\ell^2_p$, which in turn is bi-Lipschitz equivalent to $\ell^2_2$. (Use Lemmas \ref{beta-numbers-bilip} and \ref{change-resolutions} to pass between multiscale sums of beta numbers in $\ell_p^2$ and $\ell_2^2$. In particular, one may use a fixed system of dyadic squares as an intermediary between multiresolution families for $\Gamma$ in $\ell_p^2$ and in $\ell_2^2$.) Therefore, in this section, we focus on $\ell_2^2\equiv \RR^2,$ the standard Euclidean plane. We wish to emphasize that in both statements the curve is independent of $\epsilon$. We build the curves using the following procedure.

\begin{figure}
\begin{center}
\includegraphics[width=.8\textwidth]{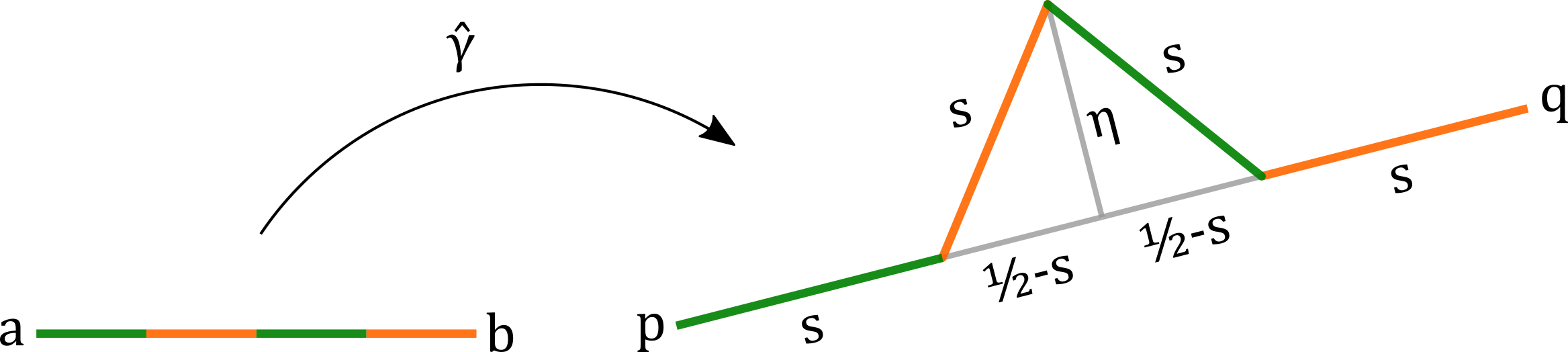}
\end{center}\caption{The snowflake map $\hat\gamma$, displayed with $|\vec v|=1$, $\eta=1/4$, and $s=5/16$.}\label{figure:snowflake}\end{figure}

\begin{algorithm}[snowflake-like curves in $\RR^2$] \label{alg:2}
Suppose that $\vec p,\vec q\in\RR^2$ and $\gamma:[a,b]\rightarrow\RR^2$ is a constant speed parameterization of $I=[\vec p,\vec q]$ from $\gamma(a)=\vec p$ to $\gamma(b)=\vec q$. That is, $$\gamma(t)=\vec p + \frac{t-a}{b-a}(\vec q-\vec p)\quad\text{for all }t\in[a,b].$$ Write $\vec v:=\vec q-\vec p$ and let $\vec y\in\RR^2$ denote the unique unit vector with $\vec y\perp \vec v$ such that $\vec y$ points to the left of the oriented line segment from $\vec p$ to $\vec q$. Given a \emph{relative height} $0 \leq  \eta \le 1/\sqrt{12}$, we define a piecewise linear path $\hat\gamma:[a,b]\rightarrow\RR^2$, as follows. Set $s:=1/4+\eta^2\in[1/4,1/3]$. Divide $[a,b]$ into quarters. For all $a\leq t \leq (3/4)a+(1/4)b$,
    $$\hat\gamma(t) = \vec p + 4\left(\frac{t-a}{b-a}\right)s\vec v$$ (see Figure \ref{figure:snowflake}).
For all $(3/4)a + (1/4)b \leq t \leq (1/2)a+(1/2)b$,
    $$\hat\gamma(t) = \vec p + s\vec v + 4\left(\frac{t-(3/4)a-(1/4)b}{b-a}\right) \big((1/2-s)\vec v + \eta|\vec v|\vec y\big).$$
For all $(1/2)a+(1/2)b\leq t\leq (1/4)a + (3/4)b$,
    $$\hat\gamma(t) = \vec p + (1/2)\vec v + \eta|\vec v|\vec y + 4\left(\frac{t-(1/2)a-(1/2)b}{b-a}\right)\big((1/2-s)\vec v-\eta|\vec v|\vec y\big).$$
Finally, for all $(1/4)a+(3/4)b\leq t\leq b$,
    $$\hat\gamma(t) = \vec p + (1-s)\vec v + 4\left(\frac{t-(1/4)a-(3/4)b}{b-a}\right)s\vec v.$$ We say that $\hat \gamma$ is obtained from $\gamma$ by \emph{adding a bump of relative height $\eta$ (on the left side)}. On each quarter of $[a,b]$, $\hat \gamma(t)$ traces a line segment of length $s|\vec v|$ at constant speed. Thus, $\hat\gamma$ has Lipschitz constant $4s|\vec v|/(b-a)=(1+4\eta^2)|\vec v|/(b-a)$. Additionally, we have $\|\gamma-\hat\gamma\|_\infty =|\gamma((a+b)/2)-\hat\gamma((a+b)/2)|=\eta|\vec v|$.

Starting from the arc length parameterization $\gamma_0:[0,1]\rightarrow\RR^2$ of the line segment $I_0=[\vec 0,\vec e_1]$ from $\gamma_0(0)=\vec 0$ to $\gamma_0(1)=\vec e_1$, we now define a sequence of piecewise linear maps $\gamma_i:[0,1]\rightarrow\RR^2$ by iteratively adding bumps of relative height $\eta_i$. Suppose that $\gamma_i$ has been defined for some $i\geq 0$ so that $\gamma_i|_{J_{i,k}}$ is a constant speed parameterization of a line segment on each interval $J_{i,k}=[a_{i,k},b_{i,k}]$ of the form $$J_{i,k}=\left[a+\frac{(k-1)}{4^i}(b-a), a + \frac{k}{4^i}(b-a)\right]\qquad (1\leq k \leq 4^i,\  k\in\ZZ).$$ For each index $1\leq k\leq 4^i$, define $\gamma_{i+1}|_{J_{i,k}} = \widehat{\gamma_i|_{J_{i,k}}}$ by adding a bump of relative height $\eta_{i+1}$. This defines a map $\gamma_{i+1}$. By induction, we obtain the full sequence $\gamma_0, \gamma_1,\gamma_2,\dots$; the image of each map $\gamma_n$ is composed of $4^n$ line segments of length $r_n:=(\Lip\gamma_n)4^{-n}$. Evidently, for all $n=0,1,2,\dots$ \begin{align}\label{gamma-lip} \Lip (\gamma_n) &= \prod_{i=1}^n (1+4\eta_i^2) \leq (4/3)^n, \\ \label{r-comparison} r_{n+1} &\leq (1/3) r_n, \\ \label{gamma-difference} \|\gamma_n-\gamma_{n+1}\|_\infty &\leq \eta_{n+1} r_n\leq 3^{-n}\eta_{n+1}<3^{-n}.\end{align} Therefore, $\gamma:[0,1]\rightarrow\RR^2$, which is defined pointwise by \begin{equation} \label{gamma-definition} \gamma(t)=\gamma_0(t)+\sum_{n=0}^\infty(\gamma_{n+1}(t)-\gamma_n(t))\quad\text{for all }t\in[0,1],\end{equation} is continuous as the uniform limit of the maps $\gamma_n$ by \eqref{gamma-difference}. We denote $\gamma([0,1])$ by $\Gamma$, and for each $n$, we denote $\gamma_n([0,1])$ by $\Gamma_n$. For all integers $n,k\geq 0$ with $0\leq k\leq 4^n$, we call the point $\gamma_n(k/4^n)$ a \emph{vertex} of $\Gamma_n$.\end{algorithm}

\begin{remark}[modulus of continuity] In fact, \eqref{gamma-lip} and \eqref{gamma-difference} imply that the maps $\gamma_n$ and $\gamma$ are uniformly $\log_4(3)$-H\"older continuous (see e.g.~ \cite[Appendix B]{BNV}). This is the optimal modulus of continuity of the von Koch snowflake curve, which corresponds to the choice of relative heights $\eta_n=1/\sqrt{12}$ for all $n$. Furthermore, if $\sum_{n=1}^\infty \eta_n^2<\infty$, then \eqref{gamma-lip} implies that $\gamma_n$ and $\gamma$ are uniformly Lipschitz with Lipschitz constant at most $\prod_{n=1}^\infty (1+4\eta_n^2)\simeq \exp(\sum_{n=1}^\infty \eta_n^2)$.
\end{remark}

\begin{remark}[vertices] \label{vertex facts} For each integer $n\geq 0$, let $V_n$ denote the set of vertices in $\Gamma_n$ and let $\widetilde V_n:=V_n\setminus V_{n-1}$ denote the set of ``new vertices" in $\Gamma_n$ with the convention that $\widetilde V_0=V_0$. For later use, we observe that \begin{equation}V_0\subset V_1\subset V_2\subset\cdots,\end{equation}
\begin{equation}V_n=\widetilde{V}_0\cup \widetilde{V}_1\cup\dots  \cup \widetilde V_{n-1}\cup \widetilde{V}_n\qquad (\widetilde{V}_i\cap \widetilde{V}_j=\emptyset\text{ when }i\neq j),\end{equation}
\begin{equation}\#\widetilde{V}_0=2, \quad \#\widetilde{V}_n = 3\cdot 4^{n-1}\ (n\geq 1),\quad \#V_n=1+4^n.\end{equation}
\end{remark}

\begin{remark}[injectivity] \label{r:injective} The restriction $\eta_i\leq 1\sqrt{12}$ on the relative heights ensures that the parameterization $\gamma$ of $\Gamma$ constructed by Algorithm \ref{alg:2} is injective. This can be shown by a geometric argument similar to the proof that the standard von Koch snowflake curve is the attractor of an iterated function system that satisfies the open set condition (i.e.~ draw equilateral triangles on the left side of each segment in $\Gamma_n$). We leave details to the dedicated reader. For other models of generalized von Koch curves, the question of injectivity of a parameterization is quite subtle, see e.g.~\cite{Keleti-Paquette}.\end{remark}

\begin{lemma}\label{l:gamma-length} If $\Gamma$ is constructed by Algorithm \ref{alg:2}, then \begin{equation}\label{length-exp} \exp\left(3.3\overline{3}\textstyle\sum_{n=1}^\infty \eta_n^2\right) \leq \Haus^1(\Gamma)=\prod_{n=1}^\infty (1+4\eta_n^2) \leq \exp\left(\textstyle4\sum_{n=1}^\infty \eta_n^2\right).\end{equation} In particular, $\Gamma$ is a rectifiable curve if and only if $\sum_{n=1}^\infty \eta_n^2<\infty$. Moreover, in that case, $\Gamma$ is Ahlfors regular with constants depending only on $\Haus^1(\Gamma)$.\end{lemma}

\begin{proof} The one-dimensional Hausdorff measure $\Haus^1$ enjoys the bound $\Haus^1(K)\geq \diam K$ for every connected set $K\subset\RR^2$ (see e.g.~\cite[Lemma 2.11]{AO-curves}).  Thus, for each $n\geq 1$, \begin{align*}\Haus^1(\Gamma) &= \sum_{k=1}^{4^n} \Haus^1(\gamma([(k-1)4^{-n},k4^{-n}])) \geq \sum_{k=1}^{4^n} \diam \gamma([(k-1)4^{-n},k4^{-n}]) \\
&\geq \sum_{k=1}^{4^n} \diam \gamma_n([(k-1)4^{-n},k4^{-n}]) = \sum_{k=1}^{4^n} \prod_{i=1}^n (1/4+\eta_i^2)= \prod_{i=1}^n (1+4\eta_i^2),
\end{align*} where the initial equality holds by Remark \ref{r:injective}. Hence $\Haus^1(\Gamma)\geq \prod_{n=1}^\infty (1+4\eta_n^2)$. Conversely, $\Haus^1(\Gamma)\leq \Lip\gamma \leq \prod_{n=1}^\infty (1+4\eta_n^2)$. Therefore, $$\Haus^1(\Gamma) = \prod_{n=1}^\infty (1+4\eta_n^2)= \exp\left(\sum_{n=1}^\infty \log(1+4\eta_n^2)\right).$$ Finally, since $0\leq 4\eta_n^2 \leq 1/3$, we may use Taylor's theorem to bound $\log(1+4\eta_n^2) \leq 4\eta_n^2$ and $\log(1+4\eta_n^2)\geq 4\eta_n^2 - \frac12(4\eta_n^2)^2 \geq (10/3)\eta_n^2.$

Suppose that $\Haus^1(\Gamma)<\infty$. By \eqref{r-comparison}, \eqref{gamma-difference}, and the bound $\eta_n\leq 1/\sqrt{12}$ for all $n\geq 1$, $$\|\gamma_n-\gamma\|_\infty \leq \sum_{j=n}^\infty \eta_{j+1}r_j \leq \frac{1}{\sqrt{12}}\sum_{k=0}^\infty r_{n+k} \leq  \frac{1}{\sqrt{12}}r_n\sum_{k=0}^\infty 3^{-k}=\frac{\sqrt{3}}{4}r_n<0.45r_n.
$$ Thus, given $x\in\Gamma$ and $n\geq 1$, we may pick $y\in \Gamma_n$ such that $|x-y|\leq (\sqrt{3}/4)r_n$. Next, let $v\in V_n$ be an endpoint of the segment in $\Gamma_n$ containing $y$ that is closest to $y$ so that $|v-y|\leq (1/2)r_n$. Hence $B(x, 0.05r_n)\subset B(v,r_n)$ and we may use \eqref{length-exp} to estimate \begin{equation}\label{gamma-AR} \Haus^1(\Gamma\cap B(x,0.05r_n)) \leq \Haus^1(\Gamma\cap B(v,r_n)) \leq 2 \exp\left(\sum_{i=n}^\infty 4\eta_i^2\right)r_n \leq 2 e^{6/5} \Haus^1(\Gamma)r_n,\end{equation} because $\Gamma_n\cap B(v,r_n)$ consists of one or two line segments of length $r_n$. From \eqref{gamma-AR}, one easily deduces that $\Haus^1(\Gamma\cap B(x,r))\simeq_{\Haus^1(\Gamma)} r$ for every $x\in \Gamma$ and every $0<r\leq 1=\diam \Gamma$. Therefore, $\Gamma$ is Ahlfors regular with constants depending only on $\Haus^1(\Gamma)$.
\end{proof}

\begin{lemma}\label{beta bound} Assume that $\Gamma$ is constructed by Algorithm \ref{alg:2}. For all  $i\geq 0$ and $j\geq i$, \begin{equation}\frac{\sqrt{3}}{4}\eta_i\leq \beta_{\Gamma_j}(B(v,r_j)) \leq 2\eta_i\quad\text{for all }v\in \widetilde{V}_i,\end{equation} where $\widetilde{V}_i$ denote the ``new vertices" in $\Gamma_i$ (see Remark \ref{vertex facts}) and $\eta_0=0$. Furthermore, \begin{equation} \frac{\sqrt{3}}{4}\eta_i \leq \beta_{\Gamma}(\overline{B(v,r_i)})\quad\text{for all $i\geq 0$ and }v\in\widetilde{V}_i.\end{equation}
\end{lemma}

\begin{proof} Note that when $v\in \widetilde{V}_0$ is an endpoint of $\Gamma_j$, the set $\Gamma_j\cap B(v,r_j)$ is a line segment. Thus, $\beta_{\Gamma_j}(B(v,r_j))=0$ for all $v\in\widetilde{V}_0$ and $j\geq 0$. Next, suppose that $i\geq 1$ and $v\in \widetilde{V}_i$. Using the line containing the segment in $\Gamma_{i-1}\cap B(v,r_i)$ to approximate the beta number, we obtain the estimate $$\beta_{\Gamma_i}(B(v,r_i)) \leq \frac{\eta_i}{2(\frac{1}{4}+\eta_i^2)}\leq 2\eta_i.$$ When $j\geq i$, the set $\Gamma_j\cap B(v,r_j)$ agrees up to a dilation centered at $v$ with $\Gamma_i\cap B(v,r_i)$. Thus, $\beta_{\Gamma_j}(B(v,r_j)) \leq 2\eta_i$ for all $v\in\widetilde{V}_i$ and $j\geq i$, as well.

To establish the lower bound, one can use symmetry to find the best fitting line for each of the two triangles formed by consecutive segments of side length $s$ in Figure \ref{figure:snowflake}. The best fitting lines are parallel to the third side of each triangle and the altitudes of the triangles are $\eta$ and $\sqrt{s}\eta=(\frac{1}{4}+\eta^2)^{1/2}\eta$. It follows that for all $v\in\widetilde{V}_i$ and $j\geq i$, $$\beta_{\Gamma_j}(B(v,r_j)) =\beta_{\Gamma_i}(B(v,r_i))\geq \frac{\frac{1}{2}(\frac14+\eta_i^2)^{1/2}\eta_i}{2(\frac{1}{4}+\eta_i^2)} =\frac{\eta_i}{4(\frac{1}{4}+\eta_i^2)^{1/2}}\geq \frac{\sqrt{3}}{4}\eta_i,$$ since $\eta_i\leq 1/\sqrt{12}$. Finally, because $\Gamma\cap \overline{B(v,r_i)}$ contains the vertices $a,b,c$ of one of the two triangles, we find that $\beta_{\Gamma}(\overline{B(v,r_i)}) \geq \beta_{\{a,b,c\}}(\overline{B(v,r_i)})\geq (\sqrt{3}/{4})\eta_i$, as well.
\end{proof}

\subsubsection{Proof of Proposition \ref{prop:2}}
Build $\Gamma$ using Algorithm \ref{alg:2} with relative heights $4\eta_i^2 := 1/(i+15)\log(i+15)$ for all $i\geq 1$. See Remark \ref{remark:log} below for an explanation of the logarithmic factor in the relative heights. Note that $\eta_{i+1}\leq \eta_i<1/8$ for all $i$. Because $\sum_{i=1}^\infty \eta_i^2=\infty$, we have $\Haus^1(\Gamma)=\infty$ by Lemma \ref{l:gamma-length}. To proceed, let $(X_k)_{k\in\ZZ}$ be any family of nested $2^{-k}$-nets for $\Gamma$, and let $\mathscr{G}=\{B(x,A2^{-k}):x\in X_k, k\in\ZZ\}$ be the associated multiresolution family with inflation factor $A>1$. Since $\Gamma$ is bounded, to prove that $S_{\Gamma,2+\epsilon}(\mathscr{G})<\infty$, it suffices to show that (cf.~Lemma \ref{l:H0}) $$S_{\Gamma,2+\epsilon}(\mathscr{G}'):=\sum_{B\in \mathscr{G}'} \beta_\Gamma(B)^{2+\epsilon}\diam B<\infty,$$ where $\mathscr{G}'$ is the subfamily of all balls starting from some initial generation $k_0$.

Recall that each intermediate curve $\Gamma_n$ consists of $4^n$ line segments of length $$r_n = 4^{-n} \prod_{i=1}^n (1+4\eta_i^2)=4^{-n} \exp\left(\sum_{i=1}^n \log(1+4\eta_i^2)\right).$$ By Taylor's theorem, we have $$4\eta_i^2- \frac{1}{2}(4\eta_i^2)^2\leq \log(1+4\eta_i^2) \leq 4 \eta_i^2\quad\text{for all }i\geq 1.$$ Combined with the elementary bounds \begin{equation*}\begin{split} \sum_{i=1}^n \frac{1}{(i+15)\log(i+15)} &\leq \frac{1}{16\log(16)}+\int_{16}^{n+15} \frac{1}{x\log(x)}\,dx \\ &<\frac{1}{16}+\log(\log(n+15))\end{split}\end{equation*} and \begin{equation*}\begin{split}\sum_{i=1}^n &\frac{1}{(i+15)\log(i+15)}-\frac{1}{2(i+15)^2\log(i+15)^2} \\ &\qquad \geq \int_{16}^{n+16} \left(\frac{1}{x\log(x)}-\frac{1}{2x^2\log(x)^2}\right)dx>\log(\log(n+16))-1.2\end{split}\end{equation*} we obtain the rough estimate (for the record, $0.3<e^{-1.2}$ and $e^{1/16}< 1.1$): \begin{equation}\label{e:rough-r} 0.3\log(n+16)4^{-n}< r_n< 1.1\log(n+15)4^{-n}\quad\text{for all }n\geq 1.\end{equation}

To continue, let $k\geq k_0$ (with $k_0$ sufficiently large depending on $A$ to be specified below) and let $x\in X_k$. Our immediate goal is to estimate $\beta_\Gamma(B(x,A 2^{-k}))$ from above in terms of $\beta_{\Gamma_m}(B(v, r_m))$ for some suitably chosen generation $m=m(A,k)\geq 1$ and vertex $v\in \Gamma_m$. Write $c:=\log_4 A>0$ so that $A2^{-k}=4^{-(\frac12 k -c)}$. We now require $\frac12k_0-c \geq 3$, which ensures $$A2^{-k_0}\leq 4^{-3} < 0.03\log(16)4^{-1}.$$ Since $k\geq k_0$, there exists a unique integer $m\geq 1$ with \begin{equation}\label{k-to-m} 0.03\log(m+16)4^{-(m+1)}< A2^{-k}\leq 0.03\log(m+15)4^{-m}.\end{equation} By \eqref{e:rough-r}, it follows that \begin{equation}\label{scale-above} r_m \lesssim A2^{-k}<\frac{1}{10}r_m.\end{equation} Next, by \eqref{r-comparison}, \eqref{gamma-difference}, and fact that $\eta_{i+1}\leq \eta_i<1/8$ for all $i$, $$\|\gamma_m-\gamma\|_\infty \leq \sum_{j=m}^\infty \eta_{j+1}r_j \leq \eta_{m+1}\sum_{l=0}^\infty r_{m+l} \leq  \eta_{m+1}r_m\sum_{k=0}^\infty 3^{-k}=(3/2)\eta_{m+1}r_m <(3/16)r_m.$$ In particular, we can find $y\in \Gamma_m$ with $|x-y|< (3/16)r_m$ and then choose a vertex $v$ in $\Gamma_m$ such that $|y-v|\leq (1/2)r_m$ (i.e.~ $v$ is an endpoint of the segment containing $y$). Since $$\frac{1}{10}+\frac{3}{16}+\frac{1}{2} <\frac{3}{4},$$  $B(x,A2^{-k})\subset B(v,(3/4)r_m)$. Invoking the bound $\|\gamma-\gamma_m\|_\infty \leq (3/2)\eta_{m+1}r_m< (3/16)r_m$ again, we conclude that \begin{equation*} E:=\excess(\Gamma\cap B(x,A2^{-k}),\Gamma_m\cap B(v,r_m))\leq (3/2)\eta_{m+1}r_m,\end{equation*} where $\excess(S,T)=\sup_{s\in S}\inf_{t\in T}|s-t|$ denotes the \emph{excess of $S$ over $T$}. Hence \begin{equation}\label{pre-beta}\beta_{\Gamma}(B(x,A2^{-k})) \leq \frac{E}{2A2^{-k}}+ \frac{r_m}{A2^{-k}}\cdot \beta_{\Gamma_m}(B(v,r_m)) \lesssim \eta_{m+1} + \beta_{\Gamma_m}(B(v,r_m))\end{equation} by the triangle inequality and \eqref{scale-above}. Note that by Lemma \ref{beta bound}, \begin{equation}\label{beta-above} \eta_{m+1} + \beta_{\Gamma_m}(B(v,r_m)) \leq \eta_{m+1}+2\eta_i\quad\text{when }v\in \widetilde{V}_i\subset V_m,\end{equation} where $\eta_0=0$. In particular, $\beta_\Gamma(B(x, A2^{-k}))\lesssim \eta_{m+1}$ when $v\in\widetilde V_0$ and $\beta_\Gamma(B(x,A2^{-k})) \lesssim \eta_i$ when $v\in\widetilde V_i$ for some $1\leq i\leq m$.

Next, we bound the number of times a scale $m$ and vertex $v\in V_m$ are associated to a point $x\in X_k$. On one hand, $\# X_k\cap B(v,r_m) \lesssim A^2$ for each vertex $v\in V_m$ by \eqref{scale-above}, since $X_k$ is $2^{-k}$-separated and our construction takes place in $\RR^2$. (We could remove the dimension dependence by using Lemma \ref{l:graph}, but do not require a sharp upper bound.) On the other hand, the scale $m$ associated to $k\geq k_0$ satisfies $$0.03\log(m+16)4^{-(m+1)}< 4^{-(\frac12 k-c)}\leq 0.03\log(m+15)4^{-m},$$ where $c=\log_4 A$. Taking logarithms and rearranging, we have $$2\big(m+c-\log_4(0.03\log (m+15))\big) \leq k < 2\big(m+1+c-\log_4(0.03\log(m+16))\big).$$ Thus, the number of integers $k$ associated to a given integer $m$ is at most $$2\big(m+1+c-\log_4(0.03\log(m+16))\big)-2\big(m+c-\log_4(0.03\log(m+15))\big) \leq 2.$$

To finish, fix a parameter $\epsilon>0$. In view the previous paragraph, \eqref{k-to-m}, and \eqref{pre-beta}, we see that \begin{equation}\begin{split}S_{\Gamma,2+\epsilon}(\mathscr{G}')&=\sum_{k=k_0}^\infty 2A2^{-k}\sum_{x\in X_k} \beta_\Gamma(B(x,A2^{-k}))^{2+\epsilon}
\\ &\lesssim_{A,\epsilon} \sum_{m=1}^\infty \log(m+15)4^{-m}\sum_{v\in V_m} (\eta_{m+1}+\beta_{\Gamma_m}(B(v,r_m))^{2+\epsilon}.\end{split}\end{equation} Decomposing $V_m=\widetilde{V}_0\cup\dots\cup\widetilde{V}_m$ (see Remark \ref{vertex facts}) and invoking \eqref{beta-above},
\begin{equation}\sum_{v\in V_m} (\eta_{m+1}+\beta_{\Gamma_m}(B(v,r_m))^{2+\epsilon} \lesssim_\epsilon \sum_{v\in \widetilde{V}_0} \eta_{m+1}^{2+\epsilon} + \sum_{i=1}^m\sum_{v\in \widetilde V_i} \eta_i^{2+\epsilon} \lesssim_\epsilon 1+\sum_{i=1}^{m} 4^{i-1} \eta_i^{2+\epsilon}.\end{equation} Combining the previous two displayed equations, it follows that \begin{equation}S_{\Gamma,2+\epsilon}(\mathscr{G}') \lesssim_{A,\epsilon} \underbrace{\sum_{m=1}^\infty \log(m+15)4^{-m}}_I + \underbrace{\sum_{m=1}^\infty \log(m+15)4^{-m}\sum_{i=1}^m 4^{i-1}\eta_i^{2+\epsilon}}_{II}.\end{equation} It is apparent that $I\lesssim 1$. To bound $II$, exchange the order of summation: \begin{equation}\begin{split} II &= \sum_{i=1}^\infty 4^{i-1}\eta_i^{2+\epsilon}\sum_{m=i}^\infty \log(m+15)4^{-m} \lesssim \sum_{i=1}^\infty 4^{i-1}\eta_i^{2+\epsilon}\log(i+15)4^{-i}\\ &\lesssim \sum_{i=1}^\infty \log(i+15)\eta_i^{2+\epsilon} = \sum_{i=1}^\infty \log(i+15)\left(\frac{1}{4(i+15)\log(i+15)}\right)^{1+\frac12\epsilon}\lesssim_\epsilon 1.\end{split}\end{equation} We conclude that $S_{\Gamma,2+\epsilon}(\mathscr{G}')\lesssim_{A,\epsilon} 1$. Therefore, by our initial discussion, $S_{\Gamma,2+\epsilon}(\mathscr{G})<\infty$ for every $\epsilon>0$ and every multiresolution family $\mathscr{G}$ for $\Gamma$. This completes the proof of Proposition \ref{prop:2}.

\begin{remark}[importance of the logarithmic factor] \label{remark:log} Seeking out examples verifying the sharpness of \eqref{e:suff-p}, it is natural to first look at snowflake curves $\Gamma$ built with relative heights $\eta_i^2 \simeq 1/i$. By carrying out the outline above with relative heights $4\eta_i^2=\delta/(i+i_0)$ with parameters $\delta>0$ and $i_0\geq 1$, one obtains $$II \lesssim \sum_{i=1}^\infty (i+i_0)^\delta \left(\frac{1}{4(i+i_0)}\right)^{1+\frac12\epsilon},$$ where the latter expression is finite precisely when $\delta<(1/2)\epsilon$. Thus, for every $\epsilon>0$, we could build a curve $\Gamma_\epsilon$ with $\Haus^1(\Gamma_\epsilon)=\infty$ and $S_{\Gamma_\epsilon,2+\epsilon}(\mathscr{G})<\infty$ by selecting $\delta=\delta(\epsilon)$ sufficiently small. The logarithmic correction used in the proof of Proposition \ref{prop:2} allows us to find a single curve $\Gamma$ such that $\Haus^1(\Gamma)=\infty$ and $S_{\Gamma,2+\epsilon}(\mathscr{G})<\infty$ for all $\epsilon>0$. \end{remark}

\subsubsection{Proof of Proposition \ref{prop:4}}

Construct $\Gamma$ using Algorithm \ref{alg:2} with relative heights $4\eta_i^2 = 1/(i+2)\log(i+2)^2$ for all $i\geq 1$. Then $\eta_{i+1}\leq\eta_i< 1/\sqrt{12}$ for all $i\geq 1$ and \begin{equation} \sum_{i=1}^\infty 4\eta_i^2 = \sum_{n=3}^\infty \frac{1}{n\log(n)^2} = 1.069...<\infty. \end{equation} Hence $\Haus^1(\Gamma)\leq \exp(\sum_{i=1}^\infty 4\eta_i^2)<\exp(1.07)<3$ by Lemma \ref{l:gamma-length}. Similarly, we have that the intermediate curves $\Gamma_n$ consist of $4^n$ segments of length $r_n$, where \begin{equation} \label{p4-r-bound} 4^{-n} < r_n \leq \exp\left(\sum_{i=1}^n 4\eta_i^2\right)4^{-n} < 3 \cdot  4^{-n}\quad\text{for all $n\geq 1$}.\end{equation} Let $(X_k)_{k\in\ZZ}$ be an arbitrary family of nested $2^{-k}$-nets for $\Gamma$, let $A>1$, and let $\mathscr{G}=\{B(x,A2^{-k}):x\in X_k, k\in\ZZ\}$ be the associated multiresolution family for $\Gamma$. We wish to show that $S_{\Gamma,2-\epsilon}(\Gamma)=\infty$ for all $\epsilon>0$. By Lemma \ref{l:gamma-length}, $\Gamma$ is Ahlfors regular with constants determined by $\Haus^1(\Gamma)$. In particular, we know that $$\#X_k\simeq 2^k\quad\text{for every $k\geq 0$}.$$ Thus, writing $\beta(k):=\inf_{x\in X_k} \beta_{\Gamma}(B(x,A2^{-k}))$, we have \begin{equation}\label{S-lower-bound} S_{\Gamma,2-\epsilon}(\mathscr{G}) \geq \sum_{k=k_1}^\infty \sum_{x\in X_k} \beta_\Gamma(B(x,A2^{-k}))^{2-\epsilon}2A2^{-k}
\gtrsim_A \sum_{k=k_1}^\infty \beta(k)^{2-\epsilon}.\end{equation} To proceed, we will bound $\beta(k)$ from below in terms of $\eta_m$ for sufficiently large $k$.

Choose $k_0$ sufficiently large such that $A 2^{-k_0}<6$. Suppose that $k\geq k_0$. Let $m(k) \geq 1$ be the unique integer such that $6\cdot 4^{-m} \leq A 2^{-k} < 6\cdot 4^{-(m-1)}$. By \eqref{p4-r-bound}, we have \begin{equation}\label{p4-k-to-m} 2r_m \leq A 2^{-k} \lesssim r_m.\end{equation} Given $x\in X_k$, choose $v\in V_m$ such that $|x-v|<r_m$ (cf.~proof of Lemma \ref{l:gamma-length}). Then $B(v,r_m)\subset B(x,2r_m)\subset B(x,A2^{-k})$. Hence $$\beta_{\Gamma}(B(x,A2^{-k})) \geq \frac{r_m}{A2^{-k}}\beta_{\Gamma}(B(v,r_m)) \gtrsim \eta_m$$ by Lemma \ref{beta bound}. As $x\in X_k$ was arbitrary, $\beta(k)\gtrsim \eta_m$. Now, $m\leq \frac{1}{2}k+\log_4(6/A)+1$. Choose $k_1\geq k_0$ sufficiently large such that $\log_4(6/A)+3\leq \frac{1}{2}k_1$. Then, for every $k\geq k_1$, we have $m+2\leq k$ and $$\beta(k) \gtrsim \eta_m \gtrsim \frac{1}{(m+2)^{1/2}\log(m+2)}\gtrsim \frac{1}{k^{1/2}\log(k)}.$$ Therefore, for every $\epsilon>0$, $$S_{\Gamma,2-\epsilon}(\mathscr{G}) \gtrsim_A \sum_{k=k_1}^\infty \beta(k)^{2-\epsilon} \gtrsim_A \sum_{k=k_1}^\infty k^{\frac{1}{2}\epsilon-1} \log(k)^{\epsilon-2}=\infty.$$ This completes the proof of Proposition \ref{prop:4}.

\subsection{Examples with critical exponent \texorpdfstring{$p\neq 2$}{p != 2}}\label{ss:exponent-p}

To complete the proof of Theorems \ref{t:main} and \ref{t:main2}, we return to the infinite-dimensional setting. Our goal is to establish:

\begin{proposition}\label{prop:3} For all $1<p < \infty$, there is a curve $\Gamma$ in $\ell_p$ such that $\Haus^1(\Gamma)=\infty$ and $S_{\Gamma, p+\epsilon}(\mathscr{G})<\infty$ for every multiresolution family $\mathscr{G}$ for $\Gamma$ and every $\epsilon>0$. In particular, the exponent $p$ in \eqref{e:suff-p} when $1<p\leq 2$ is sharp.
\end{proposition}

\begin{proposition}\label{prop:1} For all $1< p<\infty$, there is a curve $\Gamma$ in $\ell_p$ such that  $\Haus^1(\Gamma)<\infty$ and  $S_{\Gamma,p-\epsilon}(\mathscr{G})= \infty$ for every multiresolution family $\mathscr{G}$ for $\Gamma$ and every $\epsilon>0$. In particular, the exponent $p$ in \eqref{e:nec-2} when $2\leq p<\infty$ is sharp.
\end{proposition}

We construct the curves in both propositions using the following algorithm, which is inspired by examples of Edelen, Naber, and Valtorta \cite[\S5.2]{ENV-Banach} in $L^p([0,1])$. The key point is that because we are working in an infinite dimensional space, we may build each intermediate iteration of the snowflake by adding bumps in a new coordinate direction.

\begin{algorithm}[snowflake-like curves in $\ell_p$ with bumps along coordinate directions]\label{alg:p}
Let $\{e_i\}_{i=1}^\infty$ denote the standard basis in $\ell_p$, i.e.~ $e_i(j)=\delta_{ij}$. Suppose  $x,y\in\Span\{e_1,\dots,e_k\}$ and $\gamma:[a,b]\rightarrow\Span\{e_1,\dots,e_k\}\cong \ell_p^k$ is a constant speed parameterization of $I=[x,y]$ from $\gamma(a)=x$ to $\gamma(b)=y$. That is, $$\gamma(t)=x + \frac{t-a}{b-a}(y-x)\quad\text{for all }t\in[a,b].$$ Write $v:=y-x$. Given a \emph{relative height} $0 \leq  \eta <1/2$, we define a piecewise linear path $\hat\gamma:[a,b]\rightarrow\Span\{e_1,\dots,e_{k+1}\}\cong \ell_p^{k+1}$, as follows. Define $s\in[1/4,1/2)$ to be the unique solution of $s^p=(\frac12-s)^p+\eta^p$. Divide $[a,b]$ into quarters. For all $a\leq t \leq (3/4)a+(1/4)b$,
    $$\hat\gamma(t) = x + 4\left(\frac{t-a}{b-a}\right)s v$$ (cf.~Algorithm \ref{alg:2}).
For all $(3/4)a + (1/4)b \leq t \leq (1/2)a+(1/2)b$,
    $$\hat\gamma(t) = x + sv + 4\left(\frac{t-(3/4)a-(1/4)b}{b-a}\right) \big((1/2-s)v + \eta|v|_pe_{k+1}\big).$$
For all $(1/2)a+(1/2)b\leq t\leq (1/4)a + (3/4)b$,
    $$\hat\gamma(t) = x + (1/2)v + \eta|v|e_{k+1} + 4\left(\frac{t-(1/2)a-(1/2)b}{b-a}\right)\big((1/2-s)v-\eta|v|_pe_{k+1}\big).$$
Finally, for all $(1/4)a+(3/4)b\leq t\leq b$,
    $$\hat\gamma(t) = x + (1-s)v + 4\left(\frac{t-(1/4)a-(3/4)b}{b-a}\right)sv.$$ We say that $\hat \gamma$ is obtained from $\gamma$ by \emph{adding a bump of relative height $\eta$ in the direction $e_{k+1}$}. On each quarter of $[a,b]$, $\hat \gamma(t)$ traces a line segment in $\ell_p$ of length $s|v|_p$ at constant speed. Thus, $\hat\gamma$ has Lipschitz constant $4s|v|_p/(b-a)$. Additionally, we have $\|\gamma-\hat\gamma\|_\infty =|\gamma((a+b)/2)-\hat\gamma((a+b)/2)|_p=\eta|v|_p$.

Starting from the arc length parameterization $\gamma_0:[0,1]\rightarrow\Span\{e_1\}$ of the line segment $I_0=[0,e_1]$ from $\gamma_0(0)=0$ to $\gamma_0(1)=e_1$, we now define a sequence of piecewise linear maps $\gamma_i:[0,1]\rightarrow\Span\{e_1,\dots,e_{i+1}\}$  by iteratively adding bumps of relative height $\eta_i$ in the direction $e_{i+1}$. Suppose that $\gamma_i$ has been defined for some $i\geq 0$ so that $\gamma_i|_{J_{i,k}}$ is a constant speed parameterization of a line segment in $\Span\{e_1,\dots,e_{i+1}\}$ on each interval $J_{i,k}=[a_{i,k},b_{i,k}]$ of the form $$J_{i,k}=\left[a+\frac{(k-1)}{4^i}(b-a), a + \frac{k}{4^i}(b-a)\right]\qquad (1\leq k \leq 4^i,\  k\in\ZZ).$$ For each index $1\leq k\leq 4^i$, define $\gamma_{i+1}|_{J_{i,k}} = \widehat{\gamma_i|_{J_{i,k}}}$ by adding a bump of relative height $\eta_{i+1}$ in the direction $e_{i+2}$. This defines a map $\gamma_{i+1}$. By induction, we obtain the full sequence $\gamma_0, \gamma_1,\gamma_2,\dots$; the image of each map $\gamma_n$ is composed of $4^n$ line segments in $\Span\{e_1,\dots,e_{n+1}\}$ of equal length $r_n=(\Lip\gamma_n)4^{-n}=s_1\cdots s_n$, where $s_i\in[1/4,1/2)$ denotes the solution to $s_i^p=(\frac12-s_i)^p+\eta_i^p$. Moreover, \begin{align}\label{gamma-lip-p} \Lip (\gamma_n) &= \prod_{i=1}^n 4s_i<2^n, \\ \label{r-comparison-p} r_{n+1} &\leq (1/2) r_n, \\ \label{gamma-difference-p} \|\gamma_n-\gamma_{n+1}\|_\infty &\leq \eta_{n+1} r_n\leq 2^{-n}\eta_{n+1}<2^{-n}.\end{align} Therefore, $\gamma:[0,1]\rightarrow\ell_p$, which is defined pointwise by \begin{equation} \label{gamma-definition-p} \gamma(t)=\gamma_0(t)+\sum_{n=0}^\infty(\gamma_{n+1}(t)-\gamma_n(t))\quad\text{for all }t\in[0,1],\end{equation} is continuous as the uniform limit of the maps $\gamma_n$ by \eqref{gamma-difference-p}. We denote $\gamma([0,1])$ by $\Gamma$, and for each $n$, we denote $\gamma_n([0,1])$ by $\Gamma_n$. For all integers $n,k\geq 0$ with $0\leq k\leq 4^n$, we call the point $\gamma_n(k/4^n)$ a \emph{vertex} of $\Gamma_n$. Remark \ref{vertex facts} holds in this setting, as well. \end{algorithm}

\begin{remark}[modulus of continuity and injectivity] \label{mod-injective-p} In fact, \eqref{gamma-lip-p} and \eqref{gamma-difference-p} imply that the maps $\gamma_n$ and $\gamma$ are uniformly $(1/2)$-H\"older continuous. The parameterization $\gamma$ is injective and this can be verified by showing that $\pi_{1}\circ \gamma$ is injective, where $\pi_{1}:\ell_p\rightarrow\RR$ is projection onto the first coordinate. We leave details for the reader. \end{remark}

\begin{remark}[improved bounds on $r_n$] \label{r:good segment estimates} Let $\Gamma$ in $\ell_p$ ($1<p<\infty$) be constructed by Algorithm \ref{alg:p} with relative heights $0\leq \eta_i\leq \overline{\eta}$ for all $i$, for some universal constant $\overline\eta$ to be specified below. Let $s_i\in[1/4,1/2)$ be defined by $s_i^p=(\frac12-s_i)^p+\eta_i^p$. Then $$s_i\leq  \left((1/4)^p+\eta_i^p\right)^{1/p} = \frac{1}{4}\left(1+(4\eta_i)^p\right)^{1/p}.$$ Since $p>1$, Taylor's theorem with remainder gives $$ 1+\frac{1}{p}\delta-\frac{p-1}{2p^2}\delta^2\leq (1+\delta)^{1/p} \leq 1+\frac{1}{p}\delta\quad\text{for all }0\leq \delta<1.$$ On the one hand, assume that $\overline{\eta}<1/4$. Then $(4\eta_i)^p<4\overline{\eta}<1$ and we obtain $$s_i \leq \frac{1}{4}\left(1+\frac{1}{p}(4\eta_i)^p\right)=\frac{1}{4}+\frac{1}{4p}(4\eta_i)^p.$$  On the other hand, assume that $\overline{\eta}<1/8$. Then $s_i\leq1/4+(1/4p)2^{-p}<3/8$. Hence $1/8<\frac12-s_i\leq 1/4$ and using the Taylor bound we can write \begin{align*}s_i &= \left(\frac12-s_i\right)\left(1+\left(\frac12-s_i\right)^{-p}\eta_i^p\right)^{1/p} \\
&\geq \frac12 -s_i + \frac{1}{p}\left(\frac12-s_i\right)^{1-p}\eta_i^p - \frac{p-1}{2p^2}\left(\frac12-s_i\right)^{1-2p}\eta_i^{2p} \\
&\geq \frac12-s_i+\frac{1}{4p}(4\eta_i)^p-\frac{p-1}{16p^2}(8\eta_i)^{2p}.
\end{align*} Rearranging the inequality, we obtain $$s_i\geq \frac{1}{4}+\frac{1}{8p}(4\eta_i)^p-\frac{p-1}{32p^2}(8\eta_i)^{2p}=\frac{1}{4}+\left[\frac{1}{8p}-\frac{p-1}{32p^2}(16\eta_i)^p\right](4\eta_i)^p.$$ We now specify that $\overline{\eta}=1/16$ so that $$\frac{1}{8p}-\frac{p-1}{32p^2}(16\eta_i)^p\geq\frac{1}{8p} -\frac{p-1}{32p^2}\geq \frac{3}{32p}.$$ Therefore, \begin{equation} \label{sp-bound} \frac{1}{4}+\frac{3}{32p}(4\eta_i)^p \leq s_i \leq \frac{1}{4}+\frac{1}{4p}(4\eta_i)^p\quad\text{whenever }\eta_i\leq \frac{1}{16}.\end{equation} Thus, if $\eta_i\leq 1/16$ for all $i$, then the length $r_n=\prod_{i=1}^n s_i$ of each segment in $\Gamma_n$ satisfies \begin{equation}\label{rp-bound} 4^{-n}\prod_{i=1}^n \left(1+\frac{3}{8p}(4\eta_i)^p\right) \leq r_n \leq 4^{-n}\prod_{i=1}^n \left(1+\frac{1}{p}(4\eta_i)^p\right).\end{equation} The reader may verify that when $p=2$, the bound \eqref{rp-bound} is compatible with \eqref{gamma-lip}.
\end{remark}

\begin{lemma}\label{l:gamma-length-p} If $\Gamma$ is constructed by Algorithm \ref{alg:p} with relative heights $\eta_i\leq 1/16$, then \begin{equation}\label{length-exp-p} \exp\left(\frac1{4p}\textstyle\sum_{n=1}^\infty (4\eta_n)^p\right) \leq \Haus^1(\Gamma)\leq \exp\left(\frac{1}{p}\textstyle\sum_{n=1}^\infty (4\eta_n)^p\right).\end{equation} In particular, $\Gamma$ is a rectifiable curve if and only if $\sum_{n=1}^\infty \eta_n^p<\infty$. Moreover, in that case, $\Gamma$ is Ahlfors regular with constants depending only on $\Haus^1(\Gamma)$.\end{lemma}

\begin{proof} The outline is similar to the proof of Lemma \ref{l:gamma-length}. For each $n\geq 1$, \begin{align*}\Haus^1(\Gamma) &\stackrel{\text{Rk \ref{mod-injective-p}}}{=}  \sum_{k=1}^{4^n} \Haus^1(\gamma([(k-1)4^{-n},k4^{-n}])) \geq \sum_{k=1}^{4^n} \diam \gamma([(k-1)4^{-n},k4^{-n}]) \\
&\stackrel{\eqref{rp-bound}}{\geq} \sum_{k=1}^{4^n} 4^{-n} \prod_{i=1}^n \left(1+\frac38p^{-1}(4\eta_i)^p\right)= \prod_{i=1}^n \left(1+\frac38p^{-1}(4\eta_i)^p\right).
\end{align*} We conclude that $$ \prod_{n=1}^\infty \left(1+\frac{3}{8p}(4\eta_n)^p\right)\leq \Haus^1(\Gamma)\leq \Lip\gamma\leq \liminf_{m\rightarrow\infty}\Lip \gamma_m \stackrel{\eqref{rp-bound}}{\leq} \prod_{n=1}^\infty \left(1+\frac{1}{p}(4\eta_n)^p\right).$$ To derive \eqref{length-exp-p}, rewrite each infinite product as the exponential of an infinite sum and use Taylor's theorem bounds for $\log(1+x)$ with $0\leq x\leq 1/4$.

Suppose that $\Haus^1(\Gamma)<\infty$. By \eqref{r-comparison-p}, \eqref{gamma-difference-p}, and assumption $\eta_n\leq 1/16$ for all $n\geq 1$, $$\|\gamma_n-\gamma\|_\infty \leq \sum_{j=n}^\infty \eta_{j+1}r_j \leq \frac{1}{16}\sum_{k=0}^\infty r_{n+k} \leq  \frac{1}{16}r_n\sum_{k=0}^\infty 2^{-k}=\frac{1}{8}r_n.
$$ Thus, given $x\in\Gamma$ and $n\geq 1$, we may pick $y\in \Gamma_n$ such that $|x-y|\leq (1/8)r_n$. Next, let $v\in V_n$ be an endpoint of the segment in $\Gamma_n$ containing $y$ that is closest to $y$ so that $|v-y|\leq (1/2)r_n$. Hence $B(x, (3/8)r_n)\subset B(v,r_n)$ and we may use \eqref{length-exp} to estimate \begin{equation}\label{gamma-AR-p} \Haus^1(\Gamma\cap B(x,\frac{3}{8}r_n)) \leq \Haus^1(\Gamma\cap B(v,r_n)) \leq 2 \exp\left(\sum_{i=n}^\infty \frac{1}{p}(4\eta_i)^p\right)r_n \leq 2 e^{4} \Haus^1(\Gamma)r_n,\end{equation} because $\Gamma_n\cap B(v,r_n)$ consists of one or two line segments of length $r_n$. It follows that $\Gamma$ is Ahlfors regular with constants depending only on $\Haus^1(\Gamma)$.
\end{proof}

\begin{lemma}\label{beta bound-p} Assume that $\Gamma$ is constructed by Algorithm \ref{alg:p} with $\eta_i<1/8$ for all $i$. For all  $i\geq 0$ and $j\geq i$, \begin{equation}\label{beta-p1}\frac{1}{4}\eta_i\leq \beta_{\Gamma_j}(B(v,r_j)) \leq 2\eta_i\quad\text{for all }v\in \widetilde{V}_i,\end{equation} where $\widetilde{V}_i$ denote the ``new vertices" in $\Gamma_i$ (see Remark \ref{vertex facts}) and $\eta_0=0$. Furthermore, \begin{equation}\label{beta-p2} \frac{1}{4}\eta_i \leq \beta_{\Gamma}(\overline{B(v,r_i)})\quad\text{for all $i\geq 0$ and }v\in\widetilde{V}_i.\end{equation}
\end{lemma}

\begin{proof} The upper bound agrees with the case $p=2$ above. Note that when $v\in \widetilde{V}_0$ is an endpoint of $\Gamma_j$, the set $\Gamma_j\cap B(v,r_j)$ is a line segment. Thus, $\beta_{\Gamma_j}(B(v,r_j))=0$ for all $v\in\widetilde{V}_0$ and $j\geq 0$. Next, suppose that $i\geq 1$ and $v\in \widetilde{V}_i$. Using the line containing the segment in $\Gamma_{i-1}\cap B(v,r_i)$ to approximate the beta number, we obtain the estimate $$\beta_{\Gamma_i}(B(v,r_i)) \leq \frac{\eta_ir_{i-1}}{2r_i}=\frac{\eta_i}{2s_i}\leq 2\eta_i,$$ where $\sup_{x\in \Gamma_i}\dist(x,\ell)\leq \eta_ir_{i-1}$, because at stage $i$ in the construction we added the bump in the direction $e_{i+1}$ and we recall that $s_i\geq \frac14$. When $j\geq i$, the set $\Gamma_j\cap B(v,r_j)$ agrees up to a dilation centered at $v$ with $\Gamma_i\cap B(v,r_i)$. Thus, $\beta_{\Gamma_j}(B(v,r_j)) \leq 2\eta_i$ for all $v\in\widetilde{V}_i$ and $j\geq i$, as well.

Before we determine the lower bounds in $(\ref{beta-p1})$ and $(\ref{beta-p2})$, we recall two basic properties of the beta numbers. First, for any set $E \subset B(v, r)$, we have $\beta_{E}(B(v,r)) = \beta_{\overline{E}}(\overline{B(v,r)})$. (This may fail if $E \not \subset B(v, r)$!)  Second, $\beta_E(B(v,r))$ is increasing in $E$. Because $V_i \subset \Gamma_j$ for all $j \ge i$ and $V_i \subset \Gamma$, for any $v_k \in V_i$ the vertices $v_{k-1}$ and $v_{k+1}$ that are adjacent with respect to the global parametrization satisfy
$$
\{v_{k-1}, v_k, v_{k+1}\} \subset \overline{\Gamma_j \cap B(v_k, r_i)}\quad\text{for any }j\geq i
$$
and the inclusion also holds with $\Gamma$ in place of $\Gamma_j$.  Therefore, up to a translation and dilation, there are two relevant configurations of vertices: $$E_v(\eta):=\{0\} \cup \{sv\} \cup \{\tfrac12v+\eta e_{n+1}\},\quad F_v(\eta):=\{sv\} \cup \{\tfrac12v+\eta e_{n+1}\} \cup \{(1-s)v\},$$ where $v\in\Span\{e_1,\dots,e_n\}$ is an arbitrary vector with $|v|_p=1$ and $s^p=(\frac12-s)^p+\eta^p$ with $\eta<1/8$. By Remark \ref{r:good segment estimates}, $1/4 \leq s < 3/8$. The optimal lower bound in \eqref{beta-p1} and \eqref{beta-p2} is given by $\beta(\eta_i)$, where $$\beta(\eta):=\min\left\{\beta_{E_v(\eta)}(\overline{B(sv,s)}),\beta_{F_v(\eta)}(\overline{B(\tfrac12v+\eta e_{n+1},s)})\right\}.$$

We work separately for each configuration, beginning with $E_v(\eta).$ Let $L_1$ be the line containing $[0, \tfrac12v+\eta{e_{n+1}}]$. In the $v e_{n+1}$-plane, the line $L_1$ is the locus of points $(x, y)$ satisfying $y = 2\eta x$.
When $x \ge \frac{1}{4} - \frac{1}{2}\eta$,
\begin{equation*}y \ge 2\eta \left(\frac{1}{4} - \frac{1}{2}\eta\right) \ge \frac{1}{2}\eta - \frac{1}{16}\eta = \frac{7}{16}\eta.\end{equation*}
Since $s \ge \frac{1}{4}$, it follows that $B(sv, (7/16)\eta) \cap L_1 = \emptyset.$ Hence \begin{equation}\label{sv-far}\dist(sv, L_1) \ge (7/16)\eta.\end{equation} We now argue that $\beta_{E_v(\eta)}(\overline{B(sv,s)}) \ge (3/32s)\eta$ by way of contradiction.  Assume that we can find a line $L$ such that $E_v(\eta) \subset B_{(3/16)\eta}(L).$  By convexity, $L_1 \cap B(sv, s) \subset B_{(3/16)\eta}(L)$, as well. Let $x_1 \in L \cap B(0, (3/16)\eta)$ and $x_2  \in L \cap B(\frac{1}{2}v + \eta e_{n+1},(3/16)\eta).$  Once again, by convexity, the segment of $L$ that falls between $x_1$ and $x_2$ is contained in $B_{(3/16)\eta}(L_1).$  Since $\dist(x,sv)\ge \frac{1}{2}s \geq (3/16)\eta$ for all $x \in L \setminus B_{(3/16)\eta}(L_1)$, the points $x \in L$ such that $\dist(x, sv) \le (3/16)\eta$ are contained in $B_{(3/16)\eta}(L_1).$  Thus, by the triangle inequality,
\begin{equation}\label{Ev-near}
E_v(\eta) \subset B_{(6/16)\eta}(L_1).
\end{equation} Since \eqref{sv-far} and \eqref{Ev-near} are incompatible, we have reached a contradiction. Therefore, $$\beta_{E_v(\eta)}(\overline{B(sv,s)}) \ge \frac{3}{32s}\eta \ge \frac{1}{4}\eta.$$

Using symmetry, one sees the best fitting line for $F_v(\eta)$ is $\frac{1}{2}\eta e_{n+1}+\Span(v)$, whence $$\beta_{F_v(\eta)}(\overline{B(\tfrac12v+\eta e_{n+1},s)})\geq \frac{\eta}{4s} \geq \frac{2}{3}\eta.$$ Thus, the extremal configuration is given by $E_v(\eta)$ and $\beta(\eta)\geq (1/4)\eta$.
\end{proof}

\subsubsection{Proof of Proposition \ref{prop:3}}

Let $1<p<\infty$. Let $\Gamma$ be the curve constructed by Algorithm \ref{alg:p} with relative heights $\eta_i^p=\delta/(i+i_0)\log(i+i_0)$ for all $i\geq 1$, where $\delta>0$ and $i_0$ are chosen so that $\eta_1\leq 1/16$. Then $\Haus^1(\Gamma)=\infty$ in $\ell_p$ by Lemma \ref{l:gamma-length-p}. To prove that $S_{\Gamma,p+\epsilon}(\mathscr{G})<\infty$ for every multiresolution family $\mathscr{G}$ for $\Gamma$ and $\epsilon>0$, repeat the proof of Proposition \ref{prop:2}, \emph{mutatis mutandis}, using Lemma \ref{beta bound-p} \emph{in lieu} of Lemma \ref{beta bound}.

\subsubsection{Proof of Proposition \ref{prop:1}} Let $1<p<\infty$. Let $\Gamma$ be the curve constructed by Algorithm \ref{alg:p} with relative heights $\eta_i^p=\delta/(i+i_0)\log(i+i_0)^2$ for all $i\geq 1$, where $\delta>0$ and $i_0$ are chosen so that $\eta_1\leq 1/16$. Then $\Haus^1(\Gamma)<\infty$ in $\ell_p$ and $\Gamma$ is Ahlfors regular by Lemma \ref{l:gamma-length-p}. To prove that $S_{\Gamma,p-\epsilon}(\mathscr{G})=\infty$ for every multiresolution family $\mathscr{G}$ for $\Gamma$ and $\epsilon>0$, repeat the proof of Proposition \ref{prop:4}, \emph{mutatis mutandis}, again substituting Lemma \ref{beta bound-p} for Lemma \ref{beta bound}.

\subsection{Proof of Proposition \ref{t:pq}}\label{ss:pq}

Let $1< p < q < \infty$.  To begin, we verify that if $\Gamma\subset\ell_p$ is a rectifiable curve, then $\Gamma$ is also rectifiable when viewed as a curve embedded in $\ell_q$. As is well known, $\ell_p\subset \ell_q$ and $|v|_q \leq |v|_p$ for every $v \in \ell_p$. Hence the diameter of a set in $\ell_p$ does not increase when embedded into $\ell_q$.  Thus, $\mathcal{H}^s_{\ell_q}(E) \le \mathcal{H}^s_{\ell_p}(E)$ for every $s>0$ and $E\subset\ell_p$, where $\mathcal{H}^s_{\ell_r}$ denotes the $s$-dimensional Hausdorff measure in $\ell_r$. In particular, every rectifiable curve in $\ell_p$ is also a rectifiable curve in $\ell_q$, possibly with shorter length.

We now construct a curve $\Gamma$ such that $\mathcal{H}^{1}_{\ell_p}(\Gamma) = \infty$ and $\mathcal{H}^{1}_{\ell_q}(\Gamma) < \infty$ for every $q>p$. Build $\Gamma$ in $\ell_p$ using Algorithm \ref{alg:p} with relative heights $$\eta_i^p = \frac{\delta}{i \log(i+i_0)}\quad\text{for all }i\geq 1,$$ where $\delta>0$ and $i_0\in\N$ are chosen so that $0<\eta_1 \le 1/16$. Note that $\sum_{i =1}^{\infty} \eta_i^p = \infty$.
Therefore, $\mathcal{H}^{1}_{\ell_p}(\Gamma) = \infty$ by Lemma \ref{l:gamma-length-p}.

Fix an exponent $q$ with $p<q<\infty$. We break the proof that $\mathcal{H}^{1}_{\ell_q}(\Gamma) < \infty$ into a series of lemmata. First, we calculate the $\mathcal{H}^{1}_{\ell_q}$-growth of a line segment under the snowflaking procedure $\gamma\mapsto\hat\gamma$ used in Algorithm \ref{alg:p}. We emphasize that the estimate in the following lemma is independent of $n$.

\begin{lemma}\label{l:growth factor}
Given $v \in \Span\{e_1, ..., e_n\}$ with $|v|_p = 1$, let $\gamma: [0, 1]\rightarrow \Span\{e_1, ..., e_n \}$ be the unit speed parameterization of the line segment $I = [0, v]$ in $\ell_p$. Let $\hat{\gamma}$ be obtained from $\gamma$ by adding a bump of relative height $\eta$ in the direction $e_{n+1}$ (see Algorithm \ref{alg:p}). If $\eta\leq 1/16$, then $\mathcal{H}^1_{\ell_q}(I)=|v|_q$ and
\begin{equation}\label{e:growth factor}
\mathcal{H}^{1}_{\ell_q}(\hat{\gamma}([0,1])) - \mathcal{H}^{1}_{\ell_q}(I) \lesssim_q \left(\frac{\eta}{|v|_q}\right)^q|v|_q.
\end{equation}
\end{lemma}

\begin{proof}
The estimate is nearly identical to the calculation in Remark \ref{r:good segment estimates}. Since $\eta\leq 1/16$, each of the four edges of $\hat\gamma(I)$ in $\ell_p$ have length $s \in [1/4, 3/8)$ given implicitly by the equation $s^p = (\frac{1}{2} - s)^p + \eta^p$. In $\ell_q$, $\Haus^1_{\ell_q}(I)=|v|_q$ (since $I$ is a segment) and the four edges of $\hat\gamma(I)$ have total length
\begin{align*}
\Haus^{1}_{\ell_q}(\hat{\gamma}([0,1]))& = 2s |v|_q + 2\left((\tfrac{1}{2} - s)^q|v|_q^q + \eta^q \right)^{1/q}\\
&= 2s|v|_q + (1-2s)|v|_q\left(1+(\tfrac12-s)^{-q}|v|_q^{-q}\eta^q\right)^{1/q}\\
&\leq |v|_q\left(1+(\tfrac12-s)^{-q}|v|_q^{-q}\eta^q\right)^{1/q}
\end{align*} Applying $(1+\delta)^{1/q}\leq 1+\frac{1}{q}\delta$ for all $\delta\geq 0$ with $\delta=(\tfrac12-s)^{-q}|v|_q^{-q}\eta^q$, we conclude
\begin{align*}
 \Haus^1_{\ell_q}(\hat{\gamma}([0,1]))&\le |v|_q+\frac{1}{q}\left(\frac{\eta}{(\tfrac12-s)|v|_q}\right)^q|v|_q \leq |v|_q+\frac{1}{q}\left(\frac{8\eta}{|v|_q}\right)^q|v|_q. \qedhere
\end{align*}
\end{proof}

The next estimate is elementary and left as an exercise for the reader.

\begin{lemma}\label{l:finite-p-to-q} For all $n\in\N$ and $v\in\ell_p^n$, $|v|_q \geq n^{\frac{1}{q}-\frac{1}{p}} |v|_p$.
\end{lemma}

We may now give a uniform estimate on the length of the intermediate curves $\Gamma_n$ approximating $\Gamma$ in $\ell_q$.

\begin{lemma}\label{l:approximation length bound}
With the exponents and the parameter $\eta_i^p = \delta/(i\log(i+i_0))$ fixed as above, the intermediate curves produced by Algorithm \ref{alg:p} satisfy
\begin{align*}
\mathcal{H}^{1}_{\ell_q}(\Gamma_n) \leq C(p, q)<\infty\quad\text{for all }n\geq 1,
\end{align*} where $C(p,q)$ is a constant depending on $p$ and $q$ with $C(p,q)\uparrow\infty$ as $q\downarrow p$.
\end{lemma}

\begin{proof}
The initial segment $\Gamma_0 = [0,e_1]$ satisfies $\mathcal{H}^{1}_{\ell_q}(\Gamma_0) = 1$. Let $n\geq 0$, let $1\leq k\leq 4^n$, and let $J=\gamma_n([(k-1)4^{-n},k4^{-n}])$ be an edge in $\Gamma_n$. Assign $J^+:=\gamma_{n+1}([(k-1)4^{-n},k4^{-n}])$. By Lemma \ref{l:growth factor} and Lemma \ref{l:finite-p-to-q}, $$\mathcal{H}^1_{\ell_q}(J^+) - \mathcal{H}^1_{\ell_q}(J) \lesssim_q \left(\frac{\eta_{n+1}\Haus^1_{\ell_p}(J)}{\Haus^1_{\ell_q}(J)}\right)^q\Haus^1_{\ell_q}(J) \lesssim_q \left(\frac{\eta_{n+1}}{(n+1)^{\frac{1}{q}-\frac{1}{p}}}\right)^q\Haus^1_{\ell_q}(J).$$ Summing over the $4^n$ edges $J$ in $\Gamma_n$, it follows that $$\mathcal{H}^1_{\ell_q}(\Gamma_{n+1}) \leq \left(1+C_q\left(\frac{\eta_{n+1}}{(n+1)^{\frac{1}{q}-\frac{1}{p}}}\right)^q\right)\Haus^1_{\ell_q}(\Gamma_n),$$ where $C_q$ is a positive constant depending only on $q$.
Therefore, $$\mathcal{H}^{1}_{\ell_q}(\Gamma_n) \le \prod_{i=1}^n \left( 1 + C_q \left(\frac{\eta_i}{i^{\frac{1}{q}-\frac{1}{p}}}\right)^q \right) \leq \exp \left(C_q \sum_{i=1}^n \left(\frac{\eta_i}{i^{\frac{1}{q}-\frac{1}{p}}}\right)^q\right).$$ Finally, recalling our choice of $\eta_i$,
\begin{align*}
\sum_{i=1}^{n}\left(\frac{\eta_i}{i^{\frac{1}{q} - \frac{1}{p}}}\right)^q & \le \sum_{i=1}^{\infty}\left(\frac{\delta^{1/p} i^{-1/p}\log(i+i_0)^{-1/p}}{i^{\frac{1}{q} - \frac{1}{p}}}\right)^q = \sum_{i=1}^{\infty}\frac{\delta^{q/p}}{i \log(i+i_0)^{q/p}} <\infty,\end{align*} because $q/p>1$. All together, $\mathcal{H}^{1}_{\ell_q}(\Gamma_n) \leq C(p,q)<\infty$, where $C(p,q)\uparrow\infty$ when $q\downarrow p$.
\end{proof}

By \eqref{gamma-difference-p}, $\gamma_n$ converges uniformly to $\gamma$ in $\ell_p$, and thus, in $\ell_q$. Hence $\Gamma_n$ converges to $\Gamma$ in the Hausdorff metric on compact subsets of $\ell_q$. Therefore, by Go\l\c{a}b's semicontinuity theorem (see e.g.~\cite[Theorem 2.9]{AO-curves}), $\mathcal{H}^{1}_{\ell_q}(\Gamma) \le \liminf_{n \rightarrow \infty}\mathcal{H}^{1}_{\ell_q}(\Gamma_n) \leq C(p, q)<\infty$ for all $q>p$. This completes the proof of Proposition \ref{t:pq}.

\appendix

\section{Schul's prefiltration lemma in a metric space} \label{ss:prefiltration}

Before we start the proof of Lemma \ref{build-filtrations}, we record proofs of a couple of metric lemmas. Recall that a \emph{pseudometric} on $\XX$ is a function $d:\XX\times\XX\rightarrow[0,\infty)$ such that $d(x,x)=0$, $d(x,y)=d(y,x)$, and $d(x,z)\leq d(x,y)+d(y,z)$ for all $x,y,z\in\XX$.

\begin{lemma} \label{curve-pseudometric} Let $\XX$ be a metric space. If $f:[0,1]\rightarrow\overline{\Sigma}$ is a continuous parameterization of a set $\overline{\Sigma}\subset\XX$, then $d(a,b)=\diam f([a,b])$ is a pseudometric on $[0,1]$. Furthermore, with respect to this pseudometric, if $0\leq a\leq b\leq x\leq y\leq 1$, then $\diam_d([a,b])=d(a,b)$ and $\gap_d([a,b],[x,y])=d(b,x)$.\end{lemma}

\begin{proof} To verify that $d$ is a pseudometric on $[0,1]$, only the triangle inequality requires some thought. Given $a,b,c\in[0,1]$, set $$A=f([a,b]),\quad B=f([b,c]),\quad C=f([a,c]),$$ so that $\diam A=d(a,b)$, $\diam B=d(b,c)$, and $\diam C=d(a,c)$. Up to relabeling, there are two cases.

\emph{Case 1.} If $a\leq c\leq b$, then $C\subset A$ and $d(a,c)\leq d(a,b)\leq d(a,b)+d(b,c)$.

\emph{Case 2.} Suppose $a\leq b\leq c$. We now argue as in \cite[Lemma 2.18]{AO-curves}, a step in a proof of Go\l\c{a}b's semicontinuity theorem. Let $p$ and $q$ be points in $C$ such that $\dist(p,q)=d(a,c)$. Consider the 1-Lipschitz function $F:C\rightarrow\RR$ defined by $F(x)=\dist(x,p)$. Then \begin{align*}d(a,b)+d(b,c) &\geq \diam F(A)+\diam F(B) \geq \diam F(C)\geq d(a,c),\end{align*} where the first inequality holds since $\Lip(F)=1$, the second inequality holds since $F(A),F(B),F(C)$ are intervals in $\RR$ with $F(C)\subset F(A)\cup F(B)$, and the third inequality holds since $F(p)=0$ and $F(q)=\dist(p,q)=d(a,c)$.

Suppose that $I=[a,b]\subset[0,1]$. On the one hand, for any $a\leq a'\leq b'\leq b$, we know that the distance $d(a',b')=\diam f([a',b'])\leq \diam f([a,b]) = d(a,b)$, because $[a',b']\subset[a,b]$. Hence $\diam_d I \leq d(a,b)$. On the other hand, $\diam_d I\geq d(a,b)$, because $a,b\in I$.
A similar argument yields the formula for the gap between intervals.\end{proof}

\begin{remark}If the parameterization $f:[0,1]\rightarrow\overline{\Sigma}$ has positive constant speed, then the pseudometric $d$ on $[0,1]$ is a metric. The pseudometric $d$ on $[0,1]$ bears some similarity to the intrinsic diameter distance on $\overline{\Sigma}$; see e.g.~\cite{Freeman-Herron, Herron-1}.\end{remark}

\begin{example}\label{bad-gap} Let $\{e_i\}_{i=1}^\infty$ denote the standard basis in $\ell_2$. Fix a large integer $n>1$ and let $f:[0,1]\rightarrow\ell_2$ be a continuous map, which traces a constant speed, piecewise linear path from $e_1$ to $e_2$ to $e_3$ and so on until reaching $e_{n+1}$. Decompose $[0,1]$ into $n$ intervals $I_i=f^{-1}([e_i,e_{i+1}])$. With respect to the metric $d$ of Lemma \ref{curve-pseudometric}, we have $\diam_d I_i=\diam f(I_i)=\sqrt{2}$ for all $i$ \emph{and} $\diam_d [0,1]=\diam f([0,1])=\sqrt{2}$. This shows that $\diam_d$ is not superadditive. Similarly, even though there are $\gtrsim n/2$ disjoint intervals of diameter $\sqrt{2}$ lying between $I_1$ and $I_n$, we only have $\gap_d(I_1,I_n)=|e_2-e_n|_{\ell_2}=\sqrt{2}$.\end{example}

The second auxiliary lemma is related to \cite[Lemma 2.16]{Schul-AR}, but with an alternate gap condition that is motivated by the previous example.

\begin{lemma}\label{interval-lemma} Let $\XX$ be a metric space, let $f:[0,1]\rightarrow\overline{\Sigma}$ be a continuous parameterization of a set $\overline{\Sigma}\subset\XX$, let $\xi>6$, and let $0<r_-< r_+<\infty$. Suppose $\{J_i\}_{i=1}^{I}$ is a finite ($I<\infty$) or infinite ($I=\infty$) sequence of closed intervals in $[0,1]$ and $(k_i)_{i=1}^I$ is a sequence of integers bounded from below such that \begin{enumerate}
\item chain property: for all $1\leq i <I$, we have $J_i\cap J_{i+1}\neq\emptyset$;
\item geometric decay: for all $i\geq 1$, we have $\diam f(J_i)\leq \xi^{-k_i}r_+$; and
\item separation within levels: for all $i,j\geq 1$ with $i\neq j$, if $k_i=k_j=k$, then there exists at least $\lceil 3r_+/r_-\rceil$ pairwise disjoint closed intervals $K$ with $\diam f(K)\geq \xi^{-k}r_-$ and such that $K$ lies between $J_i$ and $J_j$ as subsets of $\RR$.
\end{enumerate} If $\xi>r_+/r_-$, then there exists a unique $M\geq 1$ such that $k_M=\min_{i\geq 1} k_i$; moreover, \begin{equation}\label{e:interval-growth} \sum_{i=1}^I \diam f(J_i) \leq (1+3/\xi)\xi^{-k_M}r_+.\end{equation}
\end{lemma}

\begin{proof} Without loss of generality, we may assume that $r_-=1$ and $1<r_+<\infty$. Suppose that $\xi>6$ is sufficiently large, ultimately depending only on $r_+$, to be specified below. Fix sequences $\{J_i\}_{i=1}^I$ and $(k_i)_{i=1}^I$ as above. We claim that for every integer $1\leq n\leq I$, there exists a unique $1\leq m\leq n$ such that $k_m=\min_{i=1}^n k_i$ and \begin{equation}\label{f-induction} \sum_{i=1}^n \diam f(J_i)\leq (1+2\xi^{-1}+4\xi^{-2}+8\xi^{-3}+\cdots)\xi^{-k_m}r_+.\end{equation} The base case $n=1$ is trivial. Suppose for induction that the claim holds whenever $1\leq n\leq N$. Put $n=N+1$ and choose an index $1\leq m\leq N+1$ such that $k_m=\min_{i=1}^{N+1} k_i$. (We do not yet know that $m$ is unique.) Separate $\{J_i:1\leq i\leq N+1\}\setminus\{J_m\}$ into two chains $\{J_1,\dots,J_{m-1}\}$ and $\{J_{m+1},\dots,J_{N+1}\}$, each of which contain no more than $N$ intervals. Working with the first chain (unless it is empty), the induction hypothesis and fact $\xi>6$ implies that there exists a unique $1\leq j\leq m-1$ such that $k_j=\min_{i=1}^{m-1} k_i$ and $$\sum_{i=1}^{m-1}\diam f(J_i)\leq (1+2\xi^{-1}+4\xi^{-2}+8\xi^{-3}+\cdots)\xi^{-k_j}r_+<(3/2)\xi^{-k_j}r_+.$$ We know that $k_j\geq k_m$ by our selection of $m$. Suppose for contradiction that $k_j=k_m$. By property (3), we can locate $P=\lceil 3r_+\rceil$ pairwise disjoint closed intervals $K_1,K_2,\dots, K_P$ such that each $K_p$ lies between $J_j$ and $J_m$ and $\diam f(K_p)\geq \xi^{-k_j}$. This immediately implies that $J_j$ and $J_m$ do not intersect, so $j\leq m-2$. For any $j+1\leq i\leq m-1$, we have $k_i\geq k_j+1$ by our specification of $j$; hence $\diam f(J_i) \leq \xi^{-k_i}r_+\leq \xi^{-1}\xi^{-k_j}r_+ <\xi^{-k_j}$ provided that $\xi>r_+$. This ensures that each interval $J_i$ with $j+1\leq i\leq m-1$ intersects at most two of the intervals $K_p$. (Otherwise, if some $J_i$ intersected three or more intervals, then $J_i\supset K_p$ for some $p$ and $\xi^{-k_j}>\diam f(J_i)\geq \diam f(K_p)\geq \xi^{-k_j}$. We remark that if we allowed $\diam f(J_i)\geq \diam f(K_p)$ for some $p$, then we could not control the overlap; see Example \ref{bad-gap}.) Now, $\bigcup_{p=1}^P K_p \subset \bigcup_{i=j+1}^{m-1} J_i$, because the latter intervals connect $J_j$ and $J_m$ by (1). All together, these observations yield $$P\xi^{-k_j} \leq \sum_{p=1}^P \diam f(K_p) \leq 2\sum_{i=j+1}^{m-1} \diam f(J_i) <3\xi^{-k_j}r_+.$$ Thus, $P<3r_+ \leq \lceil 3r_+\rceil=P$, which is a contradiction. Therefore, we have  $k_j\geq k_m+1$. It follows that $k_i>k_m$ for all $1\leq i\leq m-1$ and $$\sum_{i=1}^{m-1} \diam f(J_i)\leq \xi^{-1}(1+2\xi^{-1}+4\xi^{-2}+8\xi^{-3}+\cdots)\xi^{-k_m}r_+.$$ A similar conclusion holds for $m+1\leq i\leq N+1$. Adding the estimates  verifies the induction step.

To conclude, let $M\geq 1$ be an index such that $k_M=\min_{i=1}^I k_i$, which exists because the sequence $(k_i)_{i=1}^I$ is bounded from below. Because for all integers $1\leq n\leq I$, there is a unique index $1\leq m\leq n$ such that $k_m=\min_{i=1}^n k_i$, we discover that $M$ is unique and $k_m=k_M$ for all $n\geq M$. Finally, \eqref{e:interval-growth} follows from \eqref{f-induction}, since $\xi>6$.\end{proof}

\begin{proof}[{Proof of Lemma \ref{build-filtrations}}]
In the construction below, we view $\Diam(\tau)$, the diameter of the image of an arc $\tau$, as a weight assigned to $\Domain(\tau)$, a closed interval lying in the domain of the map $f$. From this vantage point, a filtration is simply a nested family of finite partitions of some set $K\subset[0,1]$ into intervals with geometrically decaying weights. This heuristic can be formalized using the pseudometric from Lemma \ref{curve-pseudometric}. As shorthand, we may speak of the intersection or union of arcs, but always mean the arc formed by taking the intersection or union in the domains of the arcs.

Let $\rho>1$, let $0<\underline{A}<A<\infty$, and let $J\geq 1$ be an integer, with $J$ sufficiently large depending on $\rho$ and $A/\underline{A}$ to be specified below. Without loss of generality, we may assume that $\underline{A}=1$ and $A>1$. Let $\mathscr{F}^0$ be an admissible family of arcs in $\overline{\Sigma}$. Our plan is to first preprocess $\mathscr{F}^0$, partitioning it into a finite number of families $\mathscr{D}^1,\mathscr{D}^2,\dots$ of well separated arcs. We then transform each family $\mathscr{D}^j$ into a nested family $\mathscr{E}^j$ of arcs satisfying \eqref{new-arc-diameters} and \eqref{new-arc-properties}. (This is where we use Lemma \ref{interval-lemma}.) Afterwards, we describe how to extend each family $\mathscr{E}^j$ to a filtration $\mathscr{F}^j$ with the same conclusions.

To begin, using the bounded overlap assumption, break $\mathscr{F}^0$ into $C$ or fewer nonempty families such that within each family, arcs in the same level $n$ are pairwise disjoint. Note that some levels of a family may be empty. By splitting each family into $$\lceil 3A\rceil+1\quad(<5A)$$ families, as necessary, we may assume that any pair distinct arcs in level $n$ of a family are separated in the domain by at least $\lceil 3A\rceil$ disjoint arcs $\sigma$ with $\Diam(\sigma)\geq \rho^{-n}$. Next, break apart each family of arcs into $J$ (or fewer) nonempty families by jumping $J$ generations at a time, i.e.~ within each family, group together all levels $n$ with $n-n_0\equiv m \pmod J$, relabeling so that each original level $n=(n_0+m)+kJ$ is assigned to level $N=(n_0+m)+k$. Denote the resulting $5ACJ$ or fewer families by $\mathscr{D}^1=\bigcup_{n=n_1}^\infty \mathscr{D}^1_{n}$, $\mathscr{D}^2=\bigcup_{n=n_2}^\infty \mathscr{D}^2_{n}$, \dots, with starting index $n_j\in\{n_0,n_0+1,\dots,n_0+J-1\}$ for all $j$. Note that level $n$ in $\mathscr{D}^j$ corresponds to level $n_j+J(n-n_j)=Jn-(J-1)n_j$ in $\mathscr{F}^0$. By our assumption on the diameters of arcs in $\mathscr{F}^0$, if $n\geq n_j$ and $\tau\in\mathscr{D}^j_n$, then \begin{equation}\label{D-diam}\rho^{(J-1)n_j}\rho^{-Jn}\leq \Diam(\tau)\leq (A\rho^{(J-1)n_j})\rho^{-Jn}.\end{equation} Further, any pair of distinct arcs in $\mathscr{D}^j_n$ are separated in the domain by at least $\lceil 3A\rceil$ disjoint intervals $\sigma$ with $\Diam(\sigma)\geq \rho^{(J-1)n_j}\rho^{-Jn}$.

Next, we transform $\mathscr{D}^1,\mathscr{D}^2,\dots$ into nested families $\mathscr{E}^1,\mathscr{E}^2,\dots$ of arcs. This will require us to join certain overlapping arcs, thereby increasing the diameter of arcs. By choosing $J$ to be sufficiently large, we can control the growth of diameters.

Fix a family $\mathscr{D}^j$. For each $n\geq n_j$, let $\mathscr{D}^j_{n+}=\mathscr{D}^j_{n+1}\cup \mathscr{D}^j_{n+2}\cup\mathscr{D}^j_{n+3}\cup\cdots$. Fix a level $n\geq n_j$ and an arc $\tau\in \mathscr{D}^j_n$. We inductively define a sequence of arcs by fusing the arc and all overlapping arcs from future generations relative to $\tau$. Formally, we set $$\tau_0=\tau,\quad \tau_{k+1}=f\,\Bigg{|}\bigcup_{\{\sigma\in\mathscr{D}^j_{n+}:\sigma\cap \tau_k\neq\emptyset\}}\Domain(\tau)\quad(k\geq 0),\quad \tau_\infty = f\,\Bigg{|}\,\overline{\bigcup_{k=0}^\infty \Domain(\tau_k)},$$ where $\sigma\cap\tau_k\neq\emptyset$ means $\Domain(\sigma)\cap\Domain(\tau_k)\neq\emptyset$ and $f|I$ denotes the restriction of $f$ to $I$. Let us upper bound $\Diam(\tau_\infty)$ in terms of $\Diam(\tau)$. Each arc $\tau_{k+1}$ can be obtained by taking the union of $\tau_k$ and up to two arcs $\sigma^-_k$ and $\sigma^+_k$ that intersect, but are not contained in $\tau_k$, one of each side. Thus, $\tau_\infty$ can be obtained by concatenating $\tau$ and two arcs, $\tau^-$ and $\tau^+$, with domains on either side of $\Domain(\tau)$, where $\tau^-$ and $\tau^+$ can each be expressed as the union of a finite or infinite chain of arcs in $\mathscr{D}^j_{n+}$. By Lemma \ref{curve-pseudometric}, $\Diam(\tau^\pm)\leq \sum_k\Diam(\sigma^\pm_k)$. To continue, we need $J$ to be sufficiently large. Assume that $\rho^J>6A\geq \max\{6,A\}$. By Lemma \ref{interval-lemma}, with $\xi=\rho^J$, $r_-=\rho^{(J-1)n_j}$, and $r_+=A\rho^{(J-1)n_j}$, \begin{equation}\label{e:tau-extensions} \Diam(\tau^\pm) \leq (1+3/\rho^J)(A\rho^{(J-1)n_j})\rho^{-J(n+1)}< \frac14\rho^{(J-1)n_j}\rho^{-Jn}\leq \frac14\Diam(\tau).\end{equation} Therefore, \begin{equation}\label{E-upper} \Diam(\tau_\infty) \leq \Diam(\tau^-)+\Diam(\tau)+\Diam(\tau^+)
<\frac{3}{2}\Diam(\tau)\leq \frac{3}{2}(A\rho^{(J-1)n_j})\rho^{-Jn}.\end{equation} Recall that distinct arcs $\tau$ and $\tilde\tau$ in $\mathscr{D}^j_n$ are separated by an arc (in fact, several arcs) $\sigma$ of $\overline{\Sigma}$ with $\Diam(\sigma)\geq \rho^{(J-1)n_j}\rho^{-Jn}$. Thus, by \eqref{e:tau-extensions} and the triangle inequality, for any distinct $\tau,\tilde\tau\in\mathscr{D}^j_n$, the domains of $\tau_\infty$ and $\tilde\tau_\infty$ are separated by some arc $\sigma$ in $\overline{\Sigma}$ with \begin{equation}\label{E-separate} \Diam(\sigma)>\frac{1}{2}\rho^{(J-1)n_j}\rho^{-Jn}.\end{equation}

For each $j$ and $n\geq n_j$, define $\mathscr{E}^j_n$ by including the arc $\tau_\infty$ in $\mathscr{E}^j_n$ for each arc $\tau\in\mathscr{D}^j_n$. By \eqref{D-diam} and \eqref{E-upper}, the families $\mathscr{E}^j=\bigcup_{n=n_j}^\infty\mathscr{E}^j_n$ satisfy---even better inequalities than---\eqref{new-arc-diameters} and \eqref{new-arc-properties}, where the roles of ``$\tau'\in\mathscr{F}^0_n$'' and ``$\tau\in\mathscr{F}^j_N$'' in the statement of the lemma are played by $\tau\in\mathscr{D}^j_n$ and $\tau_\infty\in\mathscr{E}^j_n$, respectively.  Let us verify that $\mathscr{E}^j$ is a nested family of arcs in the sense that if (the domains of) $\tau_\infty$ and $\tilde\tau_\infty$ in $\mathscr{E}^j$ intersect in more than one point, then $\Domain(\tau_\infty)$ is contained in $\Domain(\tilde\tau_\infty)$ or vice versa. Suppose that $t\geq n\geq n_j$, $\tau\in\mathscr{D}^j_n$, $\tilde\tau\in\mathscr{D}^j_t$, and $\tau_\infty$ and $\tilde\tau_\infty$ intersect in at least two points. Then the interiors of $\Domain(\tau_\infty)$ and $\Domain(\tilde\tau_\infty)$ intersect. Hence $\Domain(\tau_k)$ and $\Domain(\tilde\tau_l)$ intersect for some $k,l<\infty$, as well. Because $\tilde\tau_l$ is a finite union of overlapping arcs $\mathscr{D}^j_{(t-1)+}$ and $t\geq n$, it follows that $\tilde\tau_l\subset \tau_{k+m}$ for some sufficiently large positive integer $m$. Thus, $\tilde\tau_{l+i}\subset\tau_{k+m+i}$ for all $i\geq 0$ by the construction, and therefore, $\tilde\tau_\infty\subset\tau_\infty$. As an immediate corollary, distinct arcs in each level $\mathscr{E}^j_n$ intersect in at most one point. However, by \eqref{E-separate}, we know even more: for every $j$ and $n\geq n_j$, the arcs in $\mathscr{E}^j_n$ are pairwise disjoint. Reviewing the construction so far, we see that the map that assigns each arc $\tau\in\mathscr{F}^0_n$ to an arc $\tau_\infty\in\mathscr{E}^j_{n_j+(n-n_j)/J}$ is injective.

To build filtrations $\mathscr{F}^j$ out of the families $\mathscr{E}^j$, it remains to add arcs to each level of a family, as necessary, to ensure that the domain of any arc is the union of the domains of its children. Fix a family $\mathscr{E}^j$. We will work top-down, starting with level $n=n_j$. Let $\tau_1,\dots,\tau_k$ denote the arcs in $\mathscr{E}^j_{n}$, ordered from left to right using the natural order of their  domains as subsets of $\RR$. It is possible that this list is empty. From left to right, expand $$f=\sigma_0\cup\tau_1\cup\sigma_1\cup\tau_2\cup\cdots \cup \tau_{k-1}\cup\sigma_{k-1}\cup \tau_k\cup\sigma_k$$ with domains of consecutive arcs intersecting in exactly one point; the initial and final arc $\sigma_0$ and $\sigma_k$ may be empty or nonempty. If $\mathscr{E}^j_{n}$ is empty, write $f=\sigma_0$. Suppose that $\sigma=\sigma_i$ is nonempty. There are three alternatives.

\emph{Alternative 1.} If it happens that $\Diam(\sigma)$ is between $(1/4)\rho^{(J-1)n_j}\rho^{-n}$ and $\rho^{(J-1)n_j}\rho^{-n}$, we make no modification and include $\sigma$ in $\mathscr{F}^j_{n}$.

\emph{Alternative 2.} Another possibility is that $\Diam(\sigma)<(1/4)\rho^{(J-1)n_j}\rho^{-n}$. If $\sigma=\sigma_i$ for some $0\leq i\leq k-1$, set $\tau=\tau_{i+1}$, the first arc the right; if $\sigma=\sigma_k$, set $\tau=\tau_k$, the first arc to the left. Replace $\tau$ by $\sigma\cup\tau$. By \eqref{E-upper}, the modified arc satisfies $$\Diam(\tau)<(3/2)A\rho^{(J-1)n_j}\rho^{-n}+(1/4)\rho^{(J-1)n_j}\rho^{-n}<(7/4)A\rho^{(J-1)n_j}\rho^{-n}.$$ (Exceptionally, if $f=\sigma_0\cup\tau_1\cup\sigma_1$ and both $\sigma_0$ and $\sigma_1$ have small diameters, replace $\tau$ by $\sigma_0\cup\tau\cup\sigma_1$ instead; in this case $$\Diam(\tau)<(3/2)A\rho^{(J-1)n_j}\rho^{-n}+2(1/4)\rho^{(J-1)n_j}\rho^{-n}<2A\rho^{(J-1)n_j}\rho^{-n}.)$$ Similarly, if the original arc $\tau\in\mathscr{E}^j_{n}$ came from extending an arc $\tau'\in\mathscr{D}^j_{n}$, then the modified arc $\tau$ satisfies $\Diam(\tau)<(7/4)\Diam(\tau')$ (exceptionally, $\Diam(\tau)<2\Diam(\tau')$). Include the modified arc $\tau$ in $\mathscr{F}^j_n$.

\emph{Alternative 3.} Lastly, suppose that $\Diam(\sigma)>\rho^{(J-1)n_j}\rho^{-n}$. We will partition $\sigma$ into a finite number of shorter arcs $\zeta_i$ of diameter between $(1/4)\rho^{(J-1)n_j}\rho^{-n}$ and $\rho^{(J-1)n_j}\rho^{-n}$ (see Remark \ref{r:Diam-is-continuous}), but need to do this in an intelligent way in order to maintain nested levels. Let $\{\xi_1,\xi_2,\dots\}$ be an enumeration of the maximal arcs in $\mathscr{E}^j_{n+}=\mathscr{E}^j_{n+1}\cup\mathscr{E}^j_{n+2}\cup\cdots$ whose domains are contained in $\Domain(\sigma)$. (There may be none.) Earlier we declared that $\rho^J>6A$. Hence $\Diam(\xi_l) <(3/2)A\rho^{(J-1)n}\rho^{-(n+1)}<(1/4)\rho^{(J-1)n_j}\rho^{-n}$. We need to make sure that the endpoints of arcs in our partitions of $\sigma$ do not lie in the interior of $\Domain(\xi_l)$ for any $l$. Proceed as follows. Let $a$ be the left endpoint of $\Domain(\sigma)$. Choose $t>0$ as large as possible so that $\Diam(f|_{[a,a+t]})=(1/4)\rho^{(J-1)n_j}\rho^{-n}$, which we may do because the diameter varies continuously in $t$ and $\Diam(\sigma)$ is large. If $a+t$ does not lie in the interior of $\Domain(\xi_l)$ for any $l$, then we set $\zeta_1=f|_{[a,a+t]}$ and have $\Diam(\zeta_1)=(1/4)\rho^{(J-1)n_j}\rho^{-n}$. Otherwise, if $a+t$ lies in the interior $\xi_l$ for some $l$, set $\zeta_1=f|_{[a,a+t]}\cup \xi_l$, which satisfies $(1/4)\rho^{(J-1)n_j}\rho^{-n}\leq \Diam(\zeta_1)<(1/2)\rho^{(J-1)n_j}\rho^{-n}.$ Repeat a similar construction on $\sigma\setminus \zeta_1$, this time letting $a$ be the left endpoint of $\Domain(\sigma\setminus\zeta_1)$ and choosing $t>0$ as large as possible so that $\Diam(f|_{[a,a+t]})\leq (1/4)\rho^{(J-1)n_j}\rho^{-n}$. If $\Diam(f|_{[a,a+t]})=(1/4)\rho^{(J-1)n_j}$, then continue as before and iterate. Otherwise, at some stage, $\sigma\setminus(\zeta_1\cup\dots\cup\zeta_m)\neq\emptyset$, but $\Diam(f|_{[a,a+t]})<(1/4)\rho^{(J-1)n_j}\rho^{-n}$, where $t$ is the right endpoint of the domain of $\sigma$. Replace $\zeta_m$ by $\zeta_m\cup f|_{[a,a+t]}$. This may increase the diameter of $\zeta_m$, but in any event, the modified arc satisfies $(1/4)\rho^{(J-1)n_j}\rho^{-n}\leq \Diam(\zeta_m)<(3/4)\rho^{(J-1)n_j}\rho^{-n}$. Include each of the intervals $\zeta_1$, \dots, $\zeta_m$ from the partition of $\sigma$ in $\mathscr{F}^j_{n}$.

Carry out the indicated construction for each nonempty $\sigma=\sigma_i$. Also include any arc $\tau_i$ in $\mathscr{F}^j_n$ if it was not already included. This completes the definition of $\mathscr{F}^j_{n}$. To define the next level $\mathscr{F}^j_{n+1}$, repeat the construction from above on each arc $\tau$ in $\mathscr{F}^j_n$ independently, with $f$ replaced by $\tau$. By induction, we obtain a definition of $\mathscr{F}^j=\bigcup_{n=n_j}^\infty\mathscr{F}^j_n$ on each level. The family $\mathscr{F}^j$ is the desired filtration.
\end{proof}

\bibliography{btsp-refs}
\bibliographystyle{amsbeta}

\end{document}